\documentclass[12pt, reqno]{amsart}
\usepackage{amssymb}
\usepackage{mathrsfs}
\usepackage{mathpazo}
\usepackage[colorlinks=true, breaklinks=true]{hyperref}
\usepackage{bm}
\usepackage[latin1]{inputenc}
\usepackage{graphicx}
\newtheorem{theorem}{Theorem}

\newtheorem{corollary}[theorem]{Corollary}

\newtheorem{definition}[theorem]{Definition}

\newtheorem{lemma}[theorem]{Lemma}

\newtheorem{proposition}[theorem]{Proposition}
\newtheorem{remark}[theorem]{Remark}

\newcommand{\bea}{\begin{eqnarray}}
\newcommand{\eq}{\end{eqnarray}}
\newcommand{\eea}{\end{eqnarray}}
\newcommand{\bqn}{\begin{eqnarray*}}
\newcommand{\beaa}{\begin{eqnarray*}}
\newcommand{\eqn}{\end{eqnarray*}}
\newcommand{\eeaa}{\end{eqnarray*}}
\newcommand{\bpr}{\begin{proposition}}
\newcommand{\epr}{\end{proposition}}

\newcommand{\cal}{\mathcal}

\numberwithin{equation}{section}
\numberwithin{theorem}{section}

\begin{document}
\title{Multivariate Stochastic Volatility Models and Large Deviation Principles}
\author{Archil Gulisashvili}\let\thefootnote\relax\footnotetext{Department of Mathematics, Ohio University, Athens OH 45701; e-mail: gulisash@ohio.edu}

\date{}

\begin{abstract}
We establish a comprehensive sample path large deviation principle (LDP) for log-price processes associated with multivariate
time-inhomogeneous stochastic volatility models. Examples of models for which the new LDP holds include Gaussian models, non-Gaussian fractional models, mixed models, models with reflection, and models in which the volatility process is a solution to a Volterra type stochastic integral equation. The sample path and small-noise LDPs for log-price processes are used to obtain large deviation style asymptotic formulas for the distribution function of the first exit time of a log-price process from an open set, multidimensional binary barrier options, 
call options, Asian options, and the implied volatility. Such formulas capture leading order asymptotics of the above-mentioned important quantities arising in the theory of stochastic volatility models. We also prove a sample path LDP for solutions to Volterra type stochastic integral equations with predictable coefficients depending on auxiliary stochastic processes.

\end{abstract}

\maketitle

\vspace{0.2in}
\vspace{0.3in}
\begin{flushright}
Dedicated to the memory of Peter Carr
\end{flushright}
\vspace{0.1in}
\section{Introduction} A classical Black-Scholes-Merton model of option pricing assumes that the volatility of a financial asset is constant.
Stochastic volatility models provided corrections by taking into account random features of the volatility. In modern stochastic volatility models, the volatility is modeled by a stochastic process.

One of the main objectives in the present paper is to introduce and study 
general multivariate time-inhomogeneous stochastic volatility models. Such a model is described by the following multidimensional stochastic differential equation:
\begin{align}
dS_t&=S_t\circ[b(t,\widehat{B}_t)dt+\sigma(t,\widehat{B}_t)(\bar{C}dW_t+CdB_t)],\quad 0\le t\le T,\quad S_0=s_0\in\mathbb{R}^m
\label{E:moodk}
\end{align}
where the initial condition $s_0=(s_0^{(1)},\cdots,s_0^{(m)})$ is such that $s_i> 0$ for all $1\le i\le m$. The meanings of the symbols appearing in the previous equation will be explained in Section \ref{S:mfm}. The model in (\ref{E:moodk}) incorporates various features of numerous known stochastic volatility models (see the discussion in Section \ref{S:mfm} and the survey in Section \ref{S:unix}). The interested reader can find in \cite{G} detailed information about classical stochastic volatility models (the Hull-White, the Stein-Stein, and the Heston model). Multivariate models are discussed in \cite{AAY}, see also Chapter 11 in the monograph \cite{Ber} by L. Bergomi. This chapter is titled Multi-Asset Stohastic Volatililty. Interesting examples of time-inhomogeneous models are the rough Bergomi model 
introduced in \cite {BFG} and the super rough Bergomi model (see \cite{Gul1} and \cite{BHP}).

Stochastic volatility models are widely used in finance (see, e.g., \cite{Ber,FPS,FPSS,G,H-L,Kahl,Lew1,Lew2}). 
The stochastic model in (\ref{E:moodk}) is characterized by the drift map $b$, the volatility map $\sigma$, and the volatility process 
$\widehat{B}$. Under the conditions formulated in the next section, the equation in (\ref{E:moodk}) is uniquely solvable, and the solution 
$S_t=(S_t^{(1)},\cdots,S_t^{(m)})$, $t\in[0,T]$, is a continuous stochastic process with strictly positive components. 
The functions $S^{(i)}$, with $1\le i\le m$, can be interpreted as price processes of correlated risky assets in a portfolio or an index. 
The log-price process associated with the stochastic volatility model introduced in (\ref{E:moodk}) is defined by 
$X=(X^{(1)},\cdots,X^{(m)})$ where $X^{(i)}=\log S^{(i)}$ for
$1\le i\le m$. The initial condition for the log-price process is denoted by $x_0$. It is clear that $x_0=(\log s_0^{(1)},\cdots,s_0^{(m)})$.

A major part of this paper is devoted to asymptotic analysis of stochastic volatility models. We use sample path and small-noise large deviation principles to perform such analysis. A sample path large deviation principle (LDP) for a stochastic process characterizes logarithmic asymptotics of the probability that the path of a scaled version of the process belongs to a given set of paths. 
The theory of sample path large deviations goes back to the celebrated work 
of Varadhan \cite{V1} and Freidlin and Wentzell \cite{FW}. For more information about large deviations, see \cite{DZ,DS,DE,FK,St,V2,V3}. 
Our main goal in the present paper is to obtain a universal sample path LDP for log-price processes in multivariate stochastic volatility models that unifies the results established in the above-mentioned publications and also provides new results. Such a comprehensive LDP is formulated in Section \ref{S:lpo} (see Theorem \ref{T:nmn}).

Various sample path LDPs are known for log-price processes in stochastic volatility models (see, e.g.,\,\cite{CaP,CP,FZ,FGP,GGG,Gul1,Gul2,Gul3,GVol,JPS}). I would like to bring the attention of the interested reader to the book \cite{FGGJT} titled Large Deviations and Asymptotic Methods in Finance and also to the paper \cite{Ph} of H. Pham devoted to applications of the large deviation theory in mathematical finance.
It is also worth mentioning that there is a rich literature devoted to applications of large deviation principles to the study of the the asymptotic behavior of various quantities arising in finance (see, e.g., \cite{BC,CMD,FZ,FGGJT,GGG,Gul1,Gul2,Gul3,GVol,HJL,JPS,Ph,R}).

The unification of various large deviation principles is achieved 
in Theorem \ref{T:nmn} due to the wide variety of admissible Volterra type volatility processes used in this paper. Examples of stochastic models, for which the LDP obtained in Theorem \ref{T:nmn} holds, include multivariate Gaussian models, multivariate non-Gaussian fractional models, mixed models, multivariate models with reflection, and models in which the volatility process is a solution to a certain Volterra type stochastic integral equation. The restrictions imposed on the drift map and the volatility map in Theorem \ref{T:nmn} are rather mild 
(see Assumption A in Section \ref{S:lpo}). We also obtain large deviation style asymptotic formulas for the distribution function of the first exit time of the log-price process from an open set and similar asymptotic formulas for multidimensional binary barrier options. 

Most of the volatility processes used in the present paper are solutions to certain Volterra type integral equations (see 
the equation in (\ref{E:2a})).
The coefficients in this equation are predictable maps depending on auxiliary stochastic processes. We impose special restrictions on the equation in (\ref{E:2a}) (see Assumptions (C1) - (C7) introduced and discussed in Section \ref{S:VP}). These restrictions are based on Conditions (H1) - (H6) used in the paper \cite{CF} of Chiarini and Fischer. However, we had to adapt Conditions (H1) - (H6) to our setting since the stochastic processes employed in \cite{CF} are not of Volterra type. Moreover, we impose an additional restriction 
(Assumption (C6)) that is not needed in the case of non-Volterra stochastic differential equations studied in \cite{CF}. Under Assumptions  
(C1) - (C7), we prove a sample path LDP for the unique solution to the Volterra type stochastic integral equation in (\ref{E:2a}) (see Theorem \ref{T:beg1}). The LDP in Theorem \ref{T:beg1} generalizes various known LDPs for Volterra type processes. Theorems 
\ref{T:nmn} and \ref{T:beg1} are the main results obtained in the present paper.

The paper of Chiarini and Fischer was an important source of ideas in our work on LDPs for Volterra type volatility processes. The methods employed in 
\cite{CF} are based on variational representations of functionals of Brownian motion and the weak convergence approach to small-noise LDP problems. We use the same techniques in the proofs of our results concerning LDPs for volatility processes (see Theorems \ref{T:beg1} 
and \ref{T:beg11}). The weak convergence method was developed by Dupuis and Ellis (see \cite{DE}). Various sample path LDPs for families of functionals of Brownian motion were obtained in \cite{BD,BDa,BDb,BDM1,BDM2}. 

We will next give an overview of the contents of the present paper. In Section \ref{S:mfm}, we introduce multivariate stochastic volatility
models and their scaled versions. Section \ref{S:VP} deals with volatility processes, scaled volatility processes, controlled counterparts of Volterra type stochastic integral equations, and special restrictions that will be imposed 
on the volatility models. In Section 
\ref{S:lpo}, a sample path LDP for log-price processes is formulated and explained (see Theorem \ref{T:nmn}), while
Section \ref{S:snl} analyzes small-noise LDPs for log-price processes. 
In Section \ref{S:im}, we establish a sample path LDP for solutions to Volterra type stochastic integral equations (see Theorem \ref{T:beg1}). This theorem uses the canonical set-up. We do not know whether Theorem \ref{T:beg1} holds on any set-up (see Definition \ref{D:gens}). 
It is worth mentioning that for certain less general Volterra type stochastic integral equations, the 
LDP in Theorem \ref{T:beg1} is valid on any set-up. Exceptional examples here are volatility processes in multivariate Gaussian stochastic volatility models (see Theorem \ref{T:dert}), volatility processes in multivariate non-Gaussian fractional models 
(see Theorem \ref{T:riur}), and also Volterra type stochastic processes used in the paper \cite{NR} of Nualart and Rovira 
(see also the paper \cite{RS} of Rovira and Sanz-Sol\'{e} devoted to large deviations for stochastic Volterra equations in the plane). 

Section \ref{S:unix} of the present paper is devoted to examples of stochastic volatility models for which the sample path LDP in Theorem 
\ref{T:nmn} holds. It also provides examples of models for which Assumptions (C1) -- (C7) formulated in Section \ref{S:VP} are satisfied. Section \ref{S:unix} is split into several subsections. In Subsections \ref{SS:GM} and \ref{SS:frH}, we give a brief overview of one-factor Gaussian models and one-factor non-Gaussian fractional models, respectively. For more information about such models see \cite{GGG,Gul1,Gul3}. In Subsection \ref{SS:fu}, we merge multivariate Gaussian models and multivariate non-Gaussian fractional models, and show that Theorem \ref{T:nmn} holds for such mixtures.  Subsection \ref{SS:mGm} is devoted to LDPs for log-price processes in multivariate Gaussian stochastic volatility models on a general set-up under mild restrictions on the drift map, the volatility map, and the volatility process. Similar results are obtained in Subsection \ref{SS:mGm} for log-price processes in multivariate non-Gaussian fractional models. Note that Theorem \ref{T:dert} is a generalization of Theorem 4.2 in \cite{Gul1}. The former theorem provides an LDP for the 
log-price process in a multivariate Gaussian stochastic volatility model, while the latter one deals with the one-dimensional case. In Subsection \ref{SS:mvr}, we prove an LDP for the log-process in a multivariate stochastic volatility model with reflection. The final two subsections of Section \ref{S:unix} (Subsections \ref{SS:Wang} and 
\ref{SS:NR}) deal with volatility processes which are solutions to Volterra type stochastic integral equations. In Subsection 
\ref{SS:Wang}, we discuss the LDP obtained in \cite{Z} by Zhang, while Subsection \ref{SS:NR} is devoted to the LDP established in \cite{NR} by Nualart and Rovira.

We would also like to highlight the paper \cite{JP} of Jacquier and Pannier in which the authors prove sample path LDPs for solutions to less general Volterra type stochastic integral equations than the equation in 
(\ref{E:2a}) (see Theorems 3.8 and 3.25 in \cite{JP}). Another paper that is worth mentioning is the paper \cite{CP} of Catalini and Pacchiarotti in which a sample path LDP for multivariate time-homogeneous Gaussian models is established under stronger restrictions than those employed in our Theorem \ref{T:dert}. Interesting results were obtained in the paper \cite{BFGMS} of Bayer, Friz, Gassiat, Martin, 
and Stemper. It was established in the latter paper that a small-noise LDP holds for scaled log-processes 
in certain one-factor Volterra type stochastic volatility models (see (55), (56), and Corollary 5.5 in \cite{BFGMS}). 
The authors of \cite{BFGMS} used Hairer's regularity structures in their work. Corollary 5.5 in \cite{BFGMS}, under the restriction that the  stochastic volatility model is defined on the canonical set-up, follows from the results obtained in the present paper
(see Theorems \ref{T:rrr}, \ref{T:bed}, and Remark \ref{R:begin}). More details will be provided in Remark \ref{R:bf} in Subsection \ref{SS:Wang}.

Section \ref{S:appl} is devoted to applications of the methods developed in the present paper to mathematical finance.
Here we perform a small-noise asymptotic analysis of various important objects of study in the theory of stochastic volatility models, for instance, 
 the distribution function of the first exit time of the log-price process from an open set in $\mathbb{R}^m$ (Subsection \ref{SS:dcf}),
 barrier options (Subsection \ref{SS:bbo}), call options (Subsection \ref{SS:call}), the implied volatility (Subsection \ref{SS:sniv}), and Asian options (Subsection \ref{SS:AA}).
We obtain large deviation style asymptotic formulas for the objects mentioned in the previous sentence. These formulas capture the leading order asymptotic behavior of stochastic volatility models of our interest as the small noise parameter $\varepsilon$ tends to zero. They also provide asymptotic approximations to the above-mentioned objects. 

Assumption $B$ (see (\ref{E:um})) plays an important role in the present paper.
In Subsection \ref{SS:re}, examples of models for which Assumption $B$ holds are provided. They include multivariate Gaussian models, multivariate generalized fractional Heston models, and mixed models. 

In Subsection \ref{SS:TM} of Section \ref{S:appl}, we study a simple model (a toy model) using the methods developed in the present paper. A special uncorrelated SABR model plays the role of the toy model. We obtain various estimates for the rate function 
in the toy model (see Lemma \ref{L:use}, Corollary \ref{C:abo}, and Theorem \ref{T:kont}), and also for the small-noise limit of the implied volatility in the toy model 
(see Theorem \ref{T:on}). 

In Section \ref{S:proofs}, we prove Theorem \ref{T:beg1}, while in Section \ref{S:prov}, we include the proof of Theorem \ref{T:nmn}.

\section{Multivariate Stochastic Volatility Models}\label{S:mfm}
Let $\mathbb{R}^m$ be $m$-dimensional Euclidean space
equipped with the norm $||\cdot||_m$. For a real $(m\times m)$-matrix $M$, its Frobenius norm will be denoted by $||M||_{m\times m}$ and the symbol $M^{\prime}$ will stand for the transpose of $M$. We will next provide more details about the equation in (\ref{E:moodk}). 
This equation is defined on a probability space $(\Omega,{\cal F},\mathbb{P})$ 
carrying two independent $m$-dimensional standard Brownian motions $W$ and $B$ with respect to the measure $\mathbb{P}$, 
and the symbol $\circ$ in (\ref{E:moodk}) stands for the Hadamard (component-wise) product of vectors. 
The components of the initial condition $s_0$ of the process $S$ are strictly positive. The matrix $C=(c_{ij})$ in 
(\ref{E:moodk}) is a real $(m\times m)$-matrix such that $||C||_{m\times m}< 1$. It is clear that the matrix $\rm Id_m-C^{\prime}C$ is symmetric and positive definite, and we will denote the unique symmetric and positive definite square root of the matrix 
$\rm Id_m-C^{\prime}C$ by $\bar{C}$. The elements of the matrix $\bar{C}$ will be denoted by $\bar{c}_{ij}$. Under the previous conditions, the matrix $\bar{C}$ is invertible. By $\{\mathcal{F}_t\}_{0\le t\le T}$ is denoted the augmentation of the filtration generated by the processes $W$ and $B$ (see, e.g., Definition 7.2 in Chapter 2 of \cite{KaS}). We will also use the augmentation of the filtration generated by the process $B$, and denote it by $\{\mathcal{F}^B_t\}_{0\le t\le T}$. The symbol $b$ in (\ref{E:moodk}) stands for a continuous map defined on $[0,T]\times\mathbb{R}^d$ with values in $\mathbb{R}^m$. We call $b$ the drift map. By $\sigma$ is denoted a continuous map of $[0,T]\times\mathbb{R}^d$ into the space of $(m\times m)$ real matrices. This map will be called the volatility map. The process 
$\widehat{B}=(\widehat{B}^{(1)},\cdots,
\widehat{B}^{(d)})$ appearing in (\ref{E:moodk}) is a continuous $d$-dimensional stochastic process defined in terms of Brownian motion $B$ and adapted to the filtration $\{\mathcal{F}^B_t\}_{0\le t\le T}$. The process
$\widehat{B}$ will be called the volatility process (see Definition \ref{D:vo}). We have already mentioned in the introduction that the model in (\ref{E:moodk}) can be interpreted as a time-inhomogeneous 
stochastic volatility model describing the time-behavior of price processes of correlated risky assets. The matrix-valued process 
$\sigma(t,\widehat{B}_t)$, with $t\in[0,T]$, characterizes the joint volatility of these assets. 

The next definition introduces general set-ups (we adopt the terminology used in \cite{RW}).
\begin{definition}\label{D:gens}
The system $(\Omega,W,B,{\cal F}_T,\{\mathcal{F}_t\}_{0\le t\le T},\mathbb{P})$ is called a set-up associated with the model
in (\ref{E:moodk}), while the system $(\Omega,B,{\cal F}_T^B,\{\mathcal{F}^B_t\}_{0\le t\le T},\mathbb{P})$ is 
called a set-up associated with the volatility process in (\ref{E:moodk}). 
\end{definition}

In terms of the components $S_t^{(i)}$ of the process $S$, with $1\le i\le m$, the equation in (\ref{E:moodk}) can be rewritten as the following system of stochastic differential equations:
\begin{align}
&dS_t^{(i)}=S_t^{(i)}[b_i(t,\widehat{B}_t)dt+\sum_{k,j=1}^m\bar{c}_{jk}\sigma_{ij}(t,\widehat{B}_t)dW_t^{(k)} 
+\sum_{k,j=1}^mc_{jk}\sigma_{ij}(t,\widehat{B}_t)dB_t^{(k)}],\,\,1\le i\le m.
\label{E:erw}
\end{align}
Recall that we denoted by $X$ the log-price process associated with the model in (\ref{E:moodk}). 
We will next characterize the dynamics of the process $X$. For every $i$ with $1\le i\le m$, the equation in (\ref{E:erw}) is a linear stochastic differential equation driven by the process
$$
G_t^{(i)}=\int_0^tb_i(s,\widehat{B}_s)ds+\sum_{k,j=1}^m\bar{c}_{j,k}\int_0^t\sigma_{ij}(s,\widehat{B}_s)dW_s^{(k)} 
+\sum_{k,j=1}^mc_{j,k}\int_0^t\sigma_{ij}(s,\widehat{B}_s)dB_s^{(k)}.
$$
This process is a continuous semimatringale since
$$
\int_0^T||b(s,\widehat{B}_s)||_mdt+\int_0^T||\sigma(s,\widehat{B}_s)||_{m\times m}^2ds<\infty
$$
$\mathbb{P}$-a.s. Next, using the Dol\'{e}ans-Dade formula, we see that  
\begin{align}
S_t^{(i)}&=s_0^{(i)}\exp\{\int_0^tb_i(s,\widehat{B}_s)ds
-\frac{1}{2}\int_0^t\sum_{j=1}^m\sigma_{ij}(s,\widehat{B}_s)^2ds
+\sum_{k,j=1}^m\bar{c}_{j,k}\int_0^t\sigma_{ij}(s,\widehat{B}_s)dW_s^{(k)} \nonumber \\
&\quad+\sum_{k,j=1}^mc_{j,k}\int_0^t\sigma_{ij}(s,\widehat{B}_s)dB_s^{(k)}\}.
\label{E:fort}
\end{align}
The formula in (\ref{E:fort}) can be rewritten as follows:
\begin{align*}
&S_t^{(i)}=s_0^{(i)}\exp\{\int_0^tb_i(s,\widehat{B}_s)ds
-\frac{1}{2}\int_0^t\sum_{j=1}^m\sigma_{ij}(s,\widehat{B}_s)^2ds
+[\int_0^t\sigma(s,\widehat{B}_s)(\bar{C}dW_s+CdB_s)]_{i}\}.
\end{align*}
Recall that for any $m\times m$-matrix $D$, we denoted by $D^{\prime}$ the transpose of $D$.  We will also denote by $\mbox{diag}(D)$ the vector whose components are the diagonal elements of $D$. It is clear that the component $a_i$ of the vector 
$\mbox{diag}(\sigma(s,\widehat{B}_s)\sigma(s,\widehat{B}_s)^{\prime})$, with $1\le i\le m$, is given by 
$a_i=\sum_{j=1}^m\sigma_{ij}(s,\widehat{B}_s)^2$.
It follows that the log-price process $X$ is given by
\begin{align}
X_t&=x_0+\int_0^tb(s,\widehat{B}_s)ds-\frac{1}{2}\int_0^t\mbox{diag}(\sigma(s,\widehat{B}_s)\sigma(s,\widehat{B}_s)^{\prime})ds
\nonumber \\
&\quad+\int_0^t\sigma(s,\widehat{B}_s)(\bar{C}dW_s+CdB_s),\quad 0\le t\le T.
\label{E:mood}
\end{align}
\begin{remark}\label{R:pp}
The process obtained from the log-price process $X$ by removing one of the drift terms, more precisely, the term
$-\frac{1}{2}\int_0^t\mbox{\rm diag}(\sigma(s,\widehat{B}_s)\sigma(s,\widehat{B}_s)^{\prime})ds$, will be denoted by $\widehat{X}$. 
We have
\begin{equation}
\widehat{X}_t=x_0+\int_0^tb(s,\widehat{B}_s)ds+\int_0^t\sigma(s,\widehat{B}_s)(\bar{C}dW_s+CdB_s).
\label{E:moodd}
\end{equation}
We call the process in (\ref{E:moodd}) the modified log-price process.
\end{remark}
\begin{remark}\label{R:imm}
In the case where $m=1$, we use the correlation parameter $\rho\in(-1,1)$ and set $\bar{\rho}=\sqrt{1-\rho^2}$.
Then, the equation describing the evolution of the process $S$ is as follows:
$$
dS_t=S_t[b(t,\widehat{B}_t)dt+\sigma(t,\widehat{B}_t)(\bar{\rho}dW_t+\rho dB_t)],\quad S_0=s_0> 0.
$$
Moreover, the log-price process is given by
$$
X_t=x_0+\int_0^tb(s,\widehat{B}_s)ds-\frac{1}{2}\int_0^t\sigma(s,\widehat{B}_s)^2ds+\int_0^t\sigma(s,\widehat{B}_s)
(\bar{\rho}dW_s+\rho dB_s)
$$
where $x_0=\log s_0$. 
\end{remark}

A modulus of continuity is a nonnegative nondecreasing funciton $\omega$ on $[0,\infty)$ such that $\omega(s)\rightarrow 0$ 
as $s\rightarrow 0$.
Let $x=(t_1,v_1)$ and $y=(t_2,v_2)$ be elements of the space $[0,T]\times\mathbb{R}^d$ equipped with the Euclidean distance
$
\nu_d(x,y)=\sqrt{(t_1-t_2)^2+||v_1-v_2||_d^2}.
$
Denote by $\overline{B_d(r)}$ the closed ball centered at $(0,0)$ of radius $r> 0$ in the metric space defined above,  
and let $\omega$ be a modulus of continuity on $[0,\infty)$. 
\begin{definition}\label{D:dd}
A map $\lambda:[0,T]\times\mathbb{R}^d\mapsto\mathbb{R}^1$ is called locally $\omega$-continuous
if for every $r> 0$ there exists $L(r)> 0$ such that for all $x,y\in\overline{B_d(r)}$
the following inequality holds:
$$
|\lambda(x)-\lambda(y)|\le L(r)\omega(\nu_d(x,y)).
$$ 
\end{definition}

We will next explain what restrictions on the drift map $b$ and the volatility map $\sigma$ are imposed in the present paper. These restrictions are rather mild. \\ 
\\
\it Assumption A: \rm The components of the drift map $b$ and the elements of the volatility map $\sigma$ are locally $\omega$-continuous on the space $[0,T]\times\mathbb{R}^d$ for some modulus of continuity $\omega$. 
In addition, the elements of the volatility map $\sigma$ are not identically zero on $[0,T]\times\mathbb{R}^d$. 
\\
\\
Let $\varepsilon\in(0,1]$ be the scaling parameter. The scaled version of the log-price process $X$ is defined by
\begin{align}
X_t^{(\varepsilon)}&=x_0+\int_0^tb(s,\widehat{B}_s^{(\varepsilon)})ds
-\frac{1}{2}\varepsilon\int_0^t\mbox{diag}(\sigma(s,\widehat{B}^{(\varepsilon)}_s)\sigma(s,\widehat{B}^{(\varepsilon)}_s)^{\prime})ds
\nonumber \\
&\quad+\sqrt{\varepsilon}\int_0^t\sigma(s,\widehat{B}_s^{(\varepsilon)})
(\bar{C}dW_s+CdB_s)
\label{E:iu}
\end{align}
where $X_0^{(\varepsilon)}=x_0$ for all $s\in(0,1]$. The scaled volatility process $\widehat{B}^{(\varepsilon)}$ appearing in (\ref{E:iu}) is introduced in the next section (see Definition \ref{D:dds}). We will also use the process
\begin{equation}
\widehat{X}_t^{(\varepsilon)}=x_0+\int_0^tb(s,\widehat{B}_s^{(\varepsilon)})ds
+\sqrt{\varepsilon}\int_0^t\sigma(s,\widehat{B}_s^{(\varepsilon)})
(\bar{C}dW_s+CdB_s)
\label{E:gyt}
\end{equation} 
that is a scaled version of the modified log-price process defined in (\ref{E:moodd}).
 
\section{Volatility Processes}\label{S:VP}
Our main aim in this section is to introduce the volatility process $\widehat{B}$ that is used in (\ref{E:moodk}). We will need several definitions. For a positive integer $p\ge 1$, the symbol ${\cal W}^p$ will stand for the space of continuous $\mathbb{R}^p$-valued maps 
on $[0,T]$ equipped with the following norm:
$||f||=\max_{t\in[0,T]}||f(t)||_p$, $f\in{\cal W}^p$.
Let $B_s$, with $s\in[0,T]$, be the coordinate process on ${\cal W}^p$. Define a filtration on the space ${\cal W}^p$ by 
${\cal B}_t^p=\sigma(B_s:0\le s\le t)$, $t\in[0,T]$. The augmentation $\{\widetilde{{\cal B}}_t^p\}$ of the filtration 
$\{{\cal B}_t^p\}$ is called the canonical filtration on ${\cal W}^p$. Let $\mathbb{P}$ be the Wiener measure on $\widetilde{{\cal B}}_T^p$. 
\begin{definition}\label{D:canny}
The ordered system $({\cal W}^p,B,\widetilde{{\cal B}}_T^p,\{\widetilde{{\cal B}}_t^p\},\mathbb{P})$ is called the canonical set-up
on ${\cal W}^p$. 
\end{definition}

The set-up introduced in Definition \ref{D:canny} is a special case of a general set-up associated with the volatility process
in (\ref{E:moodk}) (see Definition \ref{D:gens}). 
The coordinate process $s\mapsto B_s$ plays the role of $p$-dimensional standard Brownian motion with respect to the measure $\mathbb{P}$. 
\begin{remark}\label{R:oneof}
One of the reasons why the canonical set-up is employed in the present paper is the following known fact.
Let $Z$ be an $\{\widetilde{{\cal B}}^m_t\}$-adapted continuous stochastic process on ${\cal W}^m$ with state space
$\mathbb{R}^d$. Then, there exists a process $\widetilde{Z}$ adapted to the filtration
$\{{\cal B}^m_t\}$ and indistinguishable from $Z$. Moreover, there is a
functional $j:{\cal W}^m\mapsto{\cal W}^d$ such that $X=j(B)$\,\,$\mathbb{P}$-a.s. on ${\cal W}^m$, 
and for every $t\in[0,T]$, the functional $j$ is ${\widetilde{\cal B}}_t^m/{\cal B}_t^d$-measurable. The functional $j$ is generated by the process $X$. 
In addition, the canonical probability space plays 
an important role in the proof of the equality in (\ref{E:opposum}).
\end{remark}

We will next define the canonical set-up on the space $\Omega=\Omega_1\times\Omega_2={\cal W}^m\times{\cal W}^m$. 
Denote the coordinate processes on $\Omega_1$ and $\Omega_2$ 
by $W$ and $B$, respectively, and consider the filtration on $\Omega$ generated by the process $t\mapsto(W_t,B_t)$, $t\in[0,T]$. Denote by 
$\{{\cal F}_t\}$ the augmentation of this filtration with respect to the measure $\mathbb{P}=\mathbb{P}_1\times\mathbb{P}_2$, where
$\mathbb{P}_1$ and $\mathbb{P}_2$ are the Wiener measures on $\Omega_1$ and $\Omega_2$, respectively. By $\{{\cal F}^B_t\}$ will be
denoted the augmentation of the filtration generated by process $t\mapsto B_t$, $t\in[0,T]$. The processes $W$ and $B$ are independent 
$m$-dimensional Brownian motions defined on the space $\Omega$. The canonical set-up on the space $\Omega={\cal W}^m\times{\cal W}^m$ is the system $(\Omega,W,B,{\cal F}_T,\{{\cal F}_t\},\{{\cal F}^B_t\},\mathbb{P})$. 

Let $Y$ be a stochastic process satisfying the following Volterra type stochastic integral equation on ${\cal W}^m$ equipped with the canonical set-up:
\begin{equation}
Y_t=y+\int_0^ta(t,s,V^{(1)},Y)ds+\int_0^tc(t,s,V^{(2)},Y)dB_s.
\label{E:2a}
\end{equation} 
In (\ref{E:2a}), $a$ is a map from the space $[0,T]^2\times{\cal W}^{k_1}\times{\cal W}^d$ into the space $\mathbb{R}^d$, while $c$ is a map 
from the space $[0,T]^2\times{\cal W}^{k_2}\times{\cal W}^d$ into the space of $(d\times m)$-matrices. For a matrix $M$ belonging to the latter space, the symbol $||M||_{d\times m}$ will stand for its Frobenius norm. 
Assumption (C1) formulated below introduces more restrictions on the maps $a$ and $c$. The vector 
$y\in\mathbb{R}^d$ in (\ref{E:2a}) is the fixed initial condition for the process $Y$. 
The processes $V^{(i)}$, with $i=1,2$, appearing in (\ref{E:2a}) are fixed auxiliary continuous stochastic processes on ${\cal W}^m$ with 
state spaces $\mathbb{R}^{k_1}$ and $\mathbb{R}^{k_2}$, respectively. They satisfy the following stochastic differential equations:
\begin{equation}
V_s^{(i)}=V_0^{(i)}+\int_0^s\bar{b}_i(r,V^{(i)})dr+\int_0^s\bar{\sigma}_i(r,V^{(i)})dB_r,\quad i=1,2.
\label{E:2b}
\end{equation}
In (\ref{E:2b}), $V_0^{(i)}\in\mathbb{R}^{k_i}$ are initial conditions, $\bar{b}_i$ are maps of $[0,T]\times{\cal W}^{k_i}$ 
into $\mathbb{R}^{k_i}$, while $\bar{\sigma}_i$ are maps of $[0,T]\times{\cal W}^{k_i}$ into the space of $k_i\times m$-matrices.
We assume that the equations in (\ref{E:2b}) satisfy Conditions (H1) - (H6) introduced
in \cite{CF} by Chiarini and Fischer. 
\begin{remark}\label{R:imp}
Examples of equations in (\ref{E:2b}), for which Conditions (H1) - (H6) hold true, are provided in 
Sections 3 and 4 of \cite{CF}. For instance, the validity of these conditions is established in \cite{CF} for 
equations with locally Lipschitz coefficients satisfying the sub-linear growth condition (see Definitions A1 and 
A2 in Section 3 of \cite{CF},
see also (11.1) on p. 128 in \cite{RW}). It is also shown in \cite{CF} that Conditions (H1) - (H6) hold true for one-dimensional 
diffusion equations with H\"{o}lder dispersion coefficient.
\end{remark} 

The next definition introduces volatility processes that are used throughout the present paper.
\begin{definition}\label{D:vo}
The volatility process $\widehat{B}$ appearing in (\ref{E:moodk}) is as follows:
$\widehat{B}=GY$, where $Y$ satisfies the equation in (\ref{E:2a}), while $G$ is a continuous map from ${\cal W}^d$ into itself that is 
$\widetilde{{\cal B}}_t^d/{\cal B}_t^d$-measurable for every $t\in[0,T]$. 
\end{definition}

An example of a map $G$ satisfying the condition in Definition \ref{D:vo} for $d=1$
is the Skorokhod map (see Subsection \ref{SS:mvr}).
The measurability condition for $G$ is included in Definition \ref{D:vo} in order the volatility process $\widehat{B}$ 
to be adapted to the filtration $\{\widetilde{{\cal B}}^d_t\}$.

Let $\varepsilon\in(0,1]$ be a small-noise parameter. A scaled version of the equation in (\ref{E:2a}) has the following form:
\begin{equation}
Y_t^{(\varepsilon)}=y+\int_0^ta(t,s,V^{1,\varepsilon},Y^{(\varepsilon)})ds
+\sqrt{\varepsilon}\int_0^tc(t,s,V^{2,\varepsilon},Y^{(\varepsilon)})dB_s.
\label{E:2ooo}
\end{equation}
In (\ref{E:2ooo}), the process $V^{i,\varepsilon}$, with $i=1,2$, is a scaled version of the process $V^{(i)}$. It satisfies the equation
\begin{equation}
V_s^{i,\varepsilon}=V_0^{(i)}+\int_0^s\bar{b}_i(r,V^{i,\varepsilon})dr+\sqrt{\varepsilon}\int_0^s\bar{\sigma}_i(r,V^{i,\varepsilon})dB_r.
\label{E:4}
\end{equation}
\begin{remark}\label{R:rert}
In \cite{CF}, more general equations than the equation appearing in (\ref{E:4}) are considered.  The coefficient maps $\bar{b}$ and $\bar{\sigma}$ in those equations may depend on the scaling parameter $\varepsilon$.
We do not study such equations in the present paper.
\end{remark}
\begin{definition}\label{D:dds}
The scaled volatility process $\widehat{B}^{(\varepsilon)}$ is given by $\widehat{B}^{(\varepsilon)}=GY^{(\varepsilon)}$ 
where $G$ is introduced in Definition \ref{D:vo}, while $Y^{(\varepsilon)}$ is the solution to (\ref{E:2ooo}).
\end{definition}

Controlled counterparts of the equations in (\ref{E:2a}) - (\ref{E:4}) will also be considered.
Let ${\cal M}^2[0,T]$ be the space of all $\mathbb{R}^m$-valued square-integrable $\{{\cal F}^B_t\}$-predictable processes. The controls will be chosen from the space ${\cal M}^2[0,T]$. Deterministic controls will be employed as well. 
They are functions belonging to the space $L^2([0,T],\mathbb{R}^m)$.
\begin{definition}\label{D:contr}
Let $N> 0$. By ${\cal M}^2_N[0,T]$ is denoted the class of controls $v\in{\cal M}^2[0,T]$ such that
\begin{equation}
\int_0^T||v_s||_m^2ds\le N\quad\mathbb{P}-\mbox{a.s.}
\label{E:cve}
\end{equation}
\end{definition}

Suppose $v\in{\cal M}^2[0,T]$. Then, controlled counterparts of the equations 
in (\ref{E:2a}) and (\ref{E:2b}) are as follows:
\begin{align}
Y_t^{(v)}&=y+\int_0^ta(t,s,V^{1,v},Y^{(v)})ds+\int_0^tc(t,s,V^{2,v},Y^{(v)})v_sds \nonumber \\
&\quad+\int_0^tc(t,s,V^{2,v},Y^{(v)})dB_s
\label{E:3ooo}
\end{align}
and
$$
V_s^{i,v}=V_0^{(i)}+\int_0^s\bar{b}_i(r,V^{i,v})dr+\int_0^s\bar{\sigma}_i(r,V^{i,v})v_rdr
+\int_0^s\bar{\sigma}_i(r,V^{i,v})dB_r,\quad i=1,2. 
$$
For $v\in{\cal M}^2[0,T]$ and $\varepsilon\in(0,1]$, controlled counterparts of the equations in (\ref{E:2ooo}) and (\ref{E:4}) satisfy
\begin{align}
Y_t^{\varepsilon,v}&=y+\int_0^ta(t,s,V^{1,\varepsilon,v},Y^{\varepsilon,v})ds+\int_0^tc(t,s,V^{2,\varepsilon,v},Y^{\varepsilon,v})v_sds
\nonumber \\
&\quad+\sqrt{\varepsilon}\int_0^tc(t,s,V^{2,\varepsilon,v},Y^{\varepsilon,v})dB_s
\label{E:asaf}
\end{align}
and
\begin{align}
V_s^{i,\varepsilon,v}&=V_0^{(i)}+\int_0^s\bar{b}_i(r,V^{i,\varepsilon,v})dr+\int_0^s\bar{\sigma}_i(r,V^{i,\varepsilon,v})v_rdr
\nonumber \\
&\quad+\sqrt{\varepsilon}\int_0^s\bar{\sigma}_i(r,V^{i,\varepsilon,v})dB_r,\quad i=1,2.
\label{E:asa}
\end{align}
\begin{remark}\label{R:cons}
In \cite{CF}, less general equations than those in (\ref{E:2ooo}) and (\ref{E:asaf}) were studied. 
These equations are as follows:
\begin{align}
&Y_t^{(\varepsilon)}=y+\int_0^ta(s,Y^{(\varepsilon)})ds+\sqrt{\varepsilon}\int_0^tc(s,Y^{(\varepsilon)})dB_s.
\label{E:ali}
\end{align}
and
\begin{align}
&Y_t^{\varepsilon,v}=y+\int_0^ta(s,Y^{\varepsilon,v})ds+\int_0^tc(s,Y^{\varepsilon,v})v_sds
+\sqrt{\varepsilon}\int_0^tc(s,Y^{\varepsilon,v})dB_s.
\label{E:5t}
\end{align}
The equations in (\ref{E:ali}) and (\ref{E:5t}) are not of Volterra type.
\end{remark}

It will be explained next what happens if $\varepsilon=0$. The equations in (\ref{E:4}) take the following form:
$V_s^{i,0}=V_0^{(i)}+\int_0^s\bar{b}_i(r,V^{i,0})dr$, $i=1,2$.
These equations can be solved pathwise, and for every $i$ all the solutions are the same by the uniqueness. 
Let us denote the solution by $v^{i,0}$. It follows that $V_s^{i,0}=v^{i,0}(s)$ for $i=1,2$ and all $s\in[0,T]$.

Suppose $\varepsilon=0$. Then, the equations in (\ref{E:asa}) and (\ref{E:asaf}) can be rewritten as follows:
$$
V_s^{i,0,v}=V_0^{(i)}+\int_0^s\bar{b}_i(r,V^{i,0,v})dr+\int_0^s\bar{\sigma}_i(r,V^{i,0,v})v_rdr,\quad i=1,2,
$$
and
\begin{equation}
Y_t^{0,v}=y+\int_0^ta(t,s,V^{1,0,v},Y^{0,v})ds+\int_0^tc(t,s,V^{2,0,v},Y^{0,v})v_sds,
\label{E:pppp}
\end{equation}
respectively. In a special case of a deterministic control $f$, we have
\begin{equation}
V_s^{i,0,f}=V_0^{(i)}+\int_0^s\bar{b}_i(r,V_r^{i,0,f})dr+\int_0^s\bar{\sigma}_i(r,V_r^{i,0,f})f(r)dr.
\label{E:rty}
\end{equation}
\begin{remark}\label{R:redi}
Under the restrictions imposed on $\bar{b}_i$ and $\bar{\sigma}_i$ in \cite{CF}, the functional equations
$$
\psi_i(s)=V_0^{(i)}+\int_0^s\bar{b}_i(r,\psi_i)dr+\int_0^s\bar{\sigma}_i(r,\psi_i)f(r)dr,\quad i=1,2,
$$
are uniquely solvable, the solutions $\psi_{i,f}$ belong to the spaces ${\cal W}^{k_i}$, and if $f_n\mapsto f$ weakly in 
$L^2([0,T],\mathbb{R}^m)$, then $\psi_{i,f_n}\mapsto\psi_{i,f}$ in ${\cal W}^{k_i}$ for $i=1,2$. Therefore, the solution 
to the equation in (\ref{E:rty}) is deterministic, and $V_s^{i,0,f}=\psi_{i,f}(s)$ for all $s\in[0,T]$ (see \cite{CF}). 
\end{remark}

Suppose $v\in{\cal M}^2_N[0,T]$ for some $N> 0$.
Then, it follows from Girsanov's Theorem that for all $0<\varepsilon\le 1$, the process 
\begin{equation}
B^{\varepsilon, v}_t=B_t+\frac{1}{\sqrt{\varepsilon}}\int_0^tv_sds,\quad t\in[0,T],
\label{E:brow}
\end{equation}
is an $m$-dimensional Brownian motion on ${\cal W}^m$ with respect to a measure $\mathbb{P}^{\varepsilon,v}$ on ${\cal F}_T^m$ that is equivalent to the 
measure $\mathbb{P}$. The process $B^{\varepsilon, v}$ is adapted to the filtration 
$\{{\cal F}_t^B\}$. 
In a special case where $\varepsilon=1$, the following notation will be used:
\begin{equation}
B^{(v)}_t=B_t+\int_0^tv_sds,\quad t\in[0,T].
\label{E:brew}
\end{equation}

We will next explain what restrictions are imposed on the model for the volatility described by the equation in (\ref{E:2a}). 
\\
\\
\it \underline{Assumption (C1)}. \rm For all $(\eta_1,\varphi)\in{\cal W}^{k_1}\times{\cal W}^d$, the map 
$(t,s)\mapsto a(t,s,\eta_1,\varphi)$ is Borel measurable, with values in the space $\mathbb{R}^d$, while for all 
$(\eta_2,\varphi)\in{\cal W}^{k_2}\times{\cal W}^d$, the map $(t,s)\mapsto c(t,s,\eta_2,\varphi)$ is Borel measurable, with values in the space of $d\times m$-matrices. Moreover, 
$a$ and $c$ are of Volterra type in the first two variables. In addition, for every $t\in[0,T]$, 
$(s,\eta_1,\varphi)\mapsto a(t,s,\eta_1,\varphi)$ and 
$(s,\eta_2,\varphi)\mapsto c(t,s,\eta_2,\varphi)$ are predictable path functionals mapping the space 
$[0,t]\times{\cal W}^{k_1}\times{\cal W}^d$ into the space 
$\mathbb{R}^d$ and the space
$[0,t]\times{\cal W}^{k_2}\times{\cal W}^d$ into the space of $d\times m$ matrices, respectively. The definition of a predictable path functional can be found in \cite{RW} (see Definition (8.3) and Remark (8.4) on p. 122). The requirement above is similar to 
Convention (8.7) on p. 123 in \cite{RW}.
\\
\\
\it \underline{Assumption (C2)}. \rm (a)\,\,Let $\eta_1\in{\cal W}^{k_1}$, $\eta_2\in{\cal W}^{k_2}$, and $\varphi\in{\cal W}^d$. Then, the following inequalities hold for all $t\in[0,T]$:
\begin{equation}
\int_0^t||a(t,s,\eta_1,\varphi)||_dds<\infty\quad\mbox{and}\quad\int_0^T||c(t,s,\eta_2,\varphi)||_{d\times m}^2ds<\infty.
\label{E:c2}
\end{equation}
(b)\,\,For all fixed $\eta_1\in{\cal W}^{k_1}$ and $\varphi\in{\cal W}^d$, the function $t\mapsto\int_0^ta(t,s,\eta_1,\varphi)ds$ 
is a continuous $\mathbb{R}^d$-valued function on $[0,T]$. In addition, for every fixed $t\in[0,T]$
the function $(\eta_1,\varphi)\mapsto\int_0^ta(t,s,\eta_1,\varphi)ds$ is continuous on the space ${\cal W}^{k_1}\times{\cal W}^d$. \\
\\
(c)\,\,Let $\eta_{2,n}\rightarrow\eta_2$ in ${\cal W}^{k_2}$ and $\varphi_n\rightarrow\varphi$ in ${\cal W}^d$ as $n\rightarrow\infty$. Then, for every $t\in[0,T]$,
$$
\int_0^t||c(t,s,\eta_{2,n},\varphi_n)-c(t,s,\eta_2,\varphi)||^2_{d\times m}ds\rightarrow 0\quad\mbox{as}\quad n\rightarrow\infty. 
$$
\vspace{0.1in}
\begin{flushleft}
\it \underline{Assumption (C3)}. \rm (a)\,\,For all $0<\varepsilon\le 1$ there exists a strong solution to 
the equation in (\ref{E:2ooo}). (b)\,\,Let $v\in M^2_N[0,T]$ for some $N> 0$. Then, any two strong solutions to the 
equation in (\ref{E:3ooo}) are $\mathbb{P}$-indistinguishable.
\end{flushleft}
\begin{remark}\label{R:Cc}
The definition of the strong solution used in the present paper includes the continuity of the solution. 
\end{remark}
\begin{remark}\label{R:rer}
Assumption (C3)(b) is weaker than the pathwise uniqueness condition employed in \cite{CF}.
\end{remark}
\it \underline{Assumption (C4)}. \rm For every function $f\in L^2([0,T],\mathbb{R}^m)$ and the functions $\psi_{1,f}$ and 
$\psi_{2,f}$ defined in Remark \ref{R:redi}, the equation
\begin{equation}
\eta(t)=y+\int_0^ta(t,s,\psi_{1,f},\eta)ds+\int_0^tc(t,s,\psi_{2,f},\eta)f(s)ds,
\label{E:t2}
\end{equation}
is uniquely solvable in ${\cal W}^d$. 
\begin{remark}\label{R:hatt}
It will be shown below that the equation in (\ref{E:t2}) is always solvable. The details can be found in Remark \ref{R:tyu}.
Therefore, only the uniqueness condition must be included in Assumption (C4).
\end{remark}
\begin{definition}\label{D:C4}
The map $\Gamma_y:L^2([0,T],\mathbb{R}^m)\mapsto{\cal W}^d$ is defined by 
$\Gamma_y f=\eta_f$ where $\eta_f$ is the unique solution to the equation in (\ref{E:t2}).
\end{definition}
\begin{flushleft}
\it \underline{Assumption (C5)}. \rm Set $D_N=\{f\in L^2([0,T],\mathbb{R}^m):\int_0^T||f(t)||^2_mdt\le N\}$. 
Then, the restriction of the map $\Gamma_y$ to $D_N$ is a continuous map from $D_N$ equipped with the weak topology 
into the space ${\cal W}^d$. 
\end{flushleft}
\vspace{0.1in}

Let $v\in M^2_N[0,T]$ for some $N> 0$. Then, there exists a map $g^{(2)}:{\cal W}^m\mapsto{\cal W}^{k_2}$ satisfying the following conditions:
(i)\,\,$g^{(2)}(B)=V^{(2)}$;\,\,(ii)\,\,$g^{(2)}(B^{(v)})=V^{2,v}$\,\,$\mathbb{P}$-a.s.;\,\,(iii)\,\,For every $t\in[0,T]$, $g^{(2)}$ is
$\widetilde{{\cal B}}_t^m/{\cal B}_t^{k_2}$-measurable (see Lemma A.1 in \cite{CF}, see also Theorem 10.4 on p. 126 in \cite{RW}). 
\begin{remark}\label{R:ghg}
Suppose Assumption (C3) holds. Then, for every $\varepsilon\in(0,1]$ there exists a map 
$h^{\varepsilon}:{\cal W}^m\mapsto{\cal W}^d$ such that the solution $Y^{(\varepsilon)}$ to the equation in (\ref{E:3ooo}) satisfies 
$Y^{(\varepsilon)}=h^{(\varepsilon)}(B)$, and the map $h^{\varepsilon}$ is $\widetilde{{\cal B}}_t^m/{\cal B}_t^d$-measurable 
for all $t\in[0,T]$ (see Remark \ref{R:oneof}). 
It is clear that $Y^{(1)}=Y$. We will denote the map $h^{(1)}$ by $h$. It follows that
$Y=h(B)$\,\,$\mathbb{P}$-a.s.
\end{remark}
\begin{flushleft}
\it \underline{Assumption (C6)}. \rm The process
$t\mapsto\int_0^tc(t,s,g^{(2)}(B^{(v)}),h(B^{(v)}))dB^{(v)}_s$, where $t\in[0,T]$,
is continuous.  
\end{flushleft}
\vspace{0.1in}
Assumption (C6) looks rather complicated. A special case, where Assumption (C6) is satisfied, is when the map $c$ does not depend on the 
variable $t$. Indeed, in such a case, the correctness of Assumption (C6) follows from the condition for the map $c$ in Assumption (C2)(a) and the continuity properties of stochastic integrals. More examples of the validity of Assumption (C6) will be provided in Section \ref{S:unix}.
\\
\\
\it \underline{Assumption (C7)}. \rm Suppose $0<\varepsilon_n< 1$, with $n\ge 1$, is a sequence of numbers such that 
$\varepsilon_n\rightarrow 0$ as 
$n\rightarrow\infty$. Let $v^{(n)}$, $n\ge 1$, be a sequence of controls satisfying the condition $v^{(n)}\in{\cal M}^2_N[0,T]$ 
for some $N> 0$ and all
$n\ge 1$ (see Definition \ref{D:contr}). Then,
the family of ${\cal W}^d$-valued random variables $Y^{\varepsilon_n,v^{(n)}}$, with $n\ge 1$, is tight in 
${\cal W}^d$. Moreover, for every $t\in[0,T]$, the following inequality is satisfied:
\begin{equation}
\sup_{n\ge 1}\int_0^t\mathbb{E}\left[||c(t,s,V^{2,\varepsilon_n,v^{(n)}},Y^{\varepsilon_n,v^{(n)}})||_{d\times m}^2\right]ds<\infty.
\label{E:cond}
\end{equation}

We will next formulate several remarks related to Assumptions (C1) - (C7).
\begin{remark}\label{R:Aa}
Assumptions (C1) and (C2)(a) guarantee that the equation in (\ref{E:2a}) exists.
\end{remark}
\begin{remark}\label{R:Bb}
It is easy to see that the first part of Assumption (C2)(b) holds if the following condition is satisfied:
Let $t,t^{\prime}\in[0,T]$, and suppose $\eta_1\in{\cal W}^{k_1}$, $\varphi\in{\cal W}^d$, and $t^{\prime}\rightarrow t$. Then,
\begin{equation}
\int_0^T||a(t^{\prime},s,\eta_1,\varphi)-a(t,s,\eta_1,\varphi)||_dds\rightarrow 0.
\label{E:ll}
\end{equation}
Similarly, the second part of Assumption (C2)(b) can be derived from the following condition:
$$
\int_0^T||a(t,s,\eta_{1,n},\varphi_n)-a(t,s,\eta_1,\varphi)||_dds\rightarrow 0
$$
as $n\rightarrow 0$
provided that $\eta_{1,n}\rightarrow\eta$ in ${\cal W}^{k_1}$ and $\varphi_n\rightarrow\varphi$ in ${\cal W}^d$ as $n\rightarrow\infty$.
\end{remark}
\begin{remark}\label{R:ut}
It follows from Assumption (C2)(a) that for every $t\in[0,T]$, the following condition holds 
$\mathbb{P}$-a.s.:
$
\int_0^t||a(t,s,V^{1,\varepsilon},Y^{(\varepsilon)})||_dds+\int_0^t||c(t,s,V^{2,\varepsilon},Y^{(\varepsilon)})||_{d\times m}^2ds<\infty.
$
\end{remark}

It was shown in \cite{CF} that if pathwise uniqueness and existence in the strong sense hold for 
the equation in (\ref{E:ali}) and the control $v$ satisfies (\ref{E:cve}), then there exists a unique strong solution to the equation in 
(\ref{E:5t}), and this solution can be represented as a measurable functional of the process $B^{\varepsilon,v}$
(see Lemma 1 in Appendix A in \cite{CF}). We will establish a similar strong solvability result for the Volterra equation in (\ref{E:asaf}) 
(see Lemma \ref{L:uni} and Remark \ref{R:rr}). 
\begin{remark}\label{R:rrr}
Using Remark \ref{R:redi} and (\ref{E:pppp}), we see that 
\begin{equation}
Y_t^{0,f}=y+\int_0^ta(t,s,\psi_{1,f},Y^{0,f})ds+\int_0^tc(t,s,\psi_{2,f},Y^{0,f})f(s)ds.
\label{E:rtys}
\end{equation}
It follows from Assumption (C4) that the equation in (\ref{E:rtys}) can be reduced to the equation in (\ref{E:t2}). Therefore, 
$Y_t^{0,f}=\Gamma_y f(t)$ for all $t\in[0,T]$. 
\end{remark}
\section{Sample Path Large Deviation Principles for Log-Price Processes}\label{S:lpo}
In this section, we formulate and discuss one of the main results of the present paper (Theorem \ref{T:nmn}). 
The proof of this theorem will be given in Section \ref{S:prov}. 

Let us denote by $\mathbb{C}_0^m$ the subspace of the space ${\cal W}^m$ consisting of the functions $f$ such that $f(0)=\vec{0}$. 
Throughout the paper, the symbol $(\mathbb{H}^1_0)^m$ will stand for the $m$-dimensional Cameron-Martin space. 
A function $f$ from the space $\mathbb{C}_0^m$ belongs to the space 
$(\mathbb{H}^1_0)^m$ if it has absolutely continuous components and the derivatives of the components on $(0,T)$ 
are square-integrable with respect to the Lebesgue measure. 
For $f\in(\mathbb{H}^1_0)^m$, we 
set $\dot{f}=(\dot{f}_1,\cdots,\dot{f}_m)$ where the symbol $\dot{f}_k$, with $1\le k\le m$, stands for the 
derivative of the component $f_k$ with respect to the variable $t$.

In the next definition, we introduce special maps $f\mapsto {\cal A}f$ and $f\mapsto\widehat{f}$. 
\begin{definition}\label{D:lip}
(i)\,\,For every function $f\in L^2([0,T];\mathbb{R}^m)$, the function ${\cal A}f\in{\cal W}^d$ is defined by 
${\cal A}{f}=G(\Gamma_yf)$
where $G$ and $\Gamma_y$ are introduced in Definitions \ref{D:vo} and \ref{D:C4}, respectively.
(ii)\,\,For every function $f\in(\mathbb{H}_0^1)^m$, the function $\widehat{f}\in{\cal W}^d$ is defined by
$\widehat{f}={\cal A}\dot{f}$.
\end{definition}

The rate function $\widetilde{Q}_T$ governing the large deviation principle for the log-price process depends on the 
measurable map $\Phi:\mathbb{C}_0^{m}\times\mathbb{C}_0^m\times{\cal W}^d\mapsto\mathbb{C}_0^m$ given by
\begin{equation} 
\Phi(l,f,h)(t)=\int_0^tb(s,\widehat{f}(s))ds+\int_0^t\sigma(s,\widehat{f}(s))\bar{C}\dot{l}(s)ds
+\int_0^t\sigma(s,\widehat{f}(s))C\dot{f}(s)ds
\label{E:disc}
\end{equation}
for all $l,f\in(\mathbb{H}^1_0)^m$, $h=\widehat{f}\in{\cal W}^d$, and $0\le t\le T$. For all the remaining triples 
$(l,f,h)$, we set $\Phi(l,f,h)(t)=0$  for $t\in[0,T]$. Let $g\in\mathbb{C}_0^m$, and define the function $\widetilde{Q}_T$ by
\begin{align}
&\widetilde{Q}_T(g) \nonumber \\
&=\inf_{l,f\in(\mathbb{H}_0^1)^m}\left[\frac{1}{2}\int_0^T||\dot{l}(s)||_m^2ds
+\frac{1}{2}\int_0^T||\dot{f}(s)||_m^2ds:\Phi(l,f,\widehat{f}(t))=g(t),\,
t\in[0,T]\right],
\label{E:QT1}
\end{align}
if the equation $\Phi(l,f,\widehat{f}(t))=g(t)$ is solvable for $l$ and $f$. If there is no solution, then we set 
$\widetilde{Q}_T(g)=\infty$. It follows that if the previous equation is solvable, then $g\in(\mathbb{H}_0^1)^m$ and
\begin{equation}
\dot{g}(t)=b(t,\widehat{f}(t))+\sigma(t,\widehat{f}(t))\bar{C}\dot{l}(t)+\sigma(t,\widehat{f}(t))C\dot{f}(t).
\label{E:rtw}
\end{equation}

The next assertion provides a sample path LDP for log-price processes in general stochastic volatility models. 
\begin{theorem}\label{T:nmn}
Suppose Assumption A and Assumptions (C1) -- (C7) hold true, and the model in (\ref{E:moodk}) is defined on the canonical set-up.
Then, the process $\varepsilon\mapsto X^{(\varepsilon)}-x_0$ with state space ${\cal W}^m$ satisfies the sample path large deviation principle with speed $\varepsilon^{-1}$ and good rate function $\widetilde{Q}_T$ defined in (\ref{E:QT1}). The validity of the large deviation principle means that for every Borel measurable subset ${\cal A}$ of ${\cal W}^m$, the following estimates hold:
\begin{align*}
&-\inf_{g\in{\cal A}^{\circ}}\widetilde{Q}_T(g)\le\liminf_{\varepsilon\downarrow 0}\varepsilon\log\mathbb{P}\left(
X^{(\varepsilon)}-x_0\in{\cal A}\right) 
\\
&\le\limsup_{\varepsilon\downarrow 0}\varepsilon\log\mathbb{P}\left(X^{(\varepsilon)}-x_0\in{\cal A}\right)
\le-\inf_{g\in\bar{{\cal A}}}\widetilde{Q}_T(g).
\end{align*}
The symbols ${\cal A}^{\circ}$ and $\bar{{\cal A}}$ in the previous estimates stand for the interior and the closure of the set 
${\cal A}$, respectively. 
\end{theorem}
\begin{remark}\label{R:csu}
In Theorem \ref{T:nmn}, the canonical set-up on the space ${\cal W}^m\times{\cal W}^m$ is employed. However, if the LDP in Theorem 
\ref{T:beg11}
is valid on the set-up associated with the volatility process in (\ref{E:moodk}), then Theorem \ref{T:nmn} holds on the set-up associated with the equation in (\ref{E:moodk}). This follows from the proof of Theorem \ref{T:nmn} given in Section \ref{S:prov}. Examples illustrating the statement formulated above are volatility processes in multivariate 
Gaussian models (see Theorem \ref{T:dert} in Subsection \ref{SS:mGm}), volatility processes in multivariate non-Gaussian fractional models
(see Theorem \ref{T:riur} in Subsection \ref{SS:mGm}), 
and the processes solving Volterra type stochastic integral equations employed in the paper \cite{NR} of Nualart and Rovira (see Theorem 1 in \cite{NR}, see also Subsection \ref{SS:NR} in the present paper).
\end{remark}

Suppose that for every $(t,u)\in[0,T]\times\mathbb{R}^d$, the matrix $\sigma(t,u)$ is invertible. Then, the expression on the 
right-hand side of (\ref{E:QT1}) can be simplified.
By taking into account (\ref{E:rtw}), we obtain the following equality:
$\dot{l}(t)=\bar{C}^{-1}\sigma(t,\widehat{f}(t))^{-1}[\dot{g}(t)-b(t,\widehat{f}(t))-\sigma(t,\widehat{f}(t))
C\dot{f}(t)]$, $t\in[0,T]$.
Hence for all $g\in(\mathbb{H}_0^1)^m$,
\begin{align}
&\widetilde{Q}_T(g) \nonumber \\
&=\frac{1}{2}\inf_{f\in(\mathbb{H}_0^1)^m}\int_0^T(||\bar{C}^{-1}\sigma(s,\widehat{f}(s))^{-1}
[\dot{g}(s)-b(s,\widehat{f}(s))-\sigma(s,\widehat{f}(s))C\dot{f}(s)]||_m^2
+||\dot{f}(s)||_m^2)ds,
\label{E:ras}
\end{align}
and $\widetilde{Q}_T(g)=\infty$ otherwise.
\begin{remark}\label{Rttx}
Let $g\in(\mathbb{H}_0^1)^m$. Then, for every fixed function $f\in(\mathbb{H}_0^1)^m$, the integral on the right-hand side of (\ref{E:ras}) is finite. Indeed, it suffices to prove that
\begin{equation}
\sup_{t\in[0,T]}||\sigma(s,\widehat{f}(s))^{-1}||_{m\times m}<\infty.
\label{E:ree}
\end{equation}
The matrix-valued function $s\mapsto\sigma(s,\widehat{f}(s))$, with $s\in[0,T]$, is continuous. Therefore, the function
$s\mapsto|\det(\sigma(s,\widehat{f}(s)))|$ is continuous on $[0,T]$ and bounded away from zero. Next, using the uniform boundedness of the adjugate matrices associated with the matrices $\sigma(s,\widehat{f}(s))$, $s\in[0,T]$, we see that the inequality in (\ref{E:ree}) holds.
\end{remark}

The next statement concerns the continuity of the rate function $\widetilde{Q}_T$.
\begin{lemma}\label{L:dff}
Suppose that for every $(t,u)\in[0,T]\times\mathbb{R}^d$, the matrix $\sigma(t,u)$ is invertible. Then, the function $\widetilde{Q}_T$
defined in (\ref{E:ras}) is continuous in the topology of the space $(\mathbb{H}_0^1)^m$.
\end{lemma}

\it Proof. \rm The lower semicontinuity of the function $\widetilde{Q}_T$ on the space $(\mathbb{H}_0^1)^m$ can be established using
the following two facts: (1)\,\,Since $\widetilde{Q}_T$ is a good rate function on the space $\mathbb{C}_0^m$, it is lower semicontinuous on this space. (2)\,\,The space $(\mathbb{H}_0^1)^m$ is continuously embedded into the space $\mathbb{C}_0^m$.

We will next show that the function $\widetilde{Q}_T$ is upper semicontinuous on the space $(\mathbb{H}_0^1)^m$. By taking into account the representation in (\ref{E:QT1}) and the fact that the greatest lower bound of any family of upper semicontinuous functions is an upper semicontinuous function, we see that it suffices to prove that for any $f\in(\mathbb{H}_0^1)^m$, the function
\begin{equation}
g\mapsto\frac{1}{2}\int_0^T(||\bar{C}^{-1}\sigma(s,\widehat{f}(s))^{-1}
[\dot{g}(s)-b(s,\widehat{f}(s))-\sigma(s,\widehat{f}(s))C\dot{f}(s)]||_m^2
+||\dot{f}(s)||_m^2)ds
\label{E:eet}
\end{equation}
is continuous on the space $(\mathbb{H}_0^1)^m$. The previous statement follows from formula (\ref{E:ree}). 

The proof of Lemma \ref{L:dff} is thus completed.

Let $m=1$, and suppose the volatility function $\sigma$ satisfies the condition $\sigma(t,u)\neq 0$ for all 
$(t,u)\in[0,T]\times\mathbb{R}^d$. Then, for every $g\in\mathbb{H}_0^1$, the formula in (\ref{E:ras}) can be rewritten as follows:
\begin{align}
&\widetilde{Q}_T(g) 
=\frac{1}{2}\inf_{f\in\mathbb{H}_0^1}\int_0^T\left[\frac{(\dot{g}(s)-b(s,\widehat{f}(s))-\rho\sigma(s,\widehat{f}(s))\dot{f}(s))^2}
{(1-\rho^2)\sigma(s,\widehat{f}(s))^2}
+\dot{f}(s)^2\right]ds.
\label{E:QT3}
\end{align}
We also have 
\begin{equation}
\widetilde{Q}_T(g)=\infty\quad\mbox{if}\quad g\in\mathbb{C}\backslash\mathbb{H}_0^1.
\label{E:qeer}
\end{equation}
\section{Small-Noise LDPs for Log-Price Processes}\label{S:snl}
Our goal in this section is to derive from Theorem \ref{T:nmn} a small-noise LDP for the process 
$\varepsilon\mapsto X^{(\varepsilon)}_T-x_0$ . In financial mathematics, such LDPs are used in small-noise asymptotic analysis of various path-independent options.

Consider a map $V:\mathbb{C}_0^m\mapsto\mathbb{R}^m$ defined by
$V(\varphi)=\varphi(T)$. It is clear that this map is continuous. It follows from Theorem \ref{T:nmn} and the contraction principle that the following assertion holds.
\begin{theorem}\label{T:rrr}
Under the restrictions in Theorem \ref{T:nmn}, the process $\varepsilon\mapsto X^{(\varepsilon)}_T-x_0$ with state space $\mathbb{R}^m$ satisfies the small-noise large deviation principle 
with speed $\varepsilon^{-1}$ 
and good rate function $\widehat{I}_T(x)$, $x\in\mathbb{R}^m$ given by
\begin{align}
&\widehat{I}_T(x)=\inf_{\{g\in(\mathbb{H}_0^1)^m:g(T)=x\}}\widetilde{Q}_T(g) 
\label{E:l}
\end{align}
where the rate function $\widetilde{Q}_T$ is defined in (\ref{E:QT1}).
\end{theorem}

The expression on the right-hand side of (\ref{E:l}) is rather complicated. Our next goal is to show that if $m=1$, then the formula in
(\ref{E:l}) can be simplified (see Theorem \ref{T:bed}). We do not know whether a similar simplification is possible for $m> 1$. Here our knowledge is fragmentary. More information will be given below.

Let $m=1$. For all $y\in\mathbb{R}$ and $f\in\mathbb{H}_0^1$, set
$$
\Psi(y,f,\widehat{f})=\int_0^Tb(s,\widehat{f}(s))ds+\rho\int_0^T\sigma(s,\widehat{f}(s))\dot{f}(s)ds
+\bar{\rho}\left\{\int_0^T\sigma(s,\widehat{f}(s))^2ds\right\}^{\frac{1}{2}}y.
$$
Define the function $\widetilde{I}_T$ on $\mathbb{R}$ by
\begin{equation}
\widetilde{I}_T(x)=\frac{1}{2}\inf_{y\in\mathbb{R},f\in\mathbb{H}_0^1}\left\{y^2+\int_0^T\dot{f}(t)^2dt:
\Psi(y,f,\widehat{f})=x\right\},
\label{E:tri}
\end{equation}
if the equation $\Psi(y,f,\widehat{f})=x$ is solvable for $y$ and $f$, and $\widetilde{I}_T(x)=\infty$ otherwise.
\begin{remark}\label{R:biu}
The equation $\Psi(y,f,\widehat{f})=x$ is as follows:
\begin{equation}
\int_0^Tb(s,\widehat{f}(s))ds+\rho\int_0^T\sigma(s,\widehat{f}(s))\dot{f}(s)ds
+\bar{\rho}\left\{\int_0^T\sigma(s,\widehat{f}(s))^2ds\right\}^{\frac{1}{2}}y=x.
\label{E:fre}
\end{equation}
\end{remark}

\begin{theorem}\label{T:bed}
Suppose $m=1$, and the conditions in Theorem \ref{T:nmn} hold. Then, for all $x\in\mathbb{R}$, $\widetilde{I}_T(x)=\widehat{I}_T(x)$.
\end{theorem}

\it Proof. \rm We will first prove that 
$\widetilde{I}_T(x)\le\widehat{I}_T(x)$ for $x\in\mathbb{R}$.
Fix $x\in\mathbb{R}$ and $g\in\mathbb{H}_0^1$ with $g(T)=x$. If the equation 
\begin{equation}
\Phi(l,f,\widehat{f})(t)=g(t),\,t\in[0,T],
\label{E:ewr}
\end{equation}
is not solvable for $l$ and $f$, then $\widetilde{Q}_T(g)=\infty$. If the equation in (\ref{E:ewr}) is solvable and $f$ is such 
that $\int_0^T\sigma(s,\widehat{f}(s))^2ds=0$,
then (\ref{E:disc}) implies the equality $\int_0^Tb(s,\widehat{f}(s))ds=x$.
Hence the equation $\Psi(y,f,\widehat{f})=x$, where $f$ is the function mentioned above, holds with any $y\in\mathbb{R}$.
Next, suppose the equation in (\ref{E:ewr}) is solvable, and the function $f$ is such that 
$\int_0^T\sigma(s,\widehat{f}(s))^2ds> 0$.
It follows that the equation
$\Psi(y,f,\widehat{f})=x$, with the same function $f$, is uniquely solvable for $y$,
and the unique solution satisfies 
$$
\int_0^T\sigma(s,\widehat{f}(s))\dot{l}(s)ds
=\left\{\int_0^T\sigma(s,\widehat{f}(s))^2ds\right\}^{\frac{1}{2}}y.
$$
Moreover, $y^2\le\int_0^T\dot{l}(s)^2ds$.
Finally, we see that the reasoning above shows that the inequality $\widetilde{I}_T(x)\le\widehat{I}_T(x)$ holds for all $x\in\mathbb{R}$.

We will next prove that the opposite inequality holds as well.
Let $x\in\mathbb{R}$, and suppose the equation 
\begin{equation}
\Psi(y,f,\widehat{f})=x
\label{E:kpl}
\end{equation} 
is not solvable for $y$ and $f$. Then, we have $\widetilde{I}_T(x)=\infty$, and hence $\widehat{I}_T(x)\le\widetilde{I}_T(x)$. 
Next, suppose the equation in (\ref{E:kpl}) is solvable for $y$ and $f$, and the function $f$ satisfies the condition
$\int_0^T\sigma(s,\widehat{f}(s))^2ds=0$. Then, $x=\int_0^Tb(s,\widehat{f}(s)ds$.
Set $g(t)=\int_0^tb(s,\widehat{f}(s))ds$. It is not hard to see that $g\in\mathbb{H}_0^1$, $g(T)=x$, and (\ref{E:ewr}) holds
for $g$, with the same function $f$ and any function $l\in\mathbb{H}_0^1$. Let us next assume that the equation in (\ref{E:kpl})
is solvable for $y$ and $f$, and the function $f$ satisfies $\int_0^T\sigma(s,\widehat{f}(s))^2ds> 0$. Then, there exists 
$l\in\mathbb{H}_0^1$ such that 
$$
\dot{l}(s)=\left\{\int_0^T\sigma(u,\widehat{f}(u))^2du\right\}^{-\frac{1}{2}}\sigma(s,\widehat{f}(s))y
$$
for all $s\in[0,T]$. It follows that $\int_0^T\dot{l}(s)^2ds=y^2$. Set
$$
g(t)=\int_0^tb(s,\widehat{f}(s))ds+\rho\int_0^t\sigma(s,\widehat{f}(s))\dot{f}(s)ds
+\bar{\rho}\int_0^t\sigma(s,\widehat{f}(s))\dot{l}(s)ds,
$$
with the same functions $l$ and $f$ as above. Then, $g\in\mathbb{H}_0^1$. It is not hard to see that
\begin{align*}
g(T)&=\int_0^Tb(s,\widehat{f}(s))ds+\rho\int_0^T\sigma(s,\widehat{f}(s))\dot{f}(s)ds 
+\bar{\rho}\left\{\int_0^T\sigma(s,\widehat{f}(s))^2ds\right\}^{\frac{1}{2}}y \\
&=\Psi(y,f,\widehat{f})=x.
\end{align*}
In addition, the functions $l$ and $f$ solve the equation in (\ref{E:ewr}). Finally, 
summarizing what was said above, we see that $\widehat{I}_T(x)\le\widetilde{I}_T(x)$ for all $x\in\mathbb{R}$.

This completes the proof of Theorem \ref{T:bed}.
\begin{remark}\label{R:begin}
The expression for the rate function $\widetilde{I}_T$ on the right-hand side of (\ref{E:tri}) can be given a simpler form. 
Note that for every $f\in\mathbb{H}_0^1$, the function $s\mapsto\sigma(s,\widehat{f}(s))$ is continuous. Define
$$
Q_1=\{f\in\mathbb{H}_0^1:\sigma(s,\widehat{f}(s))\neq 0\,\,\mbox{for at least one}\,\,s\in[0,T]\}
$$
and
$
Q_2=\{f\in\mathbb{H}_0^1:\sigma(s,\widehat{f}(s))=0\,\,\mbox{for all}\,\,s\in[0,T]\}. 
$
It is clear that $\mathbb{H}_0^1=Q_1\cup Q_2$. Set 
$
Q_3(x)=\{f\in Q_2:x=\int_0^Tb(s,\widehat{f}(s))ds\}.
$
It is not hard to see that (\ref{E:tri}) implies the following equality: 
$
\widetilde{I}_T(x)=\min\{A_1(x),A_2(x)\}
$ 
where
$$
A_1(x)=\frac{1}{2}\inf_{f\in Q_1}\left[\frac{(x-\int_0^Tb(s,\widehat{f}(s))ds
-\rho\int_0^T\sigma(s,\widehat{f}(s))\dot{f}(s)ds)^2}{\bar{\rho}^2\int_0^T\sigma(s,\widehat{f}(s))^2ds}
+\int_0^T\dot{f}(t)^2dt\right]
$$
and
$$
A_2(x)=\begin{cases}
\frac{1}{2}\inf_{\{f\in Q_3(x)\}}\int_0^T\dot{f}(t)^2dt & \text{if $Q_3(x)\neq\emptyset$} \\
\infty & \text{if $Q_3(x)=\emptyset$}.
\end{cases}
$$
In a special case, where $\sigma(s,z)\neq 0$ for all $(s,z)\in[0,T]\times\mathbb{R}^d$, we have $Q_2=\emptyset$, and hence 
\begin{equation}
\widetilde{I}_T(x)=\frac{1}{2}\inf_{f\in\mathbb{H}_0^1}\left[\frac{(x-\int_0^Tb(s,\widehat{f}(s))ds
-\rho\int_0^T\sigma(s,\widehat{f}(s))\dot{f}(s)ds)^2}{\bar{\rho}^2\int_0^T\sigma(s,\widehat{f}(s))^2ds}
+\int_0^T\dot{f}(t)^2dt\right]
\label{E:oi}
\end{equation}
for all $x\in\mathbb{R}$.
\end{remark}

The next statement follows from Theorem \ref{T:rrr}, Theorem \ref{T:bed}, and (\ref{E:oi}).
\begin{theorem}\label{T:ha}
Suppose Assumption A and Assumptions (C1) -- (C7) hold true, and the model in (\ref{E:moodk}) is defined on the canonical set-up.
Suppose also that $\sigma(s,z)\neq 0$ for all $(s,z)\in[0,T]\times\mathbb{R}^d$. Then, the process 
$\varepsilon\mapsto X^{(\varepsilon)}_T-x_0$ with state space $\mathbb{R}^1$ satisfies the small-noise large deviation principle 
with speed $\varepsilon^{-1}$ 
and good rate function $\widehat{I}_T(x)$, $x\in\mathbb{R}^m$ given by the formula in (\ref{E:oi}).
\end{theorem}

The next lemma concerns the continuity of the rate function in Theorem \ref{T:ha}.
\begin{lemma}\label{L:poi}
Suppose $\sigma(s,z)\neq 0$ for all $(s,z)\in[0,T]\times\mathbb{R}^d$. Then, the function $\widetilde{I}_T$ is continuous on $\mathbb{R}$.
\end{lemma}

\it Proof. \rm The function $\widetilde{I}_T$ is lower semicontinuous since it is a good rate function. Moreover, it is upper semicontinuous
being equal to the greatest upper bound of a family of continuous functions on $\mathbb{R}$.

This completes the proof of Lemma \ref{L:poi}.

Next, we turn our attention to the case where $m> 1$. As we have already mentioned, our knowledge in this case is incomplete. 
First of all, it is not clear how to choose the map $\Psi$. One of the acceptable candidates is as follows:
\begin{align}
\Psi(y,f,\widehat{f})&=\int_0^Tb(s,\widehat{f}(s))ds+\int_0^T\sigma(s,\widehat{f}(s))C\dot{f}(s)ds  \nonumber \\
&\quad+\left\{\int_0^T||\sigma(s,\widehat{f}(s))\bar{C}||_{*}^2ds\right\}^{\frac{1}{2}}y,\quad y\in\mathbb{R}^m,\quad
f\in(\mathbb{H}_0^1)^m,
\label{E:rere}
\end{align}
where $||\cdot||_{*}$ is the operator norm on the space of $m\times m$-matrices. The usefulness of the operator norm employed in 
(\ref{E:rere}) will be clear below (see the proof of Lemma \ref{L:bedd}). Let us also define a map $\widetilde{I}_T:\mathbb{R}^m\mapsto\mathbb{R}$ by
\begin{equation}
\widetilde{I}_T(x)=\frac{1}{2}\inf_{y\in\mathbb{R}^m,f\in(\mathbb{H}_0^1)^m}\left\{||y||_m^2+\int_0^T||\dot{f}(t)||_m^2dt:
\Psi(y,f,\widehat{f})=x\right\},
\label{E:tric}
\end{equation}
if the equation $\Psi(y,f,\widehat{f})=x$ is solvable for $y$ and $f$, and by $\widetilde{I}_T(x)=\infty$ otherwise.
\begin{lemma}\label{L:pp}
For every $x\in\mathbb{R}^m$, the inequality
$\widetilde{I}_T(x)\le\widehat{I}_T(x)$ holds for the functions defined in (\ref{E:l}) and (\ref{E:tric}), respectively.
\end{lemma}

\it Proof. \rm The proof of lemma \ref{L:pp} is similar to that of the first part of Theorem \ref{T:bed}. We include this proof for the sake of completeness.
Fix $x\in\mathbb{R}^m$ and $g\in(\mathbb{H}_0^1)^m$, with $g(T)=x$. If the equation 
\begin{equation}
\Phi(l,f,\widehat{f})(t)=g(t),\,t\in[0,T],
\label{E:ewrs}
\end{equation}
is not solvable for $l$ and $f$, then $\widetilde{Q}_T(g)=\infty$. If the equation in (\ref{E:ewrs}) is solvable, and the function $f$ is such that
$\int_0^T||\sigma(s,\widehat{f}(s))||_{*}^2ds=0$,
then for every $s\in[0,T]$, $\sigma(s,\widehat{f}(s))=0$.
It follows from (\ref{E:disc}) that $\int_0^Tb(s,\widehat{f}(s))ds=x$.
Hence the equation $\Psi(y,f,\widehat{f})=x$, where $f$ is the function mentioned above, holds with any $y\in\mathbb{R}^m$.
Next, suppose the equation in (\ref{E:ewrs}) is solvable and 
$\int_0^T||\sigma(s,\widehat{f}(s))||_{*}^2ds> 0$. 
Then, we have
$\int_0^T||\sigma(s,\widehat{f}(s))\bar{C}||_{*}^2ds> 0$.
Hence the equation 
$\Psi(y,f,\widehat{f})=x$, with the same function $f$ as above, is uniquely solvable for $y$,
and the unique solution satisfies 
\begin{equation}
\int_0^T\sigma(s,\widehat{f}(s))\bar{C}\dot{l}(s)ds
=\left\{\int_0^T||\sigma(s,\widehat{f}(s))\bar{C}||_{*}^2ds\right\}^{\frac{1}{2}}y.
\label{E:erd}
\end{equation}
Now, it is not hard to prove using H\"{o}lder's inequality that the equality in (\ref{E:erd}) implies the estimate
$||y||_m^2\le\int_0^T||\dot{l}(s)||_m^2ds$.

Finally, we see that the reasoning above shows that the inequality $\widetilde{I}_T(x)\le\widehat{I}_T(x)$ holds for all $x\in\mathbb{R}$.

This completes the proof of Lemma \ref{L:pp}.

Our next goal is to provide examples of volatility maps $\sigma$ such that $\widetilde{I}_T(x)=\widehat{I}_T(x)$
for all $x\in\mathbb{R}^m$.
\begin{lemma}\label{L:bedd}
Let $m> 1$, and suppose that for all $t\in[0,T]$ and $z\in\mathbb{R}^d$, the volatility map is given by 
$\sigma(t,z)=\xi(t,z)O(t,z)\bar{C}^{-1}$ 
where $O(t,z)$ are orthogonal $m\times m$-matrices and
$\xi$ is a real function on $[0,T]\times\mathbb{R}^d$. Let us also assume that the function $\xi$ and the map $O$ are continuous on $[0,T]\times\mathbb{R}^m$. Then, $\widetilde{I}_T(x)=\widehat{I}_T(x)$ for all $x\in\mathbb{R}^m$.
\end{lemma}

\it Proof. \rm It suffices to show that $\widehat{I}_T(x)\le\widetilde{I}_T(x)$
for $x\in\mathbb{R}^m$ (see Lemma \ref{L:pp}). Under the conditions in Lemma \ref{L:bedd}, we have
\begin{equation}
\Psi(y,f,\widehat{f})=\int_0^Tb(s,\widehat{f}(s))ds+\int_0^T\sigma(s,\widehat{f}(s))C\dot{f}(s)ds
+\left\{\int_0^T\xi(s,\widehat{f}(s))^2ds\right\}^{\frac{1}{2}}y.
\label{E:pferd}
\end{equation}
The equality in (\ref{E:pferd}) can be established by using the fact that for any orthogonal matrix $O$, $||O||_{*}=1$.

Let $x\in\mathbb{R}^m$, and suppose the equation 
$\Psi(y,f,\widehat{f})=x$ is not solvable for $y$ and $f$. Then $\widetilde{I}_T(x)=\infty$, and hence 
$\widehat{I}_T(x)\le\widetilde{I}_T(x)$. 
Now, suppose the equation $\Psi(y,f,\widehat{f})=x$ is solvable for $y$ and $f$, and the condition
$\int_0^T\xi(s,\widehat{f}(s))^2ds=0$ is satisfied. Then, we have $x=\int_0^Tb(s,\widehat{f}(s)ds$.
Set $g(t)=\int_0^tb(s,\widehat{f}(s))ds$. It is not hard to see that $g\in(\mathbb{H}_0^1)^m$, $g(T)=x$, and (\ref{E:ewrs}) holds
for $g$, with the same function $f$ and any function $l\in(\mathbb{H}_0^1)^m$. Let us next assume that the equation 
$\Psi(y,f,\widehat{f})=x$ is solvable for $y$ and $f$, and the condition $\int_0^T\xi(s,\widehat{f}(s))^2ds> 0$ is satisfied. 
Then, there exists $l\in(\mathbb{H}_0^1)^m$ such that 
\begin{equation}
\dot{l}(s)=\left\{\int_0^T\xi(u,\widehat{f}(u))^2du\right\}^{-\frac{1}{2}}\xi(s,\widehat{f}(s))O(s,\widehat{f}(s))^{\prime}y
\label{E:qwa}
\end{equation}
for all $s\in[0,T]$. It follows that
\begin{align}
&\int_0^T\sigma(s,\widehat{f}(s))\bar{C}\dot{l}(s)ds=\int_0^T\xi(s,\widehat{f}(s))O(s,\widehat{f}(s))\dot{l}(s)ds
=\left\{\int_0^T\xi(s,\widehat{f}(s))^2ds\right\}^{\frac{1}{2}}y.
\label{E:kutt}
\end{align}
Moreover, using (\ref{E:qwa}) and the fact that for every orthogonal matrix $O$ and $y\in\mathbb{R}^m$, 
the equality $||Oy||_m^2=||y||_m^2$ holds, we obtain the equality 
$\int_0^T||\dot{l}(s)||_m^2ds=||y||_m^2$. Set
$$
g(t)=\int_0^tb(s,\widehat{f}(s))ds+\int_0^t\sigma(s,\widehat{f}(s))C\dot{f}(s)ds
+\int_0^t\sigma(s,\widehat{f}(s))\bar{C}\dot{l}(s)ds
$$
where the functions $l$ and $f$ are as above. Then we have $g\in(\mathbb{H}_0^1)^m$. Next, using (\ref{E:kutt}) we obtain
\begin{align*}
g(T)&=\int_0^Tb(s,\widehat{f}(s))ds+\int_0^T\sigma(s,\widehat{f}(s))C\dot{f}(s)ds 
+\left\{\int_0^T\xi(s,\widehat{f}(s))^2ds\right\}^{\frac{1}{2}}y \\
&=\Psi(y,f,\widehat{f})=x.
\end{align*}
In addition, the functions $l$ and $f$ solve the equation $\Phi(l,f,\widehat{f})(t)=g(t)$, $t\in[0,T]$. Finally, 
by taking into account the reasoning above, we see that $\widehat{I}_T(x)\le\widetilde{I}_T(x)$ for all $x\in\mathbb{R}$.

This completes the proof of Lemma \ref{L:bedd}.
\begin{remark}\label{R:end}
Suppose the conditions in Lemma \ref{L:bedd} hold, and the function $\xi$ satisfies $\xi(t,z)> 0$ for all 
$(t,z)\in[0,T]\times\mathbb{R}^d$. Then, we have
\begin{align*}
&\widetilde{I}_T(x)  \\
&=\frac{1}{2}\inf_{f\in(\mathbb{H}_0^1)^m}\left[\frac{
||x-\int_0^Tb(s,\widehat{f}(s))ds-\int_0^T\xi(s,\widehat{f}(s))O(s,\widehat{f}(s))\bar{C}^{-1}C\dot{f}(s)ds||_m^2}
{\int_0^T\xi(s,\widehat{f}(s))^2ds}+\int_0^T||\dot{f}(t)||_m^2dt\right].
\end{align*}
The previous formula can be established by taking into account (\ref{E:tri}), (\ref{E:pferd}), and reasoning as in Remark \ref{R:begin}.
\end{remark}
\begin{remark}\label{R:uut}
Suppose the conditions in Remark \ref{R:end} hold for the stochastic model considered in Lemma \ref{L:bedd}. Suppose also that the model is uncorrelated, that is, the condition $C=0$ holds. Then, the process $S$ satisfies the following equation:
\begin{equation}
dS_t=S_t\circ[b(t,\widehat{B}_t)dt+\xi(t,\widehat{B}_t)O(t,\widehat{B}_t)dW_t],\quad 0\le t\le T,\quad S_0=s_0\in\mathbb{R}^m.
\label{E:utr}
\end{equation}
In addition, the rate function $\widetilde{I}_T$ in Remark \ref{R:end} is given by
\begin{align}
&\widetilde{I}_T(x)=\frac{1}{2}\inf_{f\in(\mathbb{H}_0^1)^m}\left[\frac{
||x-\int_0^Tb(s,\widehat{f}(s))ds||_m^2}
{\int_0^T\xi(s,\widehat{f}(s))^2ds}+\int_0^T||\dot{f}(t)||_m^2dt\right].
\label{E:ko}
\end{align}
Note that although the model in (\ref{E:utr}) depends on the family $O$ of orthogonal matrices used in Lemma \ref{L:bedd}, 
the rate function $\widetilde{I}_T$ given by the expression in (\ref{E:ko}) is independent of that family.
\end{remark}

We will next provide a special example of a model, for which the formula in Remark \ref{R:end} holds true. Similar more complicated 
models can also be constructed.
Let $d=1$, $m=2$, and choose
$
C=\begin{bmatrix}
\frac{1}{2} & 0 \\
0 & \frac{1}{2}
\end{bmatrix}.
$
Then $\bar{C}=\begin{bmatrix}
\frac{1}{2}\sqrt{3} & 0 \\
0 & \frac{1}{2}\sqrt{3}
\end{bmatrix}
$ 
and
$\bar{C}^{-1}=\begin{bmatrix}
\frac{2\sqrt{3}}{3} & 0 \\
0 & \frac{2\sqrt{3}}{3}
\end{bmatrix}.
$ 
Let $b$ be a drift map satisfying Assumption A, and suppose $\xi(t,z)$ is a real strictly positive $\omega$-continuous function on $[0,T]\times\mathbb{R}$. Consider the following family of orthogonal $2\times 2$-matrices: 
$
O(z)=\begin{bmatrix}
\cos z & -\sin z \\
\sin z & \cos z
\end{bmatrix},
$
where $z\in\mathbb{R}$, and set 
$$
\sigma(t,z)=\xi(t,z)O(z)\bar{C}^{-1},\quad(t,z)\in[0,T]\times\mathbb{R}.
$$ 
Then, the stochastic volatility model
in (\ref{E:moodk}) takes the following form:
$$
dS_t=S_t\circ\left(b(t,\widehat{B}_t)ds+\xi(t,\widehat{B}_t)\begin{bmatrix}
\cos\widehat{B}_t & -\sin\widehat{B}_t \\
\sin\widehat{B}_t & \cos\widehat{B}_t
\end{bmatrix}(dW_t+\frac{\sqrt{3}}{3}dB_t)\right),
$$
where $0\le t\le T$ and $S_0=s_0\in\mathbb{R}^2$. Here we use the equality $\bar{C}^{-1}C=\begin{bmatrix}
\frac{\sqrt{3}}{3} & 0 \\
0 & \frac{\sqrt{3}}{3}
\end{bmatrix}.
$
It is easy to see that the rate function in Remark \ref{R:uut} is given by
\begin{align*}
\widetilde{I}_T(x)&=\frac{1}{2}\inf_{f\in(\mathbb{H}_0^1)^m}\{\frac{
||x-\int_0^Tb(s,\widehat{f}(s))ds-\frac{\sqrt{3}}{3}\int_0^T\xi(s,\widehat{f}(s))\begin{bmatrix}
\cos\widehat{f}(s) & -\sin\widehat{f}(s) \\
\sin\widehat{f}(s) & \cos\widehat{f}(s)
\end{bmatrix}\dot{f}(s)ds||_m^2}
{\int_0^T\xi(s,\widehat{f}(s))^2ds} \\
&\quad+\int_0^T||\dot{f}(t)||_m^2dt\}.
\end{align*}

\section{Large Deviation Principles for Volatility Processes}\label{S:im}
The main result of this section is Theorem \ref{T:beg1} that provides a LDP for the solution to the equation in 
(\ref{E:2ooo}). Theorem \ref{T:beg1} uses the canonical set-up.
\begin{theorem}\label{T:beg1}
Suppose Assumptions (C1) -- (C7) hold, and let $Y^{(\varepsilon)}$ with $Y_0^{(\varepsilon)}=y$ be the process solving 
the equation in (\ref{E:2ooo}). Then, the process
$Y^{(\varepsilon)}$ satisfies a sample path large deviation principle with speed $\varepsilon^{-1}$ and good rate function defined on
${\cal W}^d$ by
\begin{equation}
I_y(\varphi)=\inf_{\{f\in L^2([0,T],\mathbb{R}^m):\,\Gamma_y(f)=\varphi\}}\frac{1}{2}\int_0^T||f(t)||^2_mdt
\label{E:LDP}
\end{equation}
if $\{f\in L^2([0,T],\mathbb{R}^m):\,\Gamma_y(f)=\varphi\}\neq\emptyset$, and $I_y(\varphi)=\infty$ otherwise.
\end{theorem}
\begin{remark}\label{R:rate}
The goodness of the rate function $I_y$ can be shown as follows. Consider sublevel sets of $I_y$ given by
$L_c=\{\varphi\in{\cal W}^d:\,I_y(\varphi)\le c\}$ for $c> 0$.
Then, we have $L_c=\cap_{\varepsilon> 0}\Gamma_y(D_{2c+\varepsilon})$. Every set $D_{2c+\varepsilon}$ is compact in the weak topology of the space $L^2([0,T],\mathbb{R}^m)$. It follows from Assumption (C5) that the set $\Gamma_y(D_{2c+\varepsilon})$ is compact in the space 
${\cal W}^d$. Therefore, the set $L_c$ is compact in ${\cal W}^d$ since this set can be represented as the intersection of compacts sets.
\end{remark}
\begin{corollary}\label{C:67}
Suppose Assumptions (C1) -- (C7) hold, and let $\widehat{B}^{\varepsilon}$ be the volatility process (see Definition \ref{D:vo}). Then,
the process $\widehat{B}^{\varepsilon}$ satisfies a sample path LDP with speed $\varepsilon^{-1}$ and good rate function given 
for $\varphi\in{\cal W}^d$ by
\begin{equation}
J_y(\varphi)=\inf_{\{f\in L^2([0,T],\mathbb{R}^m):\,{\cal A}f=\varphi\}}\frac{1}{2}\int_0^T||f(t)||^2_mdt
\label{E:LDPs}
\end{equation}
if $\{f\in L^2([0,T],\mathbb{R}^m):\,{\cal A}f=\varphi\}\neq\emptyset$, and $J_y(\varphi)=\infty$ otherwise. In (\ref{E:LDPs}), ${\cal A}$
is the map introduced in Definition \ref{D:lip}.
\end{corollary}

Corollary \ref{C:67} follows from Theorem \ref{T:beg1} and the contraction principle.

The following assertion can be derived from Theorem \ref{T:beg1}. It concerns a sample path large deviation principle for the process
\begin{equation}
\varepsilon\mapsto(\sqrt{\varepsilon}W,\sqrt{\varepsilon}B,\widehat{B}^{\varepsilon}),\quad \varepsilon\in(0,1]
\label{E:lab}
\end{equation} 
where $W$ and $B$ are independent $m$-dimensional Brownian motions appearing in (\ref{E:mood}), while $\widehat{B}^{\varepsilon}$ is the scaled volatility process 
(see Definition \ref{D:vo}). The state space of the process in (\ref{E:lab}) is ${\cal W}^m\times{\cal W}^m\times{\cal W}^d$.
\begin{theorem}\label{T:beg11}
Under the restrictions in Theorem \ref{T:beg1}, the process in (\ref{E:lab})
satisfies a sample path large deviation principle with speed $\varepsilon^{-1}$ and good rate function defined on 
${\cal W}^m\times{\cal W}^m\times{\cal W}^d$ by
\begin{equation}
\widetilde{I}_y(\varphi_1,\varphi_2,\varphi_3)=\frac{1}{2}\int_0^T||\dot{\varphi}_1(t)||_m^2dt
+\frac{1}{2}\int_0^T||\dot{\varphi}_2(t)||_m^2dt
\label{E:LDPa}
\end{equation}
in the case where $\varphi_1,\varphi_2\in (H_0^1)^m$ and $\varphi_3=\widehat{\varphi_2}$, and by 
$\widetilde{I}_y(\varphi_1,\varphi_2,\varphi_3)=\infty$ otherwise.
\end{theorem}
\begin{remark}\label{R:rea}
Theorem \ref{T:beg11} will be an important ingredient in the proof of the LDP in Theorem \ref{T:nmn} given in the next section.
\end{remark}

\it Proof of Theorem \ref{T:beg11}. \rm To derive Theorem \ref{T:beg11} from Theorem \ref{T:beg1}, we use the following $(d+2m)$-dimensional system:
\begin{equation}
\begin{cases}
G_t^{\varepsilon}=\sqrt{\varepsilon}W_t \\
Z_t^{\varepsilon}=\sqrt{\varepsilon}B_t \\
Y_t^{\varepsilon}=y+\int_0^ta(t,s,V^{1,\varepsilon},Y^{\varepsilon})ds
+\sqrt{\varepsilon}\int_0^tc(t,s,V^{2,\varepsilon},Y^{\varepsilon})dB_s.
\end{cases}
\label{E:ju}
\end{equation}
The third equation in (\ref{E:ju}) is the equation in (\ref{E:2ooo}). In addition, $W$ and $B$ are independent $m$-dimensional standard
Brownian motions appearing in (\ref{E:mood}). We can rewrite the model in (\ref{E:ju}) so that Theorem \ref{T:beg1} can be applied to it. The new representation depends on $2m$-dimensional standard Brownian motion $(W,B)$. We also use $2m$-dimensional deterministic controls 
$F=(f_1,f_2)\in L^2([0,T];\mathbb{R}^m)\times L^2([0,T];\mathbb{R}^m)$. It is not hard to see that the system in (\ref{E:t2}) 
becomes the following $(2m+d)$-dimensional system:
$$
\begin{cases}
\eta_1(t)=\int_0^tf_1(u)du \\
\eta_2(t)=\int_0^tf_2(u)du \\
\eta_3(t)=y+\int_0^ta(t,s,\psi_{1,f_2},\eta_3)ds+\int_0^tc(t,s,\psi_{2,f_2},\eta_3)f_2(s)ds.
\end{cases}
$$
The corresponding map $\widetilde{\Gamma}_y$ is given on ${\cal W}^m\times{\cal W}^m$ by
$$
\widetilde{\Gamma}_y(F)(t)=\left(\int_0^tf_1(u)du,\int_0^tf_2(u)du,\Gamma_y(f_2)(t)\right),\,\,t\in[0,T]
$$ 
(see Definition \ref{D:C4}).

Now let $\varphi_1,\varphi_2\in (H_0^1)^m$ and $\varphi_3\in{\cal W}^d$. Then, the equation $\widetilde{\Gamma}_y(F)(t)=
(\varphi_1,\varphi_2,\varphi_3)$ has a unique solution given by $f_1=\dot{\varphi}_1$ and $f_2=\dot{\varphi}_2$. It follows that
$\varphi_3=\Gamma_y(\dot{\varphi}_2)$. Therefore, Theorem \ref{T:beg1} shows that the process 
$\varepsilon\mapsto(\sqrt{\varepsilon}W,\sqrt{\varepsilon}B,Y^{\varepsilon})$ satisfies the large deviation principle 
with speed $\varepsilon^{-1}$ and good rate function $\widetilde{I}_y$ defined in (\ref{E:LDPa}).

Finally, we finish the proof of Theorem \ref{T:beg11} by using the previous reasoning, Definition \ref{D:lip},
and the contraction principle.

Our next goal is to discuss volatility models for which Theorems \ref{T:beg1} and \ref{T:beg11} hold on any set-up. Consider 
the following non-Volterra stochastic differential equation:
\begin{equation}
V_t^{(\varepsilon)}=y+\int_0^ta(s,V^{(\varepsilon)})ds+\sqrt{\varepsilon}\int_0^tc(s,V^{(\varepsilon)})dB_s
\label{E:22a}
\end{equation} 
where $a$ is a map from the space $[0,T]\times{\cal W}^d$ into the space $\mathbb{R}^d$, while $c$ is a map 
from the space $[0,T]\times{\cal W}^d$ into the space of $(d\times m)$-matrices. Let us assume that the maps $a$ and $c$ are locally
Lipschitz and satisfy the sub-linear growth condition (see the definitions in \cite{CF}, A1 and A2 in Section 3, or in \cite{RW}, 
(12.2) and (12.3) on p. 132). A sample path LDP for the process $\varepsilon\mapsto V^{(\varepsilon)}_{\cdot}$ was established on the canonical set-up in \cite{CF}, Theorem 3.1. 
\begin{theorem}\label{T:cans}
Under the restrictions formulated above, the LDP in \cite{CF}, Theorem 3.1 holds on any set-up 
$(\Omega,B,{\cal F}_T^B,\{\mathcal{F}^B_t\}_{0\le t\le T},\mathbb{P})$ (see Definition \ref{D:gens}).
\end{theorem}
\begin{corollary}\label{C:rre}
Suppose the conditions formulated above are satisfied for the equation in (\ref{E:22a}) defined on a general set-up 
$(\Omega,B,{\cal F}_T^B,\{\mathcal{F}^B_t\}_{0\le t\le T},\mathbb{P})$. Then, the LDP in Theorem \ref{T:beg1} holds for
the process $\varepsilon\mapsto V^{(\varepsilon)}$ solving the equation in (\ref{E:22a}), while the LDP in Theorem
\ref{T:beg11} holds for the process $\varepsilon\mapsto(\sqrt{\varepsilon}W,\sqrt{\varepsilon}B,V^{(\varepsilon)})$.
\end{corollary}

\it Proof of Theorem \ref{T:cans}. \rm We have already mentioned that Theorem \ref{T:cans} holds on the canonical set-up (Chiarini and Fischer \cite{CF}). It suffices to show that for every $\varepsilon\in(0,T]$ and $y\in\mathbb{R}^d$, the distribution of the random variable 
$V_{\cdot}^{(\varepsilon)}$, with values in the space ${\cal W}^d$, does not depend on the set-up. By taking into account Theorem 12.1 on p. 132 in \cite{RW}, we establish that the equation in 
(\ref{E:22a}) is pathwise exact (see Definition 9.4 on p. 124 in \cite{RW}). Next, we can use Theorem 10.4 in \cite{RW} to prove that for any
$\varepsilon\in(0,1]$ and any fixed initial condition $y\in\mathbb{R}^d$, there exists a measurable functional 
$F_y^{(\varepsilon)}:{\cal W}^m\mapsto {\cal W}^d$ such that $F_y^{(\varepsilon)}(B)=Y^{\varepsilon}$. The functional $F_y^{(\varepsilon)}$ does not depend on the set-up. Finally, using the previous equality and the fact that the distribution of Brownian motion $B$ with respect to the measure $\mathbb{P}$ is the same for all set-ups, we see that the distribution of $Y^{\varepsilon}$ in ${\cal W}^d$ does not depend on the set-up. It follows that since the LDPs in Theorems \ref{T:beg1} and \ref{T:beg11} hold for the processes 
$\varepsilon\mapsto V^{(\varepsilon)}_{\cdot}$ and $\varepsilon
\mapsto(\sqrt{\varepsilon}W_{\cdot},\sqrt{\varepsilon}B_{\cdot},V^{(\varepsilon)}_{\cdot})$ defined on the canonical set-up, they also hold on any set-up.

This completes the proof of Theorem \ref{T:cans}.
\begin{remark}
The fact that the distribution of the solution does not depend on the set-up was established in \cite{SV} for less general equations than that in (\ref{E:22a}) and under stronger restrictions on the coefficient maps $a$ and $c$ (see the equation in formula (1.3) and Corollary 5.1.3
in \cite{SV}).
\end{remark}

\section{Unification of Sample Path Large Deviation Principles for Stochastic Volatility Models}\label{S:unix}
Our aim in this section is to provide examples of log-price processes and volatility processes to which the LDPs in Theorems 
\ref{T:nmn} and \ref{T:beg11} can be applied. The present section is divided into several subsections. It is organized as follows. In Subsections \ref{SS:GM} and \ref{SS:frH}, we overview one-factor Gaussian models studied in \cite{Gul1} and 
one-factor non-Gaussian fractional models introduced in \cite{GGG}. These subsections are auxiliary. They provide examples of Volterra type 
kernels and processes used in Subsection \ref{SS:fu} devoted to mixtures of multivariate Gaussian models and multivariate non-Gaussian 
fractional models. Theorem \ref{T:cxc} obtained in Subsection \ref{SS:fu} is one of the main results in the present paper. It follows from  Theorem \ref{T:cxc} that the LDPs in Theorems \ref{T:nmn} and \ref{T:beg1} can be applied to log-price processes and volatility processes in mixed models (see Remark \ref{R:rrro}) defined on the canonical set-up. In Subsection \ref{SS:mGm}, we show that Theorem \ref{T:nmn} holds for the log-price process associated with a multivariate Gaussian model defined on any set-up (see Theorem \ref{T:dert}). This theorem is more general than the LDP for multivariate Gaussian models obtained in \cite{CP}. Subsection \ref{SS:mvr}
concerns stochastic volatility models with reflection. Finally, in Subsections \ref{SS:Wang} and \ref{SS:NR}, we discuss the LDPs obtained in \cite{Z} and \cite{NR}.

\subsection{Gaussian Stochastic Volatility Models}\label{SS:GM}
In this subsection, we discuss  Gaussian stochastic volatility models studied in \cite{Gul1}. Let $K$ be a real 
function on $[0,T]^2$. We call the function $K$ an admissible Hilbert-Schmidt kernel if the following conditions hold: 
(a)\,$K$ is Borel measurable on $[0,T]^2$; (b)\,$K$ is Lebesgue square-integrable over $[0,T]^2$;
(c)\,For every $t\in(0,T]$, the slice function $s\mapsto K(t,s)$, with $s\in[0,T]$, 
belongs to the space $L^2[0,T]$;
(d)\,For every $t\in(0,T]$, the slice function is not almost everywhere zero. If an admissible kernel $K$ satisfies the condition 
$K(t,s)=0$ for all $s> t$, then $K$ is called an admissible Volterra kernel. Any such kernel $K$ generates a
Hilbert-Schmidt operator 
\begin{equation}
{\cal K}(f)(t)=\int_0^tK(t,s)f(s)ds,\quad f\in L^2[0,T],\quad t\in[0,T], 
\label{E:intr}
\end{equation}
and a Volterra Gaussian process 
\begin{equation}
\widehat{B}_t=\int_0^tK(t,s)dB_s,\quad t\in[0,T].
\label{E:f11}
\end{equation}
It is clear that the process in (\ref{E:f11}) is adapted to the filtration $\{\mathcal{F}^B_t\}_{0\le t\le T}$. This process is used 
as the volatility process in a Gaussian model. The scaled volatility process is defined as follows: 
$\widehat{B}^{(\varepsilon)}_t=\sqrt{\varepsilon}\widehat{B}_t$ for $t\in[0,T]$. 

In the present paper, only continuous volatility 
processes are used. We will next formulate Fernique's condition guaranteeing that the process in (\ref{E:f11}) is a continuous Gaussian process.
Let $X_t$, $t\in[0,T]$, be a square integrable stochastic process on $(\Omega,{\cal F},\mathbb{P})$. The canonical pseudo-metric $\delta$ associated with this process is defined by the formula
$\delta^2(t,s)=\mathbb{E}[(X_t-X_s)^2]$ for $(t,s)\in[0,T]^2$.
Suppose $\eta$ is a modulus of continuity on $[0,T]$ such that
$\delta(t,s)\le\eta(|t-s|)$ for $t,s\in[0,T]$). Suppose also that for some $b> 1$, the following inequality holds:
\begin{equation}
\int_{b}^{\infty}\eta\left(u^{-1}\right)(\log u)^{-\frac{1}{2}}\frac{du}{u}<\infty.
\label{E:Fer}
\end{equation}
It was announced by Fernique in \cite{F} that a Gaussian process $X$ satisfying the previous condition 
is a continuous stochastic process. The first proof was published by Dudley in \cite{D}
(see also \cite{MS} and the references therein). By the It$\hat{\rm o}$ isometry, the following equality holds 
for the process $\widehat{B}$:
$\delta^2(t,s)=\int_0^T(K(t,u)-K(s,u))^2du$, $t,s\in[0,T]$.
The $L^2$-modulus of continuity of the kernel $K$ is defined on $[0,T]$ by
$$
M_K(\tau)=\sup_{t,s\in[0,T]:|t-s|\le\tau}\int_0^T(K(t,u)-K(s,u))^2du
$$
for all $\tau\in[0,T]$. \\
\\
\it Assumption F. \rm The kernel $K$ in (\ref{E:f11}) is an admissible Volterra kernel such that the estimate 
$M_K(\tau)\le\eta^2(\tau)$, $\tau\in[0,T]$, 
holds for some modulus of continuity $\eta$ satisfying Fernique's condition.
\\
\\
Under Assumption F, the process $\widehat{B}$ defined by (\ref{E:f11}) is a continuous Gaussian process. Moreover, the operator ${\cal K}$ introduced in (\ref{E:intr}) is compact from the space $L^2[0,T]$ into the space 
$\mathbb{C}[0,T]$. 
Important examples of Volterra Gaussian processes are classical fractional processes, e.g., fractional Brownian motion or the Riemann-Liouville fractional Brownian motion. For $0< H< 1$, fractional Brownian motion $B^H$ is a centered Gaussian process with the covariance function given by
$C_H(t,s)=\frac{1}{2}(t^{2H}+s^{2H}-|t-s|^{2H})$, $t,s\ge 0$.
The process $B^H$ was first implicitly considered by Kolmogorov in \cite{Ko}, and was studied by Mandelbrot and van Ness in \cite{MvN}. The constant $H$ is called the Hurst parameter. Fractional Brownian motion is a Volterra type process. This was established 
by Molchan and Golosov (see \cite{MG}, see also \cite{DU}). The corresponding Volterra kernel $K_H$ is as follows: 
For $\frac{1}{2}< H< 1$, 
$$
K_H(t,s)=\sqrt{\frac{H(2H-1)}{\int_0^1(1-x)^{1-2H}x^{H-\frac{3}{2}}dx}}s^{\frac{1}{2}-H}
\int_s^t(u-s)^{H-\frac{3}{2}}u^{H-\frac{1}{2}}du,\quad 0< s< t,
$$
while for $0< H<\frac{1}{2}$,
\begin{align*}
&K_H(t,s)=\sqrt{\frac{2H}{(1-2H)\int_0^1(1-x)^{-2H}x^{H-\frac{1}{2}}dx}} \\
&\quad\left[\left(\frac{t}{s}\right)^{H-\frac{1}{2}}(t-s)^{H-\frac{1}{2}}-\left(H-\frac{1}{2}\right)
s^{\frac{1}{2}-H}\int_s^tu^{H-\frac{3}{2}}(u-s)^{H-\frac{1}{2}}du\right],\quad 0< s< t.
\end{align*}
An equivalent representation of the kernel $K_H$ is the following:
$$
K_H(t,s)=C_H(t-s)^{H-\frac{1}{2}}\left(\frac{s}{t}\right)^{\frac{1}{2}-H}\,_2F_1\left(\frac{1}{2}-H,1,H+\frac{1}{2},\frac{t-s}{t}\right)
$$
where $_2F_1$ is the Gauss hypergeometric function (see, e.g, (3.11) in \cite{CJ}).

The Riemann-Liouville fractional Brownian motion
is defined by the following formula:
\begin{equation}
R^H_t=\Gamma(H+1/2)^{-1}\int_0^t(t-s)^{H-\frac{1}{2}}dB_s,\quad t\ge 0.
\label{E:RL}
\end{equation}
where $0< H< 1$, and the symbol $\Gamma$ stands for the gamma function. This stochastic process was introduced by L\'{e}vy in \cite{PL}. More information about the process $R^H$ can be found in \cite{LSi,Pi}. 
Fractional Brownian motion and the Riemann-Liouville fractional Brownian motion are continuous stochastic processes. 
In the case, where $0< H<\frac{1}{2}$, they are called rough processes since their paths 
are more rough than those of standard Brownian motion. 
In \cite{Gul1}, we introduced a new class of Gaussian stochastic volatility models (the class of super rough models). Super rough models were also considered in \cite{BHP}. In a super rough model, the modulus of continuity associated with the volatility process $\widehat{B}$ growth near zero faster than any power function. Interesting examples of super rough processes can be obtained using Gaussian processes defined in \cite{MV} by Mocioalca and Viens. It was established in \cite{MV} that 
if $\eta\in \mathbb{C}^2(0,T)$ is a modulus of continuity on $[0,T]$
such that the function $x\mapsto(\eta^2)^{\prime}(x)$ is positive and non-increasing on $(0,T)$, then the process
$\widehat{B}_t^{(\eta)}=\int_0^t\tau(t-s)dB_s$, $t\in[0,T]$, with $\tau(x)=\sqrt{(\eta^2)^{\prime}(x)}$, is a Gaussian process satisfying the following conditions: 
(i)\,\,$c_1\eta(|t-s|)\le\delta(t,s)\le c_2\eta(|t-s|)$ for some $c_1,\,c_2> 0$;\,\,(ii)\,\,$X_0=0$;\,\,
(iii)\,\,The process $X$ is adapted to the filtration $\{\mathcal{F}^B_t\}_{0\le t\le T}$. 

A typical example of a modulus of continuity that grows near zero faster than any power is the logarithmic modulus of continuity given by
$$
\eta_{\beta}(x)=\left(\log\frac{1}{x}\right)^{-\frac{\beta}{2}},\quad 0\le x< 1,\quad \beta> 0.
$$ 
It is clear that in this case, the 
function $\tau_{\beta}$ is determined from the equality 
$$
\tau_{\beta}^2(x)=\beta x^{-1}\left(\log\frac{1}{x}\right)^{-\beta-1},\quad 0\le x< 1.
$$ 
In \cite{MV}, the Volterra Gaussian process with the kernel $\tau_{\beta}$ was called logarithmic Brownian motion. If $\beta> 1$, then Fernique's condition is satisfied and the process is continuous. In \cite{Gul1}, we called a Gaussian stochastic volatility model, in which the logarithmic Brownian motion with $\beta> 1$ is the volatility process, a logarithmic Gaussian stochastic volatility model. The logarithmic model is super rough (more details can be found in \cite{Gul1}).
An interesting example of a super rough Gaussian model is the model where the volatility is described by the Wick exponential of a constant multiple of the logarithmic Brownian motion (see \cite{Gul1}). The previous model is similar in structure to the rough Bergomi model introduced in \cite{BFG}, and it may be called the super rough Bergomi model. More details can be found in \cite{Gul1} and \cite{BHP}.
A celebrated Stein and Stein model (see \cite{SS}, see also the discussion in \cite{Gul2}) was one of the first examples of a Gaussian model. 
For a super rough version of the Stein and Stein model see \cite{Gul1}.
\subsection{Non-Gaussian Fractional Stochastic Volatility Models}\label{SS:frH}
This class of one-factor stochastic volatility models was studied in the paper \cite{GGG} of Gerhold, Gerstenecker, and the author. The volatility process in such a model is given by 
$
\hat{B}_t=\int_{0}^{t}K(t,s)U(V_s)ds
$
where $U:\mathbb{R}\mapsto[0,\infty)$ is a continuous non-negative function and $K$ is an admissible kernel with the modulus of continuity in $L^2$ satisfying 
the H\"{o}lder condition. The process $V$ in the formula above is the unique solution to the stochastic differential equation
$$
dV_t=\bar{b}(V_t)dt+\bar{\sigma}(V_t)dB_t,\quad t\in[0,T],
$$ 
where $V_0=v_0> 0$. It is assumed that the following conditions hold (see Section 4.2 in \cite{CF}, see also \cite{GGG}):
\\
(i)\,\,The function $\bar{b}:\mathbb{R}\mapsto\mathbb{R}$ is locally Lipschitz on $\mathbb{R}$, satisfies the sub-linear growth condition,
and $\bar{b}(0)> 0$.
\\
(ii)\,\,The function $\bar{\sigma}:\mathbb{R}\mapsto[0,\infty)$ is locally Lipschitz continuous on $\mathbb{R}-\{0\}$,
satisfies the sub-linear growth condition, $\bar{\sigma}(0)=0$, and $\bar{\sigma}(x)\neq 0$ for all $x\neq 0$. 
\\
(iii)\,\,The function $\bar{\sigma}$ satisfies the Yamada-Watanabe
condition. 

A well-known example of the process $V$ is the CIR process for which
$\bar{b}(x)=a_1-a_2x$ and 
$\bar{\sigma}(x)=a_3\sqrt{x}$. It is assumed in \cite{GGG} that the model in (\ref{E:mood}) does not have the drift term ($b=0$), while the volatility function $\sigma$ is time-homogeneous. The drift-less fractional Heston models studied in \cite{AY,CCR,GJRS} are special cases of the models described above. In fractional Heston models, the volatility is the fractional integral operator applied to the CIR process. A different generalization of the Heston model (a rough Heston model) is due to El Euch and Rosenbaum 
(see \cite{EER}). In the rough Heston model, the fractional integral operator is applied to the CIR equation, and not to the CIR process.
We do not know whether the LDP in Theorem \ref{T:nmn} holds the log-price process in the rough Heston model. 

For a small-noise parameter $\varepsilon$, we define the scaled version $V^{(\varepsilon)}$ of the process $V$ as the solution 
to the following equation: 
$
dV^{(\varepsilon)}_t=\bar{b}(V_t^{(\varepsilon)})dt+\sqrt{\varepsilon}\bar{\sigma}(V_t^{(\varepsilon)})dB_t,
$
with the initial condition given by $V_0^{(\varepsilon)}=v_0 > 0$. The scaled volatility process in the non-Gaussian fractional 
model is given by 
$
\widehat{B}^{(\varepsilon)}_t=\int_{0}^{t}K(t, s)U(V_s^{(\varepsilon)})ds,
$
while the scaled log-price process is as follows:
$$
X^{(\varepsilon)}_t=-\frac{1}{2}\varepsilon\int_0^t\sigma(\widehat{B}^{(\varepsilon)}_t)^2dt 
+\sqrt{\varepsilon}\int_0^t\sigma(\widehat{B}^{(\varepsilon)}_t)(\bar{\rho}dW_t + \rho dB_t).
$$
In addition, the map $f\mapsto\widehat{f}$ is defined by
$
\widehat{f}(t)=\int_0^tK(t,s)U(\varphi_f(s))ds
$
where $t\in[0,T]$, $f\in\mathbb{H}_0^1$, and $\varphi_f$ is the unique solution to the ODE $\dot{v}=\bar{b}(v)+\bar{\sigma}(v)\dot{f}$, $f\in\mathbb{H}_0^1[0,T]$.
\subsection{Unification: Mixed Models}\label{SS:fu}
In this section, we introduce a new class of volatility models. A model belonging to this class may be called a mixture of a multivariate Gaussian stochastic volatility model and a multivariate non-Gaussian fractional model.

Let $K_i$, with $0\le i\le d$, and $\{K_{ij}\}$, with $1\le i\le d$ and
$1\le j\le m$, be families of admissible Volterra type Hilbert-Schmidt kernels such that Assumption F holds for them. Define an 
$(d\times k)$-matrix by ${\cal K}=(K_{ij})$. Suppose that $V$ is an auxiliary $k$-dimensional continuous process defined on the space ${\cal W}^m$ equipped with the canonical set-up. Suppose also that Conditions (H1) -- (H6) in \cite{CF} are satisfied for the equation defining the process $V$ 
(see (\ref{E:2b}) and Remark \ref{R:imp}). Let $U$ be a continuous map from 
$\mathbb{R}^k$ into $\mathbb{R}^d$, and consider the following stochastic model for the volatility process: 
$
Y_t=(Y_t^{(1)},\cdots Y_t^{(d)})
$ 
where $t\in[0,T]$ and 
\begin{equation}
Y_{t}^{(i)}=x_i+\int_0^tK_i(t,s)U_i(V_s)ds+\sum_{j=1}^m\int_0^tK_{ij}(t,s)dB^{(j)}_s,\quad 1\le i\le d.
\label{E:2d}
\end{equation} 
The volatility model introduced in (\ref{E:2d}) is a special example of the models described in (\ref{E:2a}). Indeed, we can assume that the map $a$ in (\ref{E:2a}) does not depend on the fourth variable and its components are given by $a_i(t,s,u,v)=K_i(t,s)U_i(u)$ where 
$1\le i\le d$, $t,s\in[0,T]$, and $u\in\mathbb{R}^k$. 
We can also assume that the matrix function $c$ does not depend on the third and fourth variables, and its elements are defined by 
$c_{ij}(t,s,u,v)=K_{ij}(t,s)$ where
$t,s\in[0,T]$, $1\le i\le d$, and $1\le j\le m$. The scaled version of the process in (\ref{E:2d}) is given by
\begin{equation}
Y_{t}^{i,\varepsilon}=x_i+\int_0^tK_i(t,s)U_i(V_s^{(\varepsilon)})ds+\sqrt{\varepsilon}\sum_{j=1}^m\int_0^tK_{ij}(t,s)dB^{(j)}_s,
\quad 1\le i\le d.
\label{E:2f}
\end{equation} 
The model for the volatility process described in (\ref{E:2d}) is more general than the models considered in \cite{GGG} and \cite{Gul1}. 
Unlike the latter models, the model in (\ref{E:2d}) is multidimensional, the restrictions on $V$ and $U$ are weaker than those in \cite{GGG}, and  the volatility process in (\ref{E:2d}) is a mixture of volatility processes in the above-mentioned models. 
\begin{remark}\label{R:volvol}
By assuming that $U=0$ in (\ref{E:2f}), we obtain the scaled volatility process in a multivariate Gaussian stochastic volatility model. This process is defined by
\begin{equation}
\widetilde{Y}_{t}^{i,\varepsilon}=x_i+\sqrt{\varepsilon}\sum_{j=1}^m\int_0^tK_{ij}(t,s)dB^{(j)}_s,\quad 1\le i\le d.
\label{E:volvol1}
\end{equation} 
On the other hand, the scaled volatility process in a multivariate non-Gaussian fractional model can be obtained from (\ref{E:2d}) by 
setting $K_{ij}=0$ for all $1\le i\le d$ and $1\le j\le m$. This process is given by
\begin{equation}
\widehat{Y}_{t}^{i,\varepsilon}=x_i+\int_0^tK_i(t,s)U_i(V_s^{(\varepsilon)})ds,\quad 1\le i\le d.
\label{E:volvol2}
\end{equation}
\end{remark}
\begin{theorem}\label{T:cxc}
Assumptions (C1) -- (C7) hold true for the mixed volatility model introduced in (\ref{E:2d}). 
\end{theorem}

\it Proof. \rm
It is not hard to see that Assumptions (C1) and (C2) are satisfied for the model in (\ref{E:2d}). Assumption (C3) is satisfied as well. 
Note that for every $\varepsilon\in(0,1]$, we have an equality in (\ref{E:2f}), and not an equation.  
The process $Y^{(\varepsilon)}$ defined in (\ref{E:2f}) is continuous. The previous statement follows from the continuity of the function $U$ and the process $V^{(\varepsilon)}$, and from Assumption F. Next, using (\ref{E:2f}) and the fact that pathwise uniqueness and existence in the strong sense hold for the equation defining the process $V$ (see Condition (H3) in \cite{CF}), we derive Assumption (C3).

The validity of Assumption (C4) can be checked as follows. Let $f\in L^2([0,T],\mathbb{R}^m)$. Then, the equation in (\ref{E:t2}) becomes the following equality:
\begin{equation}
\eta^{(i)}_f(t)=x_i+\int_0^tK_i(t,s)U_i(\psi_f(s))ds+\sum_{j=1}^m\int_0^tK_{ij}(t,s)f_j(s)ds,\quad 1\le i\le d
\label{E:just}
\end{equation}
where $\psi_f$ solves the equation
$
\psi(s)=v_0+\int_0^s\bar{b}(r,\psi)dr+\int_0^s\bar{\sigma}(r,\psi)f(r)dr.
$
The latter equation is uniquely solvable and the solution $\psi_f$ is in the space ${\cal W}^k$. The previous statement follows from the results obtained in \cite{CF}. 
The equality in (\ref{E:just}) implies the validity of Assumption (C4). We also have $\Gamma_xf=\eta_f$. Here $\eta_f$ is given by 
(\ref{E:just}).

We will next turn our attention to Assumption (C5). Suppose that  $f_{n}\in L^2([0,T],\mathbb{R}^m)$, with $n\ge 1$, is a sequence of control functions 
such that $f_{n}\mapsto f$ weakly in $L^2([0,T],\mathbb{R}^m)$. We have from (\ref{E:just}) that
\begin{equation}
\eta_{f_n}^{(i)}(t)=x_i+\int_0^tK_i(t,s)U_i(\psi_{f_{n}}(s))ds+\sum_{j=1}^m\int_0^tK_{ij}(t,s)f_j^{(n)}(s)ds,\quad 1\le i\le d.
\label{E:justy}
\end{equation}
It will be shown next that for every $1\le i\le d$, $\eta_{f_n}^{(i)}\mapsto\eta_f^{(i)}$ in ${\cal W}^1$ as $n\rightarrow\infty$. First, observe that
since the kernels appearing in (\ref{E:justy}) are admissible, and, moreover, the Hilbert-Schmidt operators appearing in (\ref{E:just}) 
and (\ref{E:justy}) are compact maps from $L^2([0,T],\mathbb{R}^1)$ into ${\cal W}^1$, the last term on the right-hand side
of (\ref{E:justy}) tends to the last term on the right-hand-side of (\ref{E:just}) as $n\rightarrow\infty$. We can also prove that the same conclusion holds for the second terms on the right-hand sides of (\ref{E:justy}) and (\ref{E:just}) by using the fact that
$\psi_{f_{n}}\rightarrow\psi_{f}$ in ${\cal W}^k$ (this follows from the results obtained in \cite{CF}). Hence,
for all $1\le i\le d$, $U_i(\psi_{f_{n}})\mapsto U_i(\psi_{f})$ in ${\cal W}^1$ as $n\rightarrow\infty$. Therefore,
Assumption (C5) is satisfied.

The validity of Assumption (C6) follows from the following: (a)\,Every kernel $K_{ij}$ is admissible and satisfies Fernique's condition;\,(b)\,The process $s\mapsto B_s^{v}$ is standard Brownian motion with respect to the measure
$\mathbb{P}^{v}$.

Finally, we will prove that Assumption (C7) holds. Let $\varepsilon_n\rightarrow 0$ as $n\rightarrow\infty$, and let 
$v^{(n)}\in{\cal M}^2[0,T]$, with $n\ge 1$, be such that for some $N> 0$,
\begin{equation}
\sup_{n\ge 1}\int_0^T||v^{(n)}_s||_m^2ds\le N
\label{E:nok}
\end{equation}
$\mathbb{P}$-a.s. Our first goal is to prove that the sequence $n\mapsto Y_{\cdot}^{n,v^n}$ is tight in ${\cal W}^d$. It suffices to show that the sequences of components $n\mapsto Y_{\cdot}^{i,n,v^n}$, with $1\le i\le d$, are tight in ${\cal W}^1$. For every $1\le i\le d$, we have
\begin{align}
Y_{t}^{i,n,v^n}&=x_i+\int_0^tK_i(t,s)U_i(V_s^{\varepsilon_n,v^n})ds+\sum_{j=1}^m\int_0^tK_{ij}(t,s)v^{j,n}_{s}ds \nonumber \\
&\quad+\sqrt{\varepsilon_n}\sum_{j=1}^m\int_0^tK_{ij}(t,s)dB^{(j)}_s.
\label{E:ol}
\end{align}
In the rest of the proof, we will use the fact that for a finite number of tight sequences $A_n^{(k)}$, with $1\le k\le m$ and $n\ge 1$, of random elements in a normed space, the sum $\sum_{k=1}^mA_n^{(k)}$ is a tight family of random elements.

For a fixed index $1\le i\le d$, consider the sequence $Y_{\cdot}^{i,n,v^n}$, $n\ge 1$, of random elements in
${\cal W}^1$. Our goal is to prove that this sequence is tight. By taking into account what was said above, we see that it suffices to show that the following sequences of random elements in ${\cal W}^1$ are tight:
$J_{1,n}^{(i)}(t)=\int_0^tK_i(t,s)U_i(V_s^{\varepsilon_n,v^n})ds$, 
$J_{2,n}^{(i)}(t)=\sum_{j=1}^m\int_0^tK_{ij}(t,s)v^{j,n}_{s}ds$, and
$J_{3,n}^{(i)}(t)=\sqrt{\varepsilon_n}\sum_{j=1}^m\int_0^tK_{ij}(t,s)dB^{(j)}_s$.

Let us begin with the sequence $n\mapsto J_{1,n}^{(i)}$. Using the fact that the sequence of random elements 
$n\mapsto V_{\cdot}^{\varepsilon_n,v^n}$ is tight in ${\cal W}^k$ (Condition (H6) in \cite{CF}), and $U_i$ is a 
continuous map from $\mathbb{R}^k$ into $\mathbb{R}$, we see that the sequence 
$n\mapsto U_i(V_{\cdot}^{\varepsilon_n,v^n})$ is tight in ${\cal W}^1$. The latter statement follows from Prokhorov's theorem and 
Corollary 3 on p. 9 in \cite{Bi}. Therefore, for every $\varepsilon> 0$ there exists a compact set $C_1^{(i)}$ in ${\cal W}^1$ such that
$
\mathbb{P}(U_i(V_{\cdot}^{\varepsilon_n,v^n})\in C_1^{(i)})\ge 1-\varepsilon
$
for all $n\ge 1$. Denote by $C_2^{(i)}$ the image of $C_1^{(i)}$ by the Hilbert-Schmidt operator with the kernel $K_i$. Then, $C_2^{(i)}$ is a 
compact subset of ${\cal W}^1$, by the compactness of the above-mentioned operator. Therefore,
$
\mathbb{P}(J_{1,n}^{(i)}(\cdot)\in C_2^{(i)})\ge\mathbb{P}(U_i(V_{\cdot}^{\varepsilon_n,v^n})\in C_1^{(i)})\ge 1-\varepsilon.
$
The previous estimates show that the sequence $n\mapsto J_{1,n}^{(i)}$ is tight in ${\cal W}^1$.

Let us next consider the sequence $n\mapsto J_{2,n}^{(i)}$. Since we assumed that (\ref{E:nok}) holds, there exists a set 
$\widetilde{\Omega}\subset\Omega$ of full measure such that for every $1\le j\le m$, the family $v^{j,n}_{\cdot}(\omega)$, with $n\ge 1$ and
$\omega\in\widetilde{\Omega}$, is uniformly bounded in $L^2[0,T]$. Therefore, for every $1\le i\le d$, the image of this family with respect to the Hilbert-Schmidt operator with the kernel $K_{ij}$ is a compact subset of ${\cal W}^1$. Now, it easily follows that the sequence 
of random elements $n\mapsto J_{2,n}^{(i)}$ is tight in ${\cal W}^1$. 

Finally, we turn our attention to the sequence $n\mapsto J_{3,n}^{(i)}$. The sum appearing in the definition of $J_{3,n}^{(i)}$ is a continuous stochastic process on $\Omega={\cal W}^m$ with state space ${\cal W}^1$. Therefore, the sequence of random elements $n\mapsto 
J_{3,n}^{(i)}(\cdot)$ on $\Omega$ with values in ${\cal W}^1$ converges in 
${\cal W}^1$ to the identically zero random element. Let $g$ be a bounded continuous real function on ${\cal W}^1$.
Then $\mathbb{E}[g(J_{3,n}^{(i)}(\cdot))]=0$ by the bounded convergence theorem, and hence, the sequence of random elements 
$n\mapsto J_{3,n}^{(i)}$ is weakly convergent. It follows from Prokhorov's theorem that this sequence is tight in ${\cal W}^1$. 

Summarizing what was said above, we conclude that for every $1\le i\le d$, the sequence of random elements $n\mapsto Y_{\cdot}^{i,n,v^n}$ defined in (\ref{E:ol}) is tight in ${\cal W}^1$. This establishes the first part of Assumption (C7).
It remains to show that the second part of Assumption (C7) holds. This means that we must prove the estimate in (\ref{E:cond}). 
In our special case,
this estimate reduces to the following:
$$\sup_{t\in[0,T]}\sum_{i=1}^d\sum_{j=1}^m\int_0^tK_{ij}(t,s)^2ds<\infty.$$
To obtain the previous inequality, it suffices to show that 
\begin{equation}
\sup_{t\in[0,T]}\int_0^tK_{ij}(t,s)^2ds<\infty
\label{E:pre}
\end{equation}
for all $1\le i\le d$ and $1\le j\le m$. The inequality in (\ref{E:pre}) follows from the fact that
the function $t\mapsto\int_0^tK_{ij}(t,s)^2ds$, with $t\in[0,T]$, $1\le i\le d$, and $1\le j\le m$, is the variance function 
of the continuous Gaussian process $t\mapsto \int_0^tK_{ij}(t,s)dB_s^j$, and the variance function of such a process is
continuous. Therefore, Assumption (C7) is satisfied. 

This completes the proof of Theorem \ref{T:cxc}.
\begin{remark}\label{R:rrro}
Since Assumptions (C1) -- (C7) hold for the model is (\ref{E:2d}), Theorem \ref{T:beg1} can be applied to the scaled volatility process 
defined by $\widehat{B}^{(\varepsilon)}=Y^{(\varepsilon)}$ (see (\ref{E:2f})), while Theorem \ref{T:beg11} holds for the process
$(\sqrt{\varepsilon}W,\sqrt{\varepsilon}B,\widehat{B}^{(\varepsilon)})$, provided that canonical set-up is used. Moreover, the LDP in
Theorem \ref{T:nmn} is valid for the stochastic volatility model in (\ref{E:moodk}), with the process $\widehat{B}=Y$ as the volatility process, if the model is defined on the space $\Omega={\cal W}^m\times{\cal W}^m$ equipped with the $2m$-dimensional canonical set-up.
As corollaries, we obtain sample path LDPs for multivariate Gaussian models and multivariate non-Gaussian fractional models on the canonical set-up. It will be shown in the next subsection that for multivariate Gaussian models and for certain multivariate non-Gaussian fractional models, Theorems \ref{T:beg1}, \ref{T:beg11}, and \ref{T:nmn} hold on any set-up.
\end{remark}

\subsection{LDPs for Multivariate Gaussian Models and Multivariate Non-Gaussian Fractional Models}\label{SS:mGm}
Recall that the scaled volatility process in a multivariate Gaussian stochastic volatility model is given by 
$\widehat{B}^{(\varepsilon)}=\widetilde{Y}^{(\varepsilon)}$ where the latter process is defined in (\ref{E:volvol1}). In this subsection, we will prove the following assertion.
\begin{theorem}\label{T:dert}
The LDP in Theorem \ref{T:nmn} holds for a multivariate Gaussian model defined on any set-up.
\end{theorem} 

\it Proof. \rm For one-factor Gaussian models, Theorem \ref{T:dert} was established in \cite{Gul1}, Theorem 4.2. We will only sketch the proof of Theorem \ref{T:dert} for multivariate models and leave filling in the necessary details to the interested reader. The first step in the proof of Theorem \ref{T:dert} is to establish that the LDP in Theorem \ref{T:beg11} holds on any set-up. Then, we can derive 
Theorem \ref{T:dert} by observing that the proof of Theorem \ref{T:nmn} given in Section \ref{S:prov} uses only the set-up utilized in Theorem \ref{T:beg11}. 

Let us take $U=0$ in (\ref{E:2d}). Then, the expression in (\ref{E:2d}) represents the volatility vector $\widehat{B}$ in the multivariate Gaussian model. Consider the random vector ${\cal X}=(W,B,\widehat{B})$ on $\Omega$ with values in the space 
$\Lambda={\cal W}^m\times{\cal W}^m\times{\cal W}^d$. Reasoning as in the proof of Theorem 6.8 in \cite{Gul1}, we can show that 
${\cal X}$ is a centered Gaussian random vector. Actually, ${\cal X}$ is a Gaussian random vector in a smaller space
$\widetilde{{\cal G}}$. This will be explained next (see also the proof on p. 63 of \cite{Gul1}). Let
$H=L^2([0,T],\mathbb{R}^m)\times L^2([0,T],\mathbb{R}^m)$. Then $H$ is a separable Hilbert space equipped with the norm 
$$
||(h_1,h_2)||_H=\sqrt{||h_1||_{L^2([0,T],\mathbb{R}^m)}^2+||h_2||_{L^2([0,T],\mathbb{R}^m)}^2}.
$$ 
Let $j:H\mapsto\Lambda$ be the map defined by
$j(h_0,h_1)=(g_0,g_1,g_2)$ where for every $t\in[0,T]$, $g_0(t)=\int_0^th_0(s)ds$, $g_1(t)=\int_0^th_1(s)ds$, and 
$g_2(t)=\int_0^tK(t,s)h_2(s)ds$. The map $j$ is an injection. Set $\widetilde{{\cal G}}=\overline{j(H)}$ where the closure is taken in the space $\Lambda$. Then, $\widetilde{{\cal G}}$ is a separable Banach space. It was established in \cite{Gul1} for one-dimensional Gaussian models that ${\cal X}$ is a Gaussian random vector in the space $\widetilde{{\cal G}}$ (see Theorem 6.11 in\cite{Gul1}). It is not hard to see that the same result holds in the multivariate case. We can also find the covariance operator $\widehat{K}$ by imitating the proof of Theorem 6.14 in \cite{Gul1}. Let $\zeta$ be the distribution of the random vector ${\cal X}$ on the measurable space
$(\widetilde{{\cal G}},{\cal B}(\widetilde{{\cal G}}))$. Then, $\zeta$ is a Gaussian measure. Using the same ideas as in the proof of 
Theorem 6.19 in \cite{Gul1}, we can prove that the quadruple $(\widetilde{{\cal G}},H,j,\zeta)$ is an abstract Wiener space. Finally,
applying the known LDP for abstract Wiener spaces (see the references before Theorem 6.19 in \cite{Gul1}), we establish the validity of 
the LDP in Theorem \ref{T:beg11} on any set-up.

This completes the proof of Theorem \ref{T:dert}.

We will next turn our attention to multivariate non-Gaussian fractional models. The volatility process in such a model is given by
$\widehat{B}^{(\varepsilon)}=\widehat{Y}^{(\varepsilon)}$ where the latter process is defined in (\ref{E:volvol1}). 
Recall that the process $V^{(\varepsilon)}$ appearing in (\ref{E:volvol1}) is as follows:
\begin{equation}
V_t^{(\varepsilon)}=y+\int_0^ta(s,V^{(\varepsilon)})ds+\sqrt{\varepsilon}\int_0^tc(s,V^{(\varepsilon)})dB_s.
\label{E:222a}
\end{equation} 
where $a$ is a map from the space $[0,T]\times{\cal W}^d$ into the space $\mathbb{R}^d$, while $c$ is a map 
from the space $[0,T]\times{\cal W}^d$ into the space of $(d\times m)$-matrices. It will be assumed in the rest of this subsection that the maps $a$ and $c$ are locally
Lipschitz and satisfy the sub-linear growth condition. 
\begin{theorem}\label{T:riur}
Under the conditions formulated above, the LDP in Theorem \ref{T:nmn} holds for a multivariate Gaussian model defined on any set-up.
\end{theorem}

\it Proof. \rm It was shown in Corollary \ref{C:rre} that if the maps $a$ and $c$ are locally
Lipschitz and satisfy the sub-linear growth condition, then the LDP in Theorem \ref{T:beg11}  
holds for the process $\varepsilon\mapsto(\sqrt{\varepsilon}W,\sqrt{\varepsilon}B,V^{(\varepsilon)}_{\cdot})$ defined on any set-up.
Therefore, the LDP in Theorem \ref{T:nmn} holds true on any set-up (see Remark \ref{R:csu}).

The proof of Theorem \ref{T:riur} is thus completed.

\subsection{Unification: Models with Reflection}\label{SS:mvr}
Large deviation principles for log-price processes and volatility processes in one-factor stochastic volatility models with reflection were obtained in \cite{Gul2}. 
We will next establish similar LDPs for multivariate models with reflection.

Let $O$ be an open set in $\mathbb{R}^d$, with the boundary $\partial O$ and the closure $\overline{O}$. Suppose also that a vector field $K$ of reflecting directions is given on $\partial O$. The main restrictions that we impose on the model with reflection are as follows:
(i)\,\,The unique solvability of Skorokhod's problem for $(O,K,f)$ where $f\in{\cal W}^d$ and $f(0)\in\overline{O}$. (ii)\,\,The continuity 
of the Skorokhod map $\Gamma$ in the space ${\cal W}^d$. All the necessary definitions can be found in Chapter 2 of \cite{P}. 
In a multivariate stochastic volatility model with reflection, the scaled volatility process is a reflecting diffusion in $\overline{O}$ given by $\widehat{B}_t^{(\varepsilon)}=(\Gamma U^{(\varepsilon)})(t)$, $t\in[0,T]$
where the process $t\mapsto U^{\varepsilon}_t$ is the unique continuous solution to the following $d$-dimensional stochastic differential equation:
\begin{equation}
dU_t^{(\varepsilon)}=\hat{a}(t,(\Gamma U^{(\varepsilon)})(t))dt+\sqrt{\varepsilon}
\hat{c}(t,(\Gamma U^{(\varepsilon)})(t))dB_t,\quad U_0^{\varepsilon}=y\in\overline{O}
\label{E:fse}
\end{equation}
(see \cite{P}, (2.5) and (2.6) in Section 2). In (\ref{E:fse}), $\hat{a}$ is a map from 
$[0,T]\times\mathbb{R}^d$ into $\mathbb{R}^d$, while $\hat{c}$ maps $[0,T]\times\mathbb{R}^d$
into the space of $d\times m$-matrices. We assume that the maps $a(t,\varphi)=\hat{a}(t,(\Gamma\varphi)(t)$ and $c(t,\varphi)
=\hat{a}(t,(\Gamma\varphi)(t)$, where $t\in[0,T]$ and $\varphi\in{\cal W}^d$, are locally Lipschitz continuous in the second variable, uniformly in time, and satisfy the sublinear growth condition in the second variable, uniformly in time 
(see Definitions (A1) and (A2) in Section 3 of \cite{CF}). Then, Conditions (H1) - (H6) in \cite{CF}
are satisfied (see Section 3 of \cite{CF}), and hence Theorems \ref{T:nmn} and \ref{T:beg1} hold true for the log-price process and the volatility process in a stochastic volatility model with reflection, under the restrictions mentioned above. More information about multidimensional reflecting diffusions and also examples of uniquely solvable Skorokhod's problems can be found in Chapters 2 and 3 of \cite{P}. Note that Skorokhod's problem for the half-line is uniquely solvable, and the Skorokhod map $\Gamma:{\cal W}^1\mapsto{\cal W}^1$ is defined by
$(\Gamma f)(t)=f(t)-\min_{s\in[0,t]}(f(s)\wedge 0)$, $t\in[0,T]$
(see \cite{P} for more details). The Skorokhod map $\Gamma$ is continuous from the space ${\cal W}^1$ into itself. It is also 
$\widetilde{{\cal F}}^1_t
/{\cal F}^1_t$-measurable for every $t\in[0,T]$. The previous statements follow from
Lemma 1.1.1 in \cite{P}. 

A special example of a one-factor stochastic volatility model with reflection is one of the three versions of the Stein and Stein model 
(see (\cite{Gul2})). The volatility process in this model is the instantaneously reflecting Ornstein-Uhlenbeck process.
\begin{remark}\label{R:import}
In a model with reflection, the map $G$, appearing in Definitions \ref{D:vo} and \ref{D:dds}, is the Skorokhod map $\Gamma$. In the other models considered in the present paper, the map $G$ is the identity map.
\end{remark}

\subsection{Unification: Volterra Type SDEs}\label{SS:Wang}
Consider the following multidimensional Volterra type stochastic differential equation:
\begin{equation}
Y_t=y+\int_0^ta(t,s,Y_s)ds+\int_0^tc(t,s,Y_s)dB_s.
\label{E:oo}
\end{equation}
The equation in (\ref{E:oo}) is a special case of the equation in (\ref{E:2a}). We will also use a scaled version of the equation in 
(\ref{E:oo}), that is, the equation
\begin{equation}
Y_t^{\varepsilon}=y+\int_0^ta(t,s,Y_s^{\varepsilon})ds+\sqrt{\varepsilon}\int_0^tc(t,s,Y_s^{\varepsilon})dB_s,
\label{E:sca}
\end{equation}
and a scaled controlled version given by
\begin{equation}
Y_t^{\varepsilon,v}=y+\int_0^ta(t,s,Y_s^{\varepsilon,v})ds+\int_0^tc(t,s,Y_s^{\varepsilon,v})v_sds
+\sqrt{\varepsilon}\int_0^tc(t,s,Y_s^{\varepsilon,v})dB_s.
\label{E:scac}
\end{equation}

In the papers \cite{W} of Wang and \cite{Z} of Zhang, the following conditions were formulated:
\\
\\
(H1)\,For some $p> 2$ there exists $C_T> 0$ such that for all $x,y\in\mathbb{R}^d$ and $s,t\in[0,T]$,
\begin{equation}
||a(t,s,x)-a(t,s,y)||_d\le C_TK_1(t,s)\rho^{\frac{1}{p}}(||x-y||_d^p),
\label{E:1}
\end{equation}
\begin{equation}
||c(t,s,x)-c(t,s,y)||_{d\times m}^2\le C_TK_2(t,s)\rho^{\frac{2}{p}}(||x-y||_d^p),
\label{E:2}
\end{equation}
and
\begin{equation}
\int_0^t(||a(t,s,0)||_d+||c(t,s,0)||^2_{d\times m})ds\le C_T
\label{E:3}
\end{equation}
where $K_i$, with $i=1,2$, are two positive functions on $[0,T]^2$ for which
\begin{equation} 
\int_0^t\left[K_1(t,s)^{\frac{p}{p-1}}+K_2(t,s)^{\frac{p}{p-2}}\right]ds\le C_T,\quad t\in[0,T].
\label{E:45}
\end{equation}
In addition, $\rho:\mathbb{R}^{+}\mapsto\mathbb{R}^{+}$ is a concave function satisfying
\begin{equation}
\int_{0+}^{\cdot}\rho(u)^{-1}du=\infty.
\label{E:5}
\end{equation}
\vspace{0.1in}
(H2)\,For all $t,t^{\prime},s\in[0,T]$ and $x\in\mathbb{R}^d$,
\begin{equation}
||a(t,s,x)-a(t^{\prime},s,x)||_d\le F_1(t^{\prime},t,s)(1+||x||_d),
\label{E:6}
\end{equation}
\begin{equation}
||c(t,s,x)-c(t^{\prime},s,x)||_{d\times m}^2\le F_2(t^{\prime},t,s)(1+||x||_d^2),
\label{E:7}
\end{equation}
and for some $C> 0$ and $\theta> 1$,
\begin{equation}
\int_0^t(||a(t,s,0)||_d^{\theta}+||c(t,s,0)||_{d\times m}^{2\theta})ds< C.
\label{E:8}
\end{equation}
The functions $F_i$, $i=1,2$, in (\ref{E:6}) and (\ref{E:7}) are two positive functions on $[0,T]^3$ satisfying 
the condition
\begin{equation}
\int_0^{t\wedge t^{\prime}}(F_1(t^{\prime},t,s)+F_2(t^{\prime},t,s))ds\le C|t-t^{\prime}|^{\gamma}
\label{E:9}
\end{equation}
for some $\gamma> 0$. 
\begin{remark}\label{R:mch}
It was shown in \cite{W} that if Condition (H1) holds, then there exists a unique progressively measurable solution $Y$ to the equation in
(\ref{E:oo}).
Moreover, it was established in \cite{W} that if Conditions (H1) and (H2) hold, then the unique solution $Y$ to the equation in
(\ref{E:oo}) has a $\delta$-H\"{o}lder continuous version for any
$\delta\in(0,\frac{1}{p}\wedge\frac{\theta-1}{2\theta}\wedge\frac{\gamma}{2})$. 
\end{remark}

In \cite{Z}, Zhang obtained a sample path LDP for the unique solution 
\begin{equation}
\varepsilon\mapsto Y_{\cdot}^{(\varepsilon)}(\cdot),\quad\varepsilon\in(0,1]
\label{E:pro}
\end{equation} 
to the equation in (\ref{E:sca}), under Conditions (H1) and (H2) and two extra conditions (H3) and (H4) (see Theorem 1.2 in \cite{Z}).  
Note that the initial condition 
$y\in\mathbb{R}^d$ plays the role of a variable in the process defined by (\ref{E:pro}). The state space of the process in (\ref{E:pro}) 
is the space of continuous maps from $[0,T]\times\mathbb{R}^d$
into $\mathbb{R}^d$.  Using Theorem 1.2 established in \cite{Z} and the contraction principle, we can obtain a sample path 
LDP for the process 
\begin{equation}
\varepsilon\mapsto Y_{\cdot}^{(\varepsilon)},\quad\varepsilon\in(0,1],
\label{E:eres}
\end{equation}
with the initial condition $y\in\mathbb{R}^d$ that is fixed. The state space of the process in (\ref{E:eres}) is the space ${\cal W}^d$.

It will be shown below that Conditions (H3) and (H4) are not needed in the LDP for the process in (\ref{E:eres}). 
\begin{remark}\label{R:Zhang}
In the remaining part of the present paper, we will employ a weaker condition than Condition (H2). In the new condition, 
the restriction $t,t^{\prime},s\in[0,T]$ in (\ref{E:6}) 
and (\ref{E:7}) is replaced by
the restriction $0\le s\le t,t^{\prime}\le T$. We denote the new condition by $(\widehat{H}2)$. By analyzing the main results obtained 
in the paper \cite{W} of Wang, one can see that these results hold true with Condition (H2) replaced by Condition $(\widehat{H}2)$. 
\end{remark}

We will next show that Conditions (H1) and $(\widehat{\rm H}2)$ imply Assumptions (C1) - (C7). 
\begin{lemma}\label{L:WG} 
Suppose Conditions (H1) and $(\widehat{H}2)$ hold true for the maps $a$ and $c$ appearing in (\ref{E:oo}). 
Then, Assumptions (C1) - (C7) introduced in Section \ref{S:VP} are satisfied.
\end{lemma}

The following corollary follows from Lemma \ref{L:WG}.
\begin{corollary}\label{C:Cor}
Under Conditions (H1) and $(\widehat{H}2)$, the LDPs in Theorems \ref{T:beg1} and \ref{T:beg11} hold for the process 
in (\ref{E:eres}).
\end{corollary}

\it Proof of Lemma \ref{L:WG}. \rm Some of the ideas utilized in the proof of Lemma \ref{L:WG} are borrowed from \cite{W} and \cite{Z}. 
Suppose Conditions (H1) and $(\widehat{\rm H}2)$ are satisfied.
Then, it is clear that Assumption (C1) holds. For the equation in (\ref{E:oo}), the inequalities in (C2)(a) are as follows:
$$
\int_0^t||a(t,s,\varphi(s)||_dds<\infty\quad\mbox{and}\quad \int_0^t||c(t,s,\varphi(s)||_{d\times m}^2ds<\infty
$$
for all $t\in[0,T]$ and $\varphi\in{\cal W}^d$. These inequalities can be obtained by integrating the estimates
$$
||a(t,s,\varphi(s)||_d\le||a(t,s,0||_d+CK_1(t,s)(1+||\varphi(s)||_d)
$$
and 
$$
||c(t,s,\varphi(s)||_{d\times m}^2\le
2||c(t,s,0||_{d\times m}^2+CK_2(t,s)(1+||\varphi(s)||_d^2)
$$
over the interval $[0,T]$,
and then using (\ref{E:3}) and (\ref{E:45}). The previous estimates follow from the inequalities established in \cite{W} 
(see the bottom part of p. 1064 in \cite{W}). 

We will next prove the first statement in Assumption (C2)(b). Suppose $t^{\prime}\le t$. The case where $t< t^{\prime}$ is similar.
We have
\begin{align}
&||\int_0^ta(t,s,\varphi(s))ds-\int_0^{t^{\prime}}a(t^{\prime},s,\varphi(s))ds||_d \nonumber \\
&\le\int_{t^{\prime}}^t||a(t,s,\varphi(s))-a(t,s,0)||_dds+\int_{t^{\prime}}^t||a(t,s,0)||_dds \nonumber \\
&+\int_0^{t^{\prime}}||a(t,s,\varphi(s))-a(t^{\prime},s,\varphi(s))||_dds.
\label{E:labb}
\end{align}
Suppose $t^{\prime}\rightarrow t$. Then, the three terms on the right-hand side of (\ref{E:labb}) tend to zero. For the first term, 
this follows from (\ref{E:1}), for the second one, we can use (\ref{E:8}) and H\"{o}lder's inequality, while for the third one, 
one can use the estimates in (\ref{E:6}) and (\ref{E:9}), with $0\le s\le t^{\prime}\le t\le T$. 
Note that in the previous reasoning we used Condition 
$(\widehat{\rm H}2)$. This completes the proof of the first statement in Assumption (C2)(b).

The second condition in (C2)(b) and the condition in (C2)(c)
follow from (\ref{E:1}), (\ref{E:2}), (\ref{E:45}), and the concavity of the function $\rho$. This shows that Assumption (C2) holds true.

It was established in \cite{W} that under the restrictions in (H1) and (H2) the equation in (\ref{E:oo}) is strongly uniquely solvable
(see Theorems 1.1 and 1.3 in \cite{W}). We have already mentioned that it is possible to replace (H2) by $(\widehat{\rm H}2)$ in the previous statement (see Remark \ref{R:Zhang}). Similarly, for every $\varepsilon\in(0,1]$, the equation in (\ref{E:sca}) is strongly uniquely solvable.
This proves the validity of Assumption (C3)(a). The unique solvability of the equation in (\ref{E:scac}) for any $v\in M^2_N[0,T]$ 
was mentioned on p. 2240 in \cite{Z}. This fact can be established by using Picard's iterative method. It follows that
the statement in Assumption (C3)(b) holds true.

We will next show that Assumption (C4) is satisfied. Let $\eta_1,\eta_2\in{\cal W}^d$, $f\in L^2([0,T],\mathbb{R}^m)$, 
and suppose that for all $t\in[0,T]$, 
$$
\eta_1(t)=y+\int_0^ta(t,s,\eta_1(s))ds+\int_0^tc(t,s,\eta_1(s))f(s)ds
$$ 
and
$$
\eta_2(t)=y+\int_0^ta(t,s,\eta_2(s))ds+\int_0^tc(t,s,\eta_2(s))f(s)ds.
$$
We must prove that the equality $\eta_1(t)=\eta_2(t)$ holds for all $t\in[0,T]$ (see Remark \ref{R:tyu}). We have
\begin{align}
&||\eta_1(t)-\eta_2(t)||_d^p\le\int_0^t||a(t,s,\eta_1(s))-a(t,s,\eta_2(s))||_dds \nonumber \\
&\quad+\int_0^t||c(t,s,\eta_1(s))-c(t,s,\eta_2(s))||_{d\times m}||f(s)||_mds \nonumber \\
&\le C_T\int_0^tK_1(t,s)\rho^{\frac{1}{p}}(||\eta_1(t)-\eta_2(t)||_d^p)ds \nonumber \\
&\quad+\sqrt{C_T}\int_0^tK_2(t,s)^{\frac{1}{2}}\rho^{\frac{1}{p}}(||\eta_1(t)-\eta_2(t)||_d^p)||f(s)||_mds.
\label{E:alis}
\end{align} 
It is not hard to see using H\"{o}lder's inequality that (\ref{E:alis}) implies the following:
\begin{align*}
&||\eta_1(t)-\eta_2(t)||_d^p\le C_1\left\{\int_0^tK_1(t,s)^{\frac{p}{p-1}}ds\right\}^{p-1}
\int_0^t\rho(||\eta_1(s)-\eta_2(s)||_d^p) \\
&\quad+C_1\left\{\int_0^tK_2(t,s)\rho^{\frac{2}{p}}(||\eta_1(s)-\eta_2(s)||_d^p)ds\right\}^{\frac{p}{2}}
\left\{\int_0^t||f(s)||_m^2ds\right\}^{\frac{p}{2}} \\
&\le C_2\int_0^t\rho(||\eta_1(s)-\eta_2(s)||_d^p)ds+C_2\left\{\int_0^tK_2(t,s)^{\frac{p}{p-2}}ds\right\}^{\frac{p-2}{2}}
\int_0^t\rho(||\eta_1(s)-\eta_2(s)||_d^p)ds  \\
&\le C_3\int_0^t\rho(||\eta_1(s)-\eta_2(s)||_d^p)ds.
\end{align*}
Finally, applying Bihari's inequality (see \cite{Bih}, see also Lemma 2.1 in \cite{W}) we obtain 
$$
||\eta_1(t)-\eta_2(t)||_d^p=0\quad\mbox{for all}\,\,t\in[0,T]. 
$$
This completes the proof of (C4).

Our next goal is to prove that Assumption (C5) holds. Suppose $f_n\rightarrow f$ weakly in $L^2([0,T],\mathbb{R}^m)$. Set $\eta_n=
\Gamma_yf_n$ and $\eta=\Gamma_yf$. Then, we have $\sup_{n\ge 1}||f_n||_2<\infty$. Moreover,
\begin{align*}
\eta_n(t)-\eta(t)&=\int_0^t[a(t,s,\eta_n(s))-a(t,s,\eta(s))]ds+\int_0^t[c(t,s,\eta_n(s))-c(t,s,\eta(s))]f_n(s)ds \\
&\quad+\int_0^tc(t,s,\eta(s))[f_n(s)-f(s)]ds.
\end{align*}
Using the estimates obtained in the proof of (C4), we see that
\begin{align}
&||\eta_n(t)-\eta(t)||_d^p\le C\int_0^t\rho(||\eta_n(s)-\eta(s)||_d^p)ds+A_n
\label{E:Bi}
\end{align}
where $C> 0$ is a constant independent of $t$ and $n$, and 
\begin{equation}
A_n=C\sup_{t\in[0,T]}||\int_0^tc(t,s,\eta(s))[f_n(s)-f(s)]ds||_d^p,\quad n\ge 1.
\label{E:ee}
\end{equation}

We will next prove that $\sup_{n\ge 1}A_n<\infty$. Indeed, the following estimates hold: 
\begin{align*}
A_n&\le C_1(\sup_{t\in[0,T]}\int_0^t||c(t,s,\eta(s))||_{d\times m}^2ds)^{\frac{p}{2}}
\le C_2(\sup_{t\in[0,T]}\int_0^t||c(t,s,\eta(s))-c(t,s,0)||_{d\times m}^2ds)^{\frac{p}{2}} \\
&\quad+C_2(\sup_{t\in[0,T]}\int_0^t||c(t,s,0)||_{d\times m}^2ds)^{\frac{p}{2}}
\le C_3+C_3(\sup_{t\in[0,T]}\int_0^tK_2(t,s)\rho^{\frac{2}{p}}(||\eta(s)||_d^p)ds)^{\frac{p}{2}} \\
&\le C_3+C_4\int_0^T\rho(||\eta(s)||_d^p)ds\le C_5+C_5\int_0^T||\eta(s)||_d^pds.
\end{align*}
In the previous estimates, we used (\ref{E:2}), (\ref{E:3}), (\ref{E:45}), the concavity of the function $\rho$, and the fact that
$||f||_2+\sup_{n\ge 1}||f_n||_2<\infty$. It follows that $\sup_{n\ge 1}A_n<\infty$.

Using (\ref{E:Bi}) and Bihari's inequality we obtain
$
||\eta_n(t)-\eta(t)||_d^p\le G^{-1}(G(A_n)+Ct)
$ 
for all $n\ge 1$ and $t\in[0,T]$, where $C$ is the constant in (\ref{E:Bi}), and 
$G(x)=\int_{x_0}^x\frac{dz}{\rho(z)}$ for $x\ge 0$ and $x_0> 0$.
Hence
\begin{equation}
\sup_{t\in[0,T]}||\eta_n(t)-\eta(t)||_d^p\le G^{-1}(G(A_n)+CT),\quad n\ge 1.
\label{E:rdr}
\end{equation}

Let us consider the following family ${\cal A}$ of maps: $s\mapsto c(t,s,\eta(s))\mathbb{1}_{\{s\le t\}}$, $t\in[0,T]$. The set 
${\cal A}$ is a pre-compact subset of the space $L^2([0,T],\mathbb{R}^{d\times m})$. Indeed, we can establish the previous 
statement by constructing an $\varepsilon$-net for 
the set ${\cal A}$ for every $\varepsilon> 0$. This can be done by using estimates for the map $c$ similar to those in (\ref{E:labb}), and reasoning as in the proof of the validity of the first statement in Assumption (C2)(b). We take into account the conditions in (\ref{E:2}), (\ref{E:45}), 
(\ref{E:7}), (\ref{E:8}), and (\ref{E:9}) in the proof of the existence of the $\varepsilon$-net.
Note that the estimates in (\ref{E:7}) and (\ref{E:9}) are employed in the previous proof only in the case where 
$0\le s\le t,t^{\prime}\le T$. This means that it suffices to use Condition $(\widehat{\rm H}2)$ instead of Condition (H2) in this part of the proof of Lemma \ref{L:WG}. It follows from the previous reasoning and the weak convergence of $f_n$ to $f$ that $A_n\rightarrow 0$ where $A_n$ is defined 
in (\ref{E:ee}). Next, using (\ref{E:5}) we see that $G(A_n)\rightarrow-\infty$ as $n\rightarrow\infty$. Therefore, 
$G(A_n)+CT\rightarrow-\infty$, and hence 
\begin{equation}
G^{-1}(G(A_n)+CT)\rightarrow 0,\quad\mbox{as}\quad n\rightarrow\infty.
\label{E:perv}
\end{equation} 
Finally, using (\ref{E:rdr}) and (\ref{E:perv}) we see that Assumption (C5) holds true.

To prove that Assumption (C6) is satisfied, we first observe that under Condition (H1), the following inequality holds for the solution $Y$ to the equation in (\ref{E:oo}):
$$
\sup_{t\in[0,T]}\mathbb{E}\left[||Y_t||_d^q\right]<\infty
$$
for all $q\ge 2$. This was shown in \cite{W}, Lemma 2.2. We also have $Y=h(B)$ where $h$ is a measurable map from ${\cal W}^m$ into 
${\cal W}^d$ (see Remark \ref{R:ghg}). Therefore, the following inequality holds true:
$$
\sup_{t\in[0,T]}\mathbb{E}\left[||h(B)(t)||_d^q\right]<\infty.
$$

Let $v$ be a control satisfying (\ref{E:cve}). Recall that the control $v$ defines a new Brownian motion $B^v$ with respect to the measure $\mathbb{P}^v$ that is equivalent to the measure $\mathbb{P}$ (see (\ref{E:brew})). Our next goal is to show that we can replace $B$ by $B^v$ and $\mathbb{P}$ by $\mathbb{P}^v$ in the previous inequality. Indeed, it is not hard to see that for every 
$t\in[0,T]$, the CDF of 
the random variable $||h(B)(t)||_d$  with respect to the measure $\mathbb{P}$ is equal to that of the random variable $||h(B^v)(t)||_d$ with respect to the measure $\mathbb{P}^v$. It follows from the previous remark that
\begin{equation}
\sup_{t\in[0,T]}\mathbb{E}_{\mathbb{P}^v}\left[||h(B^v)(t)||_d^q\right]<\infty
\label{E:fgf}
\end{equation}
for all $q\ge 2$. 
It was established in the proof of Theorem 1.3 on p. 1067 of \cite{W} that if (H1) and (H2) hold (one can use $(\widehat{\rm H}2)$ 
instead of (H2)), then $t\mapsto\int_0^tc(t,s,Y_s)dB_s$ is a continuous stochastic process. Hence the process 
$t\mapsto\int_0^tc(t,s,h(B)(s))dB_s$ is continuous. Actually, the same proof works if we replace the Brownian motion 
$B$ and the measure $\mathbb{P}$ in the previous stochastic integral by the Brownian motion $B^v$ 
and the measure $\mathbb{P}^v$. Here we use the inequality in (\ref{E:fgf}). 
This establishes the validity of Assumption (C6).

Our final goal is to prove that Assumption (C7) holds. Since Assumptions (C1) - (C6) hold true, the equation in (\ref{E:scac}) is strongly uniquely solvable for every $\varepsilon\in(0,1]$ and every control $v$ satisfying the condition in (\ref{E:cve})
(see Lemma \ref{L:uni} and Remark \ref{R:rr}). We will next show that if $\varepsilon_n\rightarrow 0$, and a sequence of controls
$v^n$ is such that
\begin{equation}
\sup_{n\ge 1}\int_0^T||v^n_s||_m^2ds\le N
\label{E:trt}
\end{equation}
$\mathbb{P}$-a.s., then the following inequality holds: 
\begin{equation}
\sup_{n\ge 1}\int_0^t\mathbb{E}\left[||c(t,s,Y^{\varepsilon_n,v^n}_s)||_{d\times m}^2\right]ds<\infty,\quad t\in[0,T].
\label{E:condic}
\end{equation}
 
In the next estimates, the constants $\tau> 0$ may change from line to line. It follows from (\ref{E:scac}), (\ref{E:trt}), and (\ref{E:3}) that for every $q> p$ and $t\in[0,T]$, we have
\begin{align}
&\mathbb{E}[||Y^{\varepsilon_n,v^n}_t||_d^q]\le\tau+\tau\left\{\int_0^t||a(t,s,0)||_dds\right\}^q
+\tau\mathbb{E}\left[\left\{\int_0^t||a(t,s,Y^{\varepsilon_n,v^n}_t)-a(t,s,0)||_dds\right\}^q\right] \nonumber \\
&+\tau\mathbb{E}\left[\left\{\int_0^t||c(t,s,Y^{\varepsilon_n,v^n}_t)
-c(t,s,0)||_{d\times m}^2ds\right\}^{\frac{q}{2}}\right]+\tau\mathbb{E}\left[\left\{\int_0^t||c(t,s,0)||_{d\times m}^2
ds\right\}^{\frac{q}{2}}\right] \nonumber \\
&+\tau\mathbb{E}\left[\left\{||\int_0^tc(t,s,Y_s^{\varepsilon,v})dB_s||_d\right\}^q\right]
\le \tau+\tau\mathbb{E}\left[\left\{\int_0^t||a(t,s,Y^{\varepsilon_n,v^n}_t)-a(t,s,0)||_dds\right\}^q\right] \nonumber \\
&+\tau\mathbb{E}\left[\left\{\int_0^t||c(t,s,Y^{\varepsilon_n,v^n}_t)
-c(t,s,0)||_{d\times m}^2ds\right\}^{\frac{q}{2}}\right] 
+\tau\mathbb{E}\left[\left\{||\int_0^tc(t,s,Y_s^{\varepsilon_n,v^n})dB_s||_d\right\}^q\right].
\label{E:bdg}
\end{align}
In order to estimate the last term in (\ref{E:bdg}) we use the Burkolder-Davis-Gundy inequality.
This gives
\begin{align}
&\mathbb{E}\left[\left\{||\int_0^tc(t,s,Y_s^{\varepsilon,v^n})dB_s||_d\right\}^q\right]\le 
\tau\mathbb{E}\left[\left\{\int_0^t||c(t,s,Y^{\varepsilon_n,v^n}_t||_{d\times m}^2ds\right\}^{\frac{q}{2}}\right]\nonumber \\
&\le \tau\mathbb{E}\left[\left\{\int_0^t||c(t,s,Y^{\varepsilon_n,v^n}_t)
-c(t,s,0)||_{d\times m}^2ds\right\}^{\frac{q}{2}}\right]+\tau\mathbb{E}\left[\left\{\int_0^t||c(t,s,0)||_{d\times m}^2ds\right\}^{\frac{q}{2}}\right].
\label{E:llo} 
\end{align}
By taking into account (\ref{E:bdg}), (\ref{E:llo}), H\"{o}lder's inequality, (H1), and reasoning as in the proof of Lemma 2.2 in \cite{W},
we obtain
$
\mathbb{E}[||Y^{\varepsilon_n,v^n}_t||_d^q]\le \tau+\tau\int_0^t\mathbb{E}[||Y^{\varepsilon_n,v^n}_s||_d^q]ds.
$
Now, using Gr\"{o}nwall's inequality, we establish the following estimate:
\begin{equation}
\sup_{n\ge 1}\sup_{t\in[0,T]}\mathbb{E}[||Y^{\varepsilon_n,v^n}_t||_d^q]<\infty,\quad q\ge 2.
\label{E:hoho}
\end{equation}

Our next goal is to prove (\ref{E:condic}). Using (H1), H\"{o}lder's inequality, and the inequality $\rho^{\frac{1}{p}}(u^p)\le C(1+u)$
(see \cite{W}, p. 1064), we obtain
\begin{align}
&\int_0^t\mathbb{E}\left[||c(t,s,Y^{\varepsilon_n,v^n}_s)||_{d\times m}^2\right]ds
\le 2\int_0^t\mathbb{E}\left[||c(t,s,Y^{\varepsilon_n,v^n}_s)-c(t,s,0)||_{d\times m}^2\right]ds \nonumber \\
&\quad+2\int_0^t\mathbb{E}\left[||c(t,s,0)||_{d\times m}^2\right]ds\le \tau+
\tau\left\{\int_0^t(\mathbb{E}[||Y^{\varepsilon_n,v^n}_s||_d^2])
^{\frac{p}{2}}ds\right\}^{\frac{2}{p}} \nonumber \\
&\le \tau+\tau\left\{\int_0^t\mathbb{E}[||Y^{\varepsilon_n,v^n}_s||_d^p]ds\right\}^{\frac{2}{p}}.
\label{E:al}
\end{align}
Now, it is clear that (\ref{E:hoho}) and (\ref{E:al}) imply (\ref{E:condic}).
It remains to prove that the family $Y^{\varepsilon_n,v^n}$ is tight. We will first prove that for all $t^{\prime},t\in[0,T]$ there exist constants $\alpha> 0$, $\beta> 0$, and $\gamma> 0$ such that
$$
\sup_{n\ge 1}\mathbb{E}[||Y^{\varepsilon_n,v^n}_{t^{\prime}}-Y^{\varepsilon_n,v^n}_t||_d^{\alpha}]\le\beta|t^{\prime}-t|^{1+\beta}.
$$
We have
\begin{align*}
Y^{\varepsilon_n,v^n}_{t^{\prime}}-Y^{\varepsilon_n,v^n}_t&=\left[\int_0^{t^{\prime}}a(t^{\prime},s,Y_s^{\varepsilon_n,v^n})ds
-\int_0^ta(t,s,Y_s^{\varepsilon_n,v^n})ds\right] \\
&\quad+\left[\int_0^{t^{\prime}}c(t^{\prime},s,Y_s^{\varepsilon_n,v^n})v^n_sds-\int_0^tc(t,s,Y_s^{\varepsilon_n,v^n})v^n_sds\right] \\
&\quad+\sqrt{\varepsilon_n}\left[\int_0^{t^{\prime}}c(t^{\prime},s,Y_s^{\varepsilon_n,v^n})dB_s-\int_0^tc(t,s,Y_s^{\varepsilon_n,v^n})dB_s\right]
\\
&=J_1^{(n)}(t^{\prime},t)+J_2^{(n)}(t^{\prime},t)+J_3^{(n)}(t^{\prime},t).
\end{align*}
To estimate $J_3^{(n)}(t^{\prime},t)$, we use (\ref{E:hoho}) and reason as in the proof of Theorem 1.3 in \cite{W}. This gives the following:
\begin{equation}
\sup_{n\ge 1}\mathbb{E}[||J_3^{(n)}(t^{\prime},t)||_d^q]\le C|t^{\prime}-t|^{\delta}
\label{E:ses}
\end{equation}
for $q>\max(p,\frac{2}{\gamma})$ and $\delta=\min\{\frac{1}{p},\frac{\theta-1}{2\theta},\frac{\gamma}{2}\}$. The constant $C$ in (\ref{E:ses}) does not depend on $n$. The last statement follows from (\ref{E:hoho}). The inequality similar to that in (\ref{E:ses}) can also be established
for $J_1^{(n)}(t^{\prime},t)$ and $J_2^{(n)}(t^{\prime},t)$. Therefore, the conditions in Kolmogorov's tightness criterion are satisfied for the sequence $Y^{\varepsilon_n,v^n}$, $n\ge 1$. It follows that Assumption (C7) holds true. 

The proof of Lemma \ref{L:WG} is thus completed.
\begin{remark}\label{R:mgo}
The following example of a model satisfying Conditions (H1) and (H2) was provided in \cite{Z}, p. 2242. Consider a stochastic volatility model introduced in (\ref{E:moodk}), with the volatility process that is the unique solution to the equation in
(\ref{E:oo}), under the following restrictions: \\
(i)\,\,The coefficient maps in (\ref{E:oo}) are given by 
\begin{equation}
a(t,s,x)=K_H(t,s)U_1(s,x),\quad c(t,s,x)=K_H(t,s)U_2(s,x), 
\label{E:cd}
\end{equation}
where $0\le s\le t\le 1$ and $x\in\mathbb{R}^d$. \\
(ii)\,\,In (\ref{E:cd}),
$K_H$ is the kernel in the Molchan-Golosov representation of fractional Brownian motion with Hurst parameter $H\in(0,1)$
(see Subsection \ref{SS:GM}). \\
(iii) $U_1$ and $U_2$ are maps from $[0,t]\times\mathbb{R}^d$ into 
$\mathbb{R}^d$ and from $[0,t]\times\mathbb{R}^d$ into $\mathbb{R}^d\times\mathbb{R}^m$, respectively. They satisfy the following
conditions: There exist $\eta> 0$ and $C> 0$ such that
\begin{equation}
||U_1(s,x)-U_1(s^{\prime},x)||_d+||U_2(s,x)-U_2(s^{\prime},x)||_{d\times m}\le C|s-s^{\prime}|^{\eta}(1+||x||_d)
\label{E:cd0}
\end{equation}
and
\begin{equation}
||U_1(s,x)-U_1(s,y)||_d+||U_2(s,x)-U_2(s,y)||_{d\times m}\le C||x-y||_d
\label{E:cd1}
\end{equation}
for all $x,y\in\mathbb{R}^d$ and $s,s^{\prime}\in[0,1]$ (see p. 2242 in \cite{Z}). 

Since the model described above satisfies Conditions (H1) and (H2), Theorem \ref{T:nmn} holds for the log-price process in the model.
\end{remark}
\begin{remark}\label{R:pte}
The equation in (\ref{E:oo}), with the coefficient functions $a$ and $c$ similar to those in (\ref{E:cd}), was studied in \cite{CD} 
in the case where $d=m=1$. 
It was established in \cite{CD} that there exists a unique continuous strong solution to this equation if $U_1$ and $U_2$ satisfy 
(\ref{E:cd1}), and the kernel $K$ is such that the operator $f\mapsto\int_0^tK(t,s)f(s)ds$ has certain smoothing properties on $L^1[0,T]$ 
and $L^2[0,T]$.
\end{remark}
\begin{remark}\label{R:bf}
In \cite{BFGMS}, the following one-factor drift-less stochastic volatility model was considered (we use the notation adopted in the present paper):
\begin{align*}
&dS_t=S_t\sigma(Y_t)(\bar{\rho}dW_t+\rho dB_t) \\
&Y_t=y+\int_0^TK_H(t,s)U_1(Y_s)ds+\int_0^tK_H(t,s)U_2(Y_s)dB_s
\end{align*}
where the function $f$ is smooth, while the functions $U_1$ and $U_2$ are from the space $\mathbb{C}^{(3)}_b$ (see \cite{BFGMS}, p. 813). 
Here the symbol $\mathbb{C}^{(3)}_b$ stands for the space of three times continuously differentiable functions on $\mathbb{R}$, 
with bounded derivatives. It follows that
the functions $U_1$ and $U_2$ satisfy the estimates in (\ref{E:cd0}) and (\ref{E:cd1}). Therefore, if the canonical set-up is used, then 
the small-noise LDP for the scaled log-price process obtained in Corollary 5.5 in \cite{BFGMS} is a special case of small-noise LDPs established in Section \ref{S:snl} of the present paper. One can use Lemma \ref{L:WG}, Theorems \ref{T:rrr}, \ref{T:bed}, and Remark \ref{R:begin}
to justify the previous statement.
\end{remark}

\subsection{Unification: More Volterra Type SDEs}\label{SS:NR}
In \cite{NR}, Nualart and Rovira formulated the following restrictions on the coefficients in (\ref{E:oo}):
\\
$(H_1)$\,The map $a$ is measurable from $\{0\le s\le t\le T\}\times\mathbb{R}^d$ to $\mathbb{R}^d$, while the map $c$ is 
measurable from $\{0\le s\le t\le T\}\times\mathbb{R}^d$ to $\mathbb{R}^{d\times m}$. \\
$(H_2)$\,The maps $a$ and $c$ are Lipschitz in $x$ uniformly in the other variables, that is, 
$$
||c(t,s,x)-c(t,s,y)||_{d\times m}+||a(t,s,x)-a(t,s,y)||_d\le K||x-y||_d
$$
for some constant $K> 0$, all $x,y\in\mathbb{R}^d$, and all $0\le s\le t\le T$. \\
$(H_3)$\,The maps $a$ and $c$ are $\alpha$-H\"{o}lder continuous in $t$ on $[s,T]$ uniformly in the other variables. 
This means that there exists a constant $K> 0$ such that
$$
||c(t,s,x)-c(r,s,x)||_{d\times m}+||a(t,s,x)-a(r,s,x)||_d\le K|t-r|^{\alpha}
$$
for all $x\in\mathbb{R}^d$ and $s\le t,r\le T$ where $0<\alpha\le 1$. \\
$(H_4)$\,There exists a constant $K> 0$ such that
$$
||c(t,s,x)-c(r,s,x)-c(t,s,y)+c(r,s,y)||_{d\times m}\le K|t-r|^{\gamma}||x-y||_d
$$
for all $x,y\in\mathbb{R}^d$ and $T\ge t,r\ge s$ where $0<\gamma\le 1$. \\
$(H_5)$\,$a(t,s,x_0)$ and $c_j(t,s,x_0)$ are bounded. \\
\\
A sample path LDP was established in \cite{NR} for the unique solution to the equation in (\ref{E:sca}) under Conditions $(H_1) - (H_5)$ (see Theorem 1 in \cite{NR}). We will next show that only Conditions $(H_1) - (H_3)$ and $(H_5)$ are needed in order Theorem 1 in \cite{NR} to be true.
\begin{lemma}\label{L:NRW}
Conditions $(H_1)$ - $(H_3)$ and $(H_5)$ used in \cite{NR} imply Conditions (H1) and $(\widehat{H}2)$ formulated in Subsection \ref{SS:Wang}. 
\end{lemma}

\it Proof. \rm Let us assume that Conditions $(H_1)$ - $(H_3)$ and $(H_5)$ hold. Choose $\rho(u)=u$, any $p> 2$, and 
$K_1(t,s)=K_2(t,s)=\mathbb{1}_{\{0\le s\le t\le T\}}$. 
Then, there exists a constant $C> 0$ depending on the constant $K$ in Condition $(H_2)$ and such that 
(\ref{E:1}) and (\ref{E:2}) hold. The equality in (\ref{E:5}) clearly holds, while (\ref{E:45}) also holds with some constant $C> 0$. The validity of (\ref{E:3}) follows from the fact that Conditions $(H_2)$ and $(H_5)$ imply the boundedness of $a(t,s,0)$ and $c_j(t,s,0)$. 
This completes the proof of Condition (H1) used in \cite{W}.

Next, observe that for
$$
F_1(t^{\prime},t,s)=K|t-t^{\prime}|^{\alpha}\mathbb{1}_{\{0\le s\le t\le T\}}
\mathbb{1}_{\{0\le s\le t^{\prime}\le T\}}
$$
and
$$
F_2(t^{\prime},t,s)=K^{\prime}|t-t^{\prime}|^{2\alpha}\mathbb{1}_{\{0\le s\le t\le T\}}\mathbb{1}_{\{0\le s\le t^{\prime}\le T\}},
$$
the estimates in (\ref{E:6}) and (\ref{E:7}) hold true for $s\le t,t^{\prime}\le T$, with $K^{\prime}> 0$ depending on the constant $K$ appearing in Condition $(H_3)$. The previous statement follows from Condition $(H_3)$. 
In addition, it is easy to see that the estimates in (\ref{E:8}) and (\ref{E:9}) are valid. Therefore, Condition $(\rm\widehat{H}2)$
in Remark \ref{R:Zhang} is satisfied.

The proof of Lemma \ref{L:NRW} is thus completed.
\begin{remark}\label{R:oy} 
It follows from Lemma \ref{L:WG}, Lemma \ref{L:NRW}, and Theorem \ref{T:beg1} that the sample path LDP formulated in Theorem 1 in \cite{NR} holds for the solution to the equation in (\ref{E:oo}) under Conditions
$(H_1)$ - $(H_3)$ and $(H_5)$ provided that the canonical set-up is used. Hence, under the previous restrictions, Condition $(H_4)$ can be removed from 
Theorem 1 in \cite{NR}. We do not know whether the same conclusion can be reached 
if the model in \cite{NR} is defined on a general set-up.
\end{remark}
\section{Applications}\label{S:appl}
\subsection{First Exit Times}\label{SS:dcf}
In this subsection, we obtain a large deviation style formula for the distribution function of the first exit time of the log-price process from an open set in $\mathbb{R}^m$. Such formulas go back to the known results on first exit time due to Freidlin and Wentzell  (see \cite{FW,VF1,VF2}). Suppose $X^{(\varepsilon)}_t$, with $t\in[0,T]$, is the scaled log-price process (see (\ref{E:iu})) starting at $x_0\in\mathbb{R}^m$. It will be assumed in the rest of this section that the conditions in Theorem \ref{T:nmn} are satisfied, and, moreover, for all $(t,u)\in[0,T]\times\mathbb{R}^d$, the matrix $\sigma(t,u)$ is invertible. Then, formula (\ref{E:ras}) holds for the rate function $\widetilde{Q}_T$ in Theorem \ref{T:nmn}.
Let $O$ be a proper open subset of $\mathbb{R}^m$ such that $x_0\in O$. 
\begin{definition}\label{D:dfd}
(i)\,\,For every $\varepsilon\in(0,1]$, the first exit time of the scaled log-price process from the set $O$ is defined by 
$
\tau^{(\varepsilon)}=\inf\{s\in(0,T]:X^{(\varepsilon)}_s\notin O\}
$
if the previous set is not empty, and by $\tau^{(\varepsilon)}=\infty$ otherwise. 
(ii)\,\,For every $\varepsilon\in(0,1]$, the first exit time probability function is defined by 
$v_{\varepsilon}(t)=\mathbb{P}(\tau^{(\varepsilon)}\le t)$, $t\in(0,T]$.
\end{definition}

In the book \cite{FW} of Freidlin and Wentzell, the following restriction on an open set $O\subset\mathbb{R}^m$ was used: There exist interior points of the complement of $O$ arbitrarily close to every point of the boundary of $O$ (see \cite{FW}, Example 3.5). The previous condition can be formulated as follows: 
\begin{equation}
\partial O=\partial(\rm ext(O))
\label{E:sds}
\end{equation} 
where $\rm ext(O)$ is the set of interior points of the complement of $O$, 
and, for a set $D\subset\mathbb{R}^m$, the symbol $\partial D$ stands for the boundary of $D$. 

Let us fix $t\in(0,T]$, and put 
$
{\cal A}_t=\{f\in\mathbb{C}_0^m:f(s)\notin O-x_0\,\,\mbox{for some}\,\,s\in(0,t]\}.
$
Then, ${\cal A}_t$ is a closed subset of the space $\mathbb{C}_0^m$. Its interior ${\cal A}_t^{\circ}$ 
consists of the maps $f\in{\cal A}_t$
for which there exists $s< t$ with $f(s)\notin\rm cl\,(O)-x_0$. Here the symbol $\rm cl\,(O)$ stands for the closure 
of the set $O$ in the space $\mathbb{R}^m$. The boundary $\rm bd\,({\cal A}_t)$ of the set ${\cal A}_t$ in the space 
$\mathbb{C}_0^m$ consists of 
the maps $f\in\mathbb{C}_0^m$ which hit the set $\partial O-x_0$ before $t$, or at $s=t$, but never exit the set
$\rm cl\,(O)-x_0$ before $t$. 

A Borel set $A\subset\mathbb{C}_0^m$ is called a set of continuity for the rate function $\widetilde{Q}_T$ if
the following equality holds: 
$$\inf_{g\in A^{\circ}}\widetilde{Q}_T(g)=\inf_{g\in\bar{A}}\widetilde{Q}_T(g).
$$ 
It follows from the LDP in Theorem \ref{T:nmn} that for any set of continuity $A$, 
\begin{equation}
\varepsilon\log\mathbb{P}\left(X^{(\varepsilon)}-x_0\in A\right)=-\inf_{g\in A}\widetilde{Q}_T(g)+o(1)\,\,\mbox{as}\,\,
\varepsilon\rightarrow 0.
\label{E:o1}
\end{equation}

The next theorem provides a sample path large deviation style formula for the first exit time probability function.
\begin{theorem}\label{T:sco}
Suppose an open set $O\subset\mathbb{R}^m$ is such that the condition in (\ref{E:sds}) holds. Then, the set 
${\cal A}_t$ is a set of continuity for the rate function $\widetilde{Q}_T$, and hence
\begin{equation}
\varepsilon\log\mathbb{P}(\tau^{(\varepsilon)}\le t)=-\inf_{g\in{\cal A}_t}\widetilde{Q}_T(g)+o(1)\,\,\,\mbox{as}\,\,
\varepsilon\rightarrow 0.
\label{E:treb}
\end{equation}
\end{theorem}

\it Proof. \rm  In the proof of Theorem \ref{T:sco}, we borrow some ideas from the proof of Theorem 2.16 in \cite{Gul1} and also take into account Example 3.5 in \cite{FW}. Our first goal is to provide a sufficient condition for a Borel set $A\subset\mathbb{C}_0^m$ to be a set of continuity for 
the rate function $\widetilde{Q}_T$. The following statement can be obtained using the continuity of the function 
$\widetilde{Q}_T$ on the space $(\mathbb{H}_0^1)^m$ (see Lemma \ref{L:dff}).
\begin{lemma}\label{L:sce}
Suppose a Borel set $A\subset\mathbb{C}_0^m$ is such that for every $h\in\rm bd\,(A)\cap(\mathbb{H}_0^1)^m$, there exists a sequence 
$h_n\in A^{\circ}\cap(\mathbb{H}_0^1)^m$ for which $\displaystyle{\lim_{n\rightarrow\infty}h_n=h}$ in the space 
$(\mathbb{H}_0^1)^m$. Then, the set $A$ is a set of continuity for the rate function $\widetilde{Q}_T$.
\end{lemma}

Let us continue the proof of Theorem \ref{T:sco}. It is not hard to see that 
\begin{equation}
\{\tau^{(\varepsilon)}\le t\}=\{X^{(\varepsilon)}-x_0\in{\cal A}_t\}.
\label{E:scet}
\end{equation} 
We will next prove that the set ${\cal A}_t$ is a set of continuity for $\widetilde{Q}_T$. Let
$
f\in\rm bd\,({\cal A}_t)\cap(\mathbb{H}_0^1)^m,
$ 
and let $t_0\in(0,t]$ be such that $f(t_0)\in\partial O-x_0$ and $f(u)\in\rm cl\,(O)-x_0$ 
for all $u\in(0,t]$.
Using (\ref{E:sds}), we see that for every $n\ge 1$ there is a point $x_n\in\rm ext(O)-x_0$ such that
$||x_n-f(t_0)||_m\le\frac{1}{n}$ for all $n\ge 1$. Define a sequence of functions on $[0,T]$ by the following formula:
$f_n(t)=f(t)+\frac{t}{t_0}(x_n-f(t_0))$, $n\ge 1$.
It is easy to see that for every $n\ge 1$, we have $f_n\in(\mathbb{H}_0^1)^m$. Moreover, $f_n\rightarrow f$ in  
$(\mathbb{H}_0^1)^m$ as $n\rightarrow\infty$. Since $f_n(t_0)=x_n\in\rm ext(O)-x_0$, the map $f_n$ exits the set $\rm cl\,(O)-x_0$
before $t$. This is clear for $t_0< t$, while for $t_0=t$, we can use the continuity of $f$. It follows that $f_n\in{\cal A}_t
^{\circ}\cap(\mathbb{H}_0^1)^m$ for all $n\ge 1$. Next, using Lemma \ref{L:sce}, we see that ${\cal A}_t$ is a set of continuity 
for the rate function $\widetilde{Q}_T$. Finally, the equality in (\ref{E:treb}) follows from (\ref{E:o1}). 
\subsection{Binary Barrier Options}\label{SS:bbo}
Our goal in the present subsection is to obtain a large deviation style formula in the small-noise regime for multidimensional binary barrier options. 
Suppose that the model in (\ref{E:moodk}) describes the dynamics of price processes associated with a portfolio of correlated assets. 
Let $\varepsilon\in(0,1]$, and consider the scaled $m$-dimensional asset price process $t\mapsto S^{(\varepsilon)}_t$  
and the scaled log-price process $t\mapsto X^{(\varepsilon)}_t$. The latter process is given by the expression in (\ref{E:iu}).

We will study the small-noise asymptotic behavior of binary up-and-in barrier options. 
Similar results can be obtained for up-and-out, down-and-in, and down-and-out options. We refer the reader to \cite{Gul1} where one-dimensional Gaussian models are considered. 

Denote by $\mathbb{R}^m_{+}$ the subset of $\mathbb{R}^m$ consisting of all the vectors $s=(s_1,\cdots,s_m)\in\mathbb{R}^m$ such that $s_i> 0$ for all $1\le i\le m$, and let 
$O\subset\mathbb{R}^m_{+}$ be an open set satisfying the condition in (\ref{E:sds}). 
The boundary $\partial O$ of the set $O$ will play the role of the barrier. Throughout the present section we assume that
the model in (\ref{E:moodk}) satisfies the restrictions imposed in Theorem \ref{T:nmn}.

Let us suppose that for every $\varepsilon\in(0,1]$ the initial condition $s_0$ for the process $t\mapsto S_t^{(\varepsilon)}$
is such that $s_0\in O$. 
\begin{definition}\label{D:ddd}
Let $O$ be an open set in $\mathbb{R}^m_{+}$ satisfying the condition in (\ref{E:sds}).
In a small-noise setting, a binary up-and-in barrier option pays a fixed amount of cash, say one dollar, if the $m$-dimensional asset price process $S^{(\varepsilon)}$
hits the barrier $\partial O$ at some time during the life of the option. 
\end{definition}

For every $\varepsilon\in(0,1]$, the payoff of a binary up-and-in barrier option at the maturity is $\{S_t^{(\varepsilon)}\in\partial O
\,\mbox{for some}\,t\in[0,T]\}$.
Therefore, the price $B(\varepsilon)$ 
of the option at $t=0$ is given by
\begin{equation}
B(\varepsilon)=e^{-rT}\mathbb{P}(S_t^{(\varepsilon)}\in\partial O
\,\mbox{for some}\,t\in[0,T])
\label{E:pric}
\end{equation}
where $r> 0$ is the interest rate. 
It is not hard to see using (\ref{E:pric}) that
\begin{equation}
B(\varepsilon)=e^{-rT}\mathbb{P}(S_t^{(\varepsilon)}\notin O\,\mbox{for some}\,t\in[0,T]).
\label{E:yry}
\end{equation}
It is clear that a binary up-and-in barrier option contract depends on $r$, $s_0$, the maturity $T$ of the option, 
and the barrier $\partial O$. We will study the asymptotic behavior of the price of the barrier option as $\varepsilon\rightarrow 0$. 
Barrier options are path-dependent options.

Our next goal is to rewrite the expression on the right-hand side of (\ref{E:yry}) in terms of the log-price process. 
The resulting equality is as follows:
\begin{equation}
B(\varepsilon)=e^{-rT}\mathbb{P}(X_t^{(\varepsilon)}-x_0\notin\widetilde{O}-x_0\,\mbox{for some}\,t\in[0,T])
\label{E:ugugu}
\end{equation}
where $\widetilde{O}$ is the open subset of $\mathbb{R}^m$ defined by 
$$
\widetilde{O}=\{x=(x_1,\cdots,x_m)\in\mathbb{R}^m:(e^{x_1},\cdots,e^{x_m})\in O\}.
$$
It is easy to see that the set $\widetilde{O}-x_0$ satisfies the condition in (\ref{E:sds}).

The price of a binary up-and-in barrier option is related to the exit time probability function $\tau_{\widetilde{O}}$ 
of the log-price process 
$X^{(\varepsilon)}$ from the set $\widetilde{O}$. Indeed, it was shown in the proof of Theorem \ref{T:sco} that 
\begin{equation}
\{\tau^{(\varepsilon)}_{\widetilde{O}}\le T\}=\{X^{(\varepsilon)}-x_0\in{\cal A}_T\},
\label{E:gtu}
\end{equation} 
where
\begin{equation}
{\cal A}_T=\{f\in\mathbb{C}_0^m:f(s)\notin\widetilde{O}-x_0\,\,\mbox{for some}\,\,s\in[0,T]\}.
\label{E:ljl}
\end{equation}
(see (\ref{E:scet})).

The next assertion provides a large deviation style formula in the small-noise regime for the price $B(\varepsilon)$ of the 
binary up-and-in barrier option given by (\ref{E:yry}). 
\begin{theorem}\label{T:333}
The following asymptotic formula holds:
$$
\varepsilon\log B(\varepsilon)=-\inf_{g\in{\cal A}_T}\widetilde{Q}_T(g)+o(1)\,\,\,\mbox{as}\,\,
\varepsilon\rightarrow 0,
$$
where $\widetilde{Q}_T$ is the rate function in Theorem \ref{T:nmn}, and the set ${\cal A}_T$ is defined by (\ref{E:ljl}). 
\end{theorem}

Theorem \ref{T:333} can be easily derived from (\ref{E:ugugu}), (\ref{E:gtu}), (\ref{E:ljl}), and Theorem \ref{T:sco}.
\subsection{Call Options}\label{SS:call}
In this subsection, we consider the model in (\ref{E:mood}) with $m=1$ and $d\ge 1$. It will be assumed that $b(t,u)=r$ where $r\ge 0$ is the interest rate. This means that we turn our attention to $d$-factor stochastic volatility models of financial mathematics. More precisely, the models that we study in this subsection are the following:
\begin{equation}
\frac{dS_t}{S_t}=rdt+\sigma(t,\widehat{B}_t)(\bar{\rho}dW_t+\rho dB_t),\quad t\in[0,T].
\label{E:form}
\end{equation}
In (\ref{E:form}), $\rho\in(-1,1)$ is the correlation parameter, $\bar{\rho}=\sqrt{1-\rho^2}$, and $\widehat{B}$ is the $d$-dimensional volatility process introduced in Definition \ref{D:vo}. The scaled version of the model in (\ref{E:form}) is as follows:
\begin{equation}
\frac{dS_t^{(\varepsilon)}}{S_t^{(\varepsilon)}}=rdt+\sqrt{\varepsilon}
\sigma(t,\widehat{B}_t^{(\varepsilon)})(\bar{\rho}dW_t+\rho dB_t),\quad t\in[0,T]
\label{E:forme}
\end{equation}
where $\varepsilon\in(0,1]$ and $\widehat{B}^{(\varepsilon)}$ is the scaled volatility process (see Definition \ref{D:dds}). 

The price of the call option in the small-noise regime is the following function of the strike $K> 0$, the maturity $T> 0$, and the small-noise parameter $\varepsilon\in(0,1]$:
\begin{equation}
C(\varepsilon,T,K)=\mathbb{E}[(S_T^{(\varepsilon)}-K)^{+}]
\label{E:call}
\end{equation}
 We assume that $K$ and $T$ are fixed and study the asymptotic 
behavior of the call price as $\varepsilon\rightarrow 0$. The following assumption will be used in the sequel: \\
\\
\it Assumption $B$. \rm For every $\alpha> 0$ there exists $\varepsilon_0\in(0,1]$ depending only on $\alpha$
and such that for all 
$0<\varepsilon<\varepsilon_0$, the following estimate holds true:
\begin{equation}
\mathbb{E}\left[\exp\left\{\alpha\int_0^T\sigma(s,\hat{B}^{(\varepsilon)}_s)^2ds\right\}\right]\le M
\label{E:um}
\end{equation}
where $M> 0$ is a constant depending only on $\alpha$. 
\begin{lemma}\label{L:mar}
Suppose Assumption $B$ holds. Then, for every $C> 0$ there exists $\varepsilon_1\in(0,1]$ depending on $C$ and such that 
for every $\widehat{C}$, with $0<\widehat{C}\le C$, the stochastic exponential
$$
{\cal E}(t,\varepsilon,\widehat{C})=\exp\left\{-\frac{1}{2}\widehat{C}^2\int_0^t\sigma(s,\hat{B}^{(\varepsilon)}_s)^2ds
+\widehat{C}\int_0^t\sigma(s,\hat{B}^{(\varepsilon)}_s)
(\bar{\rho}dW_t+\rho dB_t)\right\},\quad t\in[0,T]
$$ 
is an $\{{\cal F}_t\}$-martingale for all $\varepsilon\le\varepsilon_1$. 
\end{lemma}

\it Proof. \rm Set $\alpha=\frac{1}{2}C^2$ and $U(t,\varepsilon,C)=C\sigma(t,\hat{B}^{(\varepsilon)}_t)$, $t\in[0,T]$. 
Then, Assumption B implies that there exists $\varepsilon_1\in(0,1]$ such that Novikov's condition is satisfied for the 
process $t\mapsto U(t,\varepsilon,C)$, with $0<\varepsilon\le\varepsilon_1$. It follows that the stochastic exponential in Lemma 
\ref{L:mar} is an $\{{\cal F}_t\}$-martingale for all $\varepsilon\le\varepsilon_1$. The same proof works for any $\widehat{C}< C$
since if the inequality in (\ref{E:um}) holds for some $\alpha> 0$ and $\varepsilon<\varepsilon_0$, it also holds for any 
$0<\hat{\alpha}<\alpha$ and $\varepsilon<\varepsilon_0$.

The proof of Lemma \ref{L:mar} is thus completed.

Using the Dol\'{e}ans-Dade formula, we see that
for every $\varepsilon\in(0,1]$, the discounted scaled asset price process is given by
\begin{equation}
e^{-rt}S_t^{(\varepsilon)}=s_0\exp\left\{-\frac{\varepsilon}{2}\int_0^t\sigma(s,\widehat{B}_s^{(\varepsilon)})^2ds 
 +\sqrt{\varepsilon}\int_0^t\sigma(s,\widehat{B}_s^{(\varepsilon)})(\bar{\rho}dW_t+\rho dB_t)\right\}
\label{E:lab1}
\end{equation}
where $t\in[0,T]$.
\begin{remark}\label{E:por}
Applying Lemma \ref{L:mar}, with $C=\sqrt{\varepsilon}$, we see that if Assumption B holds, then for small enough values of $\varepsilon$, 
the discounted asset price process in (\ref{E:lab1}) is a martingale. It follows that the measure $\mathbb{P}$ in the model defined in 
(\ref{E:forme}), with small enough values of $\varepsilon$, is risk-neutral. 
\end{remark}
\begin{remark}\label{R:ldf}
It is not hard to see, using the formula in (\ref{E:oi}), that for the model in (\ref{E:form}), with 
\begin{equation}
\sigma(s,z)\neq 0\quad\mbox{for all}\quad (s,z)\in[0,T]\times\mathbb{R}^d,
\label{E:rew}
\end{equation}
the rate function $\widetilde{I}_T$ can be represented as follows:
\begin{equation}
\widetilde{I}_T(x)=\frac{1}{2}\inf_{f\in\mathbb{H}_0^1}\left[\frac{(x-r
-\rho\int_0^T\sigma(s,\widehat{f}(s))\dot{f}(s)ds)^2}{\bar{\rho}^2\int_0^T\sigma(s,\widehat{f}(s))^2ds}
+\int_0^T\dot{f}(t)^2dt\right].
\label{E:oit}
\end{equation}
Moreover, the function in (\ref{E:oit}) is continuous on $\mathbb{R}$ (see Lemma \ref{L:poi}).
\end{remark}
\begin{remark}\label{R:ldd}
Suppose we do not impose the restriction in (\ref{E:rew}) on the model in (\ref{E:form}). Then, the formula in (\ref{E:oit})
holds for all $x\in\mathbb{R}$ such that the set $Q_2$ is empty (see Remark \ref{R:begin}). If $Q_2\neq\emptyset$, then (\ref{E:fre}) implies that the set $Q_3(x)$ can be nonempty only if $x=r$. Therefore, the formula in (\ref{E:oit}) holds for all $x> r$. In addition, the function
$\widetilde{I}_T$ is continuous on the set $x> r$. The previous statement can be established exactly as Lemma \ref{L:poi}.
\end{remark}

The next assertion provides an LDP-style formula for the price of the scaled call option.
\begin{theorem}\label{T:ttr}
Suppose Assumption A and Assumptions (C1) -- (C7) hold true for the model in (\ref{E:form}) defined on the canonical set-up. Further suppose that the condition in (\ref{E:rew}) holds and Assumption $B$ is satisfied. Then, 
the following asymptotic formula is valid:
\begin{equation}
\varepsilon\log C(\varepsilon,T,K)=-\inf_{x\ge k}\widetilde{I}_T(x)+o(1)
\label{E:tote}
\end{equation}
as $\varepsilon\rightarrow 0$ where $k$ is the log-moneyness defined by $k=
\log\frac{K}{s_0}$, and $\widetilde{I}_T$ is the good rate function given by (\ref{E:oi}).
\end{theorem}

\it Proof. \rm We will first establish a lower large deviation estimate for the scaled call price. It is not hard to see that for every 
$\delta> 0$,
$$
C(\varepsilon,T,K)\ge\delta\mathbb{P}(S_T^{(\varepsilon)}> K+\delta)
=\delta\mathbb{P}\left(X_T^{(\varepsilon)}-x_0> \log\frac{K+\delta}{s_0}\right).
$$
It follows from the LDP in Theorem \ref{T:ha} that
\begin{equation}
\liminf_{\varepsilon\rightarrow 0}\varepsilon\log C(\varepsilon,T,K)\ge-\inf_{x>\log\frac{K+\delta}{s_0}}\widetilde{I}_T(x).
\label{E:cri}
\end{equation}
Since (\ref{E:cri}) holds for all $\delta> 0$ and the function $\widetilde{I}_T$ is continuous on $\mathbb{R}$ (see Remark \ref{R:ldf}), 
we derive the following estimate from (\ref{E:cri}):
\begin{equation}
\liminf_{\varepsilon\rightarrow 0}\varepsilon\log C(\varepsilon,T,K)\ge-\inf_{x>k}\widetilde{I}_T(x).
\label{E:rer}
\end{equation}
The formula in (\ref{E:rer}) is a lower LDP-style estimate for the scaled call price.

We will next obtain a similar upper estimate. We borrow some ideas from the proof on page 1131 in \cite{GVol}. 
Let $p> 1$ and $q=\frac{p}{p-1}$. Using H\"{o}lder's inequality, we get
$$
C(\varepsilon,T,K)\le\left\{\mathbb{E}\left[(S_T^{(\varepsilon)})^p\right]\right\}^{\frac{1}{p}}
\mathbb{P}\left(S_T^{(\varepsilon)}\ge K\right)^{\frac{1}{q}}.
$$
Therefore,
\begin{align}
\limsup_{\varepsilon\rightarrow 0}\varepsilon\log C(\varepsilon,T,K)&\le\frac{1}{p}\limsup_{\varepsilon\rightarrow 0}
\varepsilon\log\mathbb{E}\left[(S_T^{(\varepsilon)})^p\right] \nonumber \\
&\quad+\frac{1}{q}\limsup_{\varepsilon\rightarrow 0}
\varepsilon\log\mathbb{P}\left(S_T^{(\varepsilon)}\ge K\right).
\label{E:eds}
\end{align}
The second term on the right-hand side of (\ref{E:eds}) can be estimated using the LDP in Theorem \ref{T:ha} 
as in the proof of (\ref{E:rer}). This gives
\begin{equation}
\limsup_{\varepsilon\rightarrow 0}\varepsilon\log\mathbb{P}\left(S_T^{(\varepsilon)}\ge K\right)\le-\inf_{x\ge k}\widetilde{I}_T(x).
\label{E:gg}
\end{equation}
Note that, by the continuity of the rate function, we have
\begin{equation}
\inf_{x>k}\widetilde{I}_T(x)=\inf_{x\ge k}\widetilde{I}_T(x)
\label{E:xx}
\end{equation}

Our next goal is to estimate the first term on the right-hand side of (\ref{E:eds}). 
Using (\ref{E:lab1}), we see that for every $p> 1$,
all $t\in[0,T]$, and all $\varepsilon\in(0,1]$,
$$
\mathbb{E}\left[(S_t^{(\varepsilon)})^p\right]\le s_0^pe^{prt}\mathbb{E}\left[{\cal E}(t,\varepsilon, 2p)\right]^{\frac{1}{2}}
\mathbb{E}\left[\exp\left\{(2p^2-p)\int_0^t\sigma(s,\widehat{B}_s^{(\varepsilon)})^2ds\right\}\right]^{\frac{1}{2}}.
$$
It follows from Assumption B and Lemma \ref{L:mar} that there exist $\varepsilon_p\in(0,1]$ and $M_p> 0$ depending on $p$ and such that
\begin{equation}
\mathbb{E}\left[(S_t^{(\varepsilon)})^p\right]\le s_0^pe^{pr}M_p
\label{E:trtt}
\end{equation}
for all $t\in[0,T]$ and $0<\varepsilon\le\varepsilon_p$. It is easy to see that (\ref{E:trtt}) implies the following equality:
\begin{equation}
\limsup_{\varepsilon\rightarrow 0}\varepsilon\log\mathbb{E}\left[(S_T^{(\varepsilon)})^p\right]=0
\label{E:sw}
\end{equation}
for every $p> 1$. Now, using (\ref{E:eds}), (\ref{E:gg}), and (\ref{E:sw}), and assuming that $q\rightarrow 1$, we get the upper estimate
\begin{equation}
\limsup_{\varepsilon\rightarrow 0}\varepsilon\log C(\varepsilon,T,K)\le-\inf_{x\ge k}\widetilde{I}_T(x).
\label{E:now}
\end{equation} 
Finally, using (\ref{E:rer}), (\ref{E:xx}), and (\ref{E:now}), we obtain formula (\ref{E:tote}).

This completes the proof of Theorem \ref{T:ttr}.
 
The restrictions in Theorem \ref{T:ttr} include the condition in (\ref{E:rew}).
The class of stochastic volatility models for which this condition is satisfied includes stochastic volatility models, with exponential volatility function, e.g., the Scott model (see \cite{Scott}), the rough Bergomi model (see \cite{BFG}), and the super rough Bergomi model 
(see \cite{Gul1} and \cite{BHP}). On the other hand, the condition in (\ref{E:rew}) is not satisfied for certain classical stochastic volatility models, for instance, the Stein and Stein model (see \cite{SS}) and the Heston model (see \cite{Hes}). In the former model, the volatility function is $\sigma(x)=x$, for all
$x\in\mathbb{R}$, while in the latter one, the volatility function is given by $\sigma(x)=\sqrt{x}$, for all $x\ge 0$. 

The next assertion explains what happens if we remove the condition in (\ref{E:rew}) from Theorem \ref{T:ttr}.
\begin{theorem}\label{T:las}
Suppose Assumption A and Assumptions (C1) -- (C7) hold true for the model in (\ref{E:form}) defined on the canonical set-up. Further suppose 
that Assumption $B$ is satisfied. Then, the formula 
$$
\varepsilon\log C(\varepsilon,T,K)=-\inf_{x\ge k}\widetilde{I}_T(x)+o(1)
$$
as $\varepsilon\rightarrow 0$ holds for all $K>s_0e^{r}$.
\end{theorem}

\it Proof. \rm Theorem \ref{T:las} can be established using the same methods as in the proof of Theorem \ref{T:ttr} and by taking into account Remark \ref{R:ldd}. This gives the formula in (\ref{E:tote}) for all $k> r$ that is equivalent to $K> s_0e^{r}$.

The proof of Theorem \ref{T:las} is thus completed.

\subsection{Implied Volatility in the Small-Noise Regime}\label{SS:sniv}
Our main objective in the present subsection is to use fundamental results of Gao and Lee (see \cite{GL}) which provide relations between the asymptotic behavior of the logarithm of the call price with respect to various parametizations and the asymptotic behavior of the implied volatility.

Let us consider the scaled stochastic volatility model described in (\ref{E:forme}).
 We assume that $s_0=1$ and $r=0$. The previous normalization is employed in \cite{GL}. It follows that $k=\log K$. 
Recall that the small-noise call option price 
$C(\varepsilon,T,K)$ is defined by (\ref{E:call}). Set
\begin{equation}
\widetilde{C}(\varepsilon,T,k)=C(\varepsilon,T,e^k)=\mathbb{E}[(S_T^{(\varepsilon)}-e^k)^{+}],
\label{E:set}
\end{equation}
with $k\in\mathbb{R}$. 

It follows from Theorem \ref{T:las}
that for every $k> 0$,
\begin{equation}
\log\frac{1}{\widetilde{C}(\varepsilon,T,k)}=\frac{1}{\varepsilon}\inf_{x\ge k}\widetilde{I}_T(x)+o(\frac{1}{\varepsilon})
\label{E:pou}
\end{equation}
as $\varepsilon\rightarrow 0$. Recall that in the present subsection, we assume that $r=0$.

In \cite{GL}, the following formula is used for the call price $C_{BS}$ in the Black-Scholes model: (see (3.1) in \cite{GL}):
\begin{equation}
C_{BS}(v,k)=\frac{1}{\sqrt{2\pi}}\int_{-\infty}^{d_1}e^{-\frac{y^2}{2}}dy-\frac{e^k}{\sqrt{2\pi}}\int_{-\infty}^{d_2}
e^{-\frac{y^2}{2}}dy
\label{E:ddd}
\end{equation}
where 
$$ 
d_1=\frac{-k+\frac{1}{2}v^2}{v}\quad\mbox{and}\quad
d_2=\frac{-k-\frac{1}{2}v^2}{v}.
$$
The formula in (\ref{E:ddd}) represents the call price in the Black-Scholes model as a function of the log-strike $k\in\mathbb{R}$ and the dimensionless implied volatility $v> 0$ (see formula (3.1) in \cite{GL}). 

Let $\sigma> 0$ be the volatility parameter in the classical Black-Scholes model. It is clear that if we replace $v$ by 
$\sigma\sqrt{T}$ in (\ref{E:ddd}), the formula in (\ref{E:ddd}) becomes the classical Black-Scholes formula for the call price.

A generally accepted small-noise parametrization of the implied volatility $v$ is as follows: $v(\varepsilon)
=\sqrt{\varepsilon}\sigma$. Our next goal is to take into account the scaling in the Black-Scholes formula. Set
\begin{equation}
\widehat{C}_{BS}(\varepsilon,T,k,\sigma)=\frac{1}{\sqrt{2\pi}}\int_{-\infty}^{d_1(\varepsilon,k,\sigma)}e^{-\frac{y^2}{2}}dy-\frac{e^k}{\sqrt{2\pi}}
\int_{-\infty}^{d_2(\varepsilon,k,\sigma)}e^{-\frac{y^2}{2}}dy
\label{E:ddu}
\end{equation}
where
$$
d_1(\varepsilon,T,k,\sigma)=\frac{-k+\frac{1}{2}\varepsilon T\sigma^2}{\sqrt{\varepsilon T}\sigma}\quad\mbox{and}\quad
d_2(\varepsilon,T,k,\sigma)=\frac{-k-\frac{1}{2}\varepsilon T\sigma^2}{\sqrt{\varepsilon T}\sigma}.
$$ven
\begin{definition}\label{D:sc}
Let $\widetilde{C}$ be the scaled call price function defined in (\ref{E:set}). Given $k\in\mathbb{R}$, $T> 0$, and $\varepsilon\in(0,1]$, the implied volatility 
in the small-noise setting associated with the function $\widetilde{C}$ is the value of the volatility parameter $\sigma$ in (\ref{E:ddu}) for which 
$\widetilde{C}(\varepsilon,T,k)=\widehat{C}_{BS}(\varepsilon,T,k,\sigma)$. The implied volatility will be denoted by $V(\varepsilon,T,k)$. 
\end{definition}

The formulas above show that we consider the case of fixed $k$ and $T$, and the dimensionless implied volatility $v$ parametrized as follows:
\begin{equation}
v(\varepsilon)=\sqrt{\varepsilon T}V(\varepsilon,T,k),\,\,\varepsilon\in(0,1]. 
\label{E:ts}
\end{equation}

The next statement provides an asymptotic formula for the implied volatility $V$ as $\varepsilon\rightarrow 0$.
\begin{theorem}\label{T:ter}
Suppose Assumption A and Assumptions (C1) -- (C7) hold true for the model in (\ref{E:form}) defined on the canonical set-up. Further suppose that Assumption $B$ is satisfied. Then, the following formula holds for every $k> 0$ and $T> 0$:
\begin{equation}
\lim_{\varepsilon\rightarrow 0}V(\varepsilon,T,k)=\frac{k}{\sqrt{2T\inf_{x\ge k}\widetilde{I}_T(x)}}.
\label{E:pr} 
\end{equation}
\end{theorem}

\it Proof. \rm It follows from (\ref{E:pou}), (\ref{E:ddu}), and Definition \ref{D:sc} that
\begin{equation}
\log\frac{1}{C_{BS}(v(\varepsilon),k)}=\frac{1}{\varepsilon}\inf_{x\ge k}\widetilde{I}_T(x)+o\left(\frac{1}{\varepsilon}\right),\quad
\varepsilon\rightarrow 0,
\label{E:lses}
\end{equation}
where $v(\varepsilon)$ is given by (\ref{E:ts}). It is not hard to see that the conditions in Colollary 7.2 in \cite{GL} hold. It follows from (\ref{E:lses})
and formulas (7.6) and (7.8) in \cite{GL} that as $\varepsilon\rightarrow 0$,
$$
v(\varepsilon)=\frac{\sqrt{\varepsilon}k}{\sqrt{2\inf_{x\ge k}\widetilde{I}_T(x)}}+o(\sqrt{\varepsilon}).
$$
Next, using (\ref{E:ts}), we see that the formula in (\ref{E:pr}) can be obtained from the previous equality.
\begin{corollary}\label{C:unco}
Suppose the conditions in Theorem \ref{T:ter} hold. Suppose also that the restriction in (\ref{E:rew}) is satisfied. If $r=0$, 
and the model is uncorrelated ($\rho=0$), then for all $k> 0$ and $T> 0$,
\begin{equation}
\lim_{\varepsilon\rightarrow 0}V(\varepsilon,T,k)=\frac{k}{\sqrt{2T\widetilde{I}_T(k)}}.
\label{E:proo} 
\end{equation}
\end{corollary}

\it Proof. \rm It follows from the conditions in Corollary \ref{C:unco} that the rate function in (\ref{E:oit}) is increasing on $(0,\infty)$.
Since the rate function is continuous (see Lemma \ref{L:poi}), we have 
\begin{equation}
\inf_{x\ge k}\widetilde{I}_T(x)=\widetilde{I}_T(k).
\label{E:ooo}
\end{equation} 
Now, it is clear that (\ref{E:proo}) follows from (\ref{E:pr}) and (\ref{E:ooo}).

\subsection{A Toy Model}\label{SS:TM}
In this subsection, we consider a simple model (a toy model) and discuss the applicability of Theorem \ref{T:ter} to the toy model. We choose
a special uncorrelated SABR model as the toy model.

The SABR model was introduced in \cite{HKLW}. The toy model analyzed in the present subsection is the following special case of the SABR model:
\begin{align}
&dS_t=X_tS_tdW_t \nonumber \\
&dX_t=X_tdB_t
\label{E:S}
\end{align}
where $0\le t\le T$, $S_0=1$, $X_0=1$, and $W$ and $B$ are independent standard Brownian motions. The toy model in (\ref{E:S}) is the SABR model with
$\nu=1$, $\beta=1$, and $\rho=0$ (see the notation in \cite{HKLW}).
 
The toy model described in (\ref{E:S}) is one of the Gaussian models used in the present paper. Indeed, it follows from (\ref{E:S}) that
\begin{equation}
dS_t=\sigma(t,B_t)S_tdW_t,\quad 0\le t\le T
\label{E:SA}
\end{equation}
where
\begin{equation}
\sigma(t,u)=\exp\{-\frac{1}{2}t+u\},\,\,(t,u)\in[0,T]\times\mathbb{R}^1.
\label{E:bel}
\end{equation}
The Gaussian model in (\ref{E:SA}) is drift-less and uncorrelated ($\rho=0$). Moreover, the volatility process $\widehat{B}$ in the toy model is the standard Brownian motion $B$, and the volatility function $\sigma$ is given by (\ref{E:bel}). The scaled volatility process is as follows: 
$\widehat{B}^{(\varepsilon)}_t=\sqrt{\varepsilon}B_t$, with $0\le t\le T$ and $\varepsilon\in(0,1]$.
\begin{remark}\label{R:HW}
The SABR model in (\ref{E:S}) is a special case of the Hull-White model studied in \cite{HW} (see also \cite{G}).
\end{remark}

Our next goal is to apply Theorem \ref{T:ter} to the toy model.
It is easy to see that Assumption A is satisfied. The reader can find the formulation of Assumption A after 
Definition \ref{D:dd}. It is also not hard to check the validity of the estimate in (\ref{E:um}). Therefore, Assumption B is satisfied. Moreover, 
Assumptions (C1) -- (C7) formulated in Section \ref{S:VP} are also satisfied. The previous statement follows from the fact that the toy model is Gaussian and from Theorem \ref{T:cxc}.

We will next describe the mapping $f\mapsto\widehat{f}$ associated with the toy model. By taking into account the equalities
$\widehat{B}_t=B_t=\int_0^tdB_s$, we see that for the toy model, the equation in (\ref{E:t2}) has the following form: 
$\eta_f=\int_0^tf(s)ds$. The previous statement can be derived from the fact that in the equation in (\ref{E:2a}), $c(t,s)=1$ for $s\le t$, 
$y=0$, and $a=0$. 
Now, using Definition \ref{D:lip}, we see that $\widehat{f}=f$, for all $f\in\mathbb{H}_0^1$.

Since in the toy model, $\sigma(t,u)\neq 0$, for all $(t,u)\in[0,T]\times\mathbb{R}$, the rate function 
$\widetilde{I}_T$ appearing in Theorem \ref{T:ter} has the following form:
\begin{equation}
\widetilde{I}_T(x)=\frac{1}{2}\inf_{f\in\mathbb{H}_0^1}\left[\frac{x^2}
{\int_0^T\exp\{-t+2f(t)\}dt}
+\int_0^T\dot{f}(t)^2dt\right],\quad x\in\mathbb{R}.
\label{E:oioi}
\end{equation}
The previous formula can be obtained from (\ref{E:oi}). It follows from Corollary \ref{C:unco} that the formula in (\ref{E:proo})
holds for the toy model.

\begin{remark}\label{R:oh}
It would be interesting to find a simple explicit representation for the rate function $\widetilde{I}_T(k)$, $k> 0$. However, we do not know how to obtain such a representation. In the remaining part of the present subsection, we will establish estimates from above and below for 
the rate function in the toy model. 
\end{remark}

We will first obtain an estimate from below for the rate function. It is not difficult to prove that for any $f\in\mathbb{H}_0^1$ the following inequality holds:
$$
T\int_0^T\dot{f}(t)^2dt\ge(\max_{t\in[0,T]}|f(t)|)^2. 
$$
Therefore, (\ref{E:oioi}) implies that for every $k> 0$ we have
$$
\widetilde{I}_T(k)\ge\frac{1}{2}\min_{a\ge 0}\left[\frac{k^2}{(1-e^{-T})e^{2a}}
+\frac{1}{T}a^2\right].
$$
The minimization problem in the formula above can be easily solved. The resulting estimate is as follows:
\begin{equation}
\widetilde{I}_T(k)\ge\frac{1}{2}\left[\frac{k^2}
{(1-e^{-T})e^{2a(k)}}+\frac{1}{T}a(k)^2\right],
\label{E:ak1} 
\end{equation}
where $a(k)> 0$ is such that
\begin{equation}
a(k)e^{2a(k)}=\frac{Tk^2}{1-e^{-T}}.
\label{E:ak2}
\end{equation}
It is easy to prove that the equalities in (\ref{E:ak1}) and (\ref{E:ak2}) imply the following estimate:
\begin{equation}
\widetilde{I}_T(k)\ge\frac{1}{2T}\left[a(k)+a(k)^2\right].
\label{E:ak7} 
\end{equation}

Define a function on the set $[0,\infty)$ by $h(u)=ue^{u}$. The function $h$ is strictly increasing and continuous. 
It follows from (\ref{E:ak2}) that 
\begin{equation}
a(k)=\frac{1}{2}h^{-1}\left(\frac{2Tk^2}{1-e^{-T}}\right).
\label{E:ak3}
\end{equation}
In (\ref{E:ak3}), the symbol $h^{-1}$ stands for the inverse function of the function $h$.
\begin{remark}\label{R:Lag}
Using the Lagrange inversion formula, we can represent the function $h^{-1}(y)$ by its Taylor series at $y=0$. This gives
\begin{equation}
h^{-1}(y)=\sum_{n=1}^{\infty}(-1)^{n-1}\frac{n^{n-1}}{n!}y^n.
\label{E:lag}
\end{equation}
It follows from the ratio test that the radius of convergence $R$ of the alternating power series appearing in (\ref{E:lag}) is given by
$R=\frac{1}{e}$. Moreover, it is not hard to prove that for $0< y<\frac{1}{e}$, the absolute values of the terms of the series
in (\ref{E:lag}) decrease. The proof is based on the fact that the sequence $(1+\frac{1}{n})^n$ increases and its limit as 
$n\rightarrow\infty$ is equal to $e$. It follows from what was said above that 
\begin{equation}
h^{-1}(y)\ge y-y^2,
\label{E:opt}
\end{equation}
provided that $0< y<\frac{1}{e}$.
\end{remark}

\begin{remark}\label{R:pto}
By taking into account the previous reasoning, especially the estimates in (\ref{E:ak7}), (\ref{E:ak3}), and (\ref{E:opt}), we obtain an estimate from below for the rate function $\widetilde{I}_T$ in the toy model.
\end{remark}
\begin{lemma}\label{L:use}
Suppose $T> 0$ and $0< k<\sqrt{\frac{1-e^{-T}}{2Te}}$. Then, the following estimate holds true: 
\begin{equation}
\widetilde{I}_T(k)\ge\frac{k^2(e-1)}{2e(1-e^{-T})}.
\label{E:oho}
\end{equation}
\end{lemma}

\it Proof. \rm The condition in Lemma (\ref{L:use}) is equivalent to the following: $\frac{2Tk^2}{1-e^{-T}}<\frac{1}{e}$. Now, (\ref{E:ak3})
and (\ref{E:opt}) imply that 
$$
a(k)\ge\frac{1}{2}\left(\frac{2Tk^2}{1-e^{-T}}-\frac{4T^2k^4}{(1-e^{-T})^2}\right)\ge\frac{Tk^2}{1-e^{-T}}
\left(1-\frac{1}{e}\right)=\frac{Tk^2(e-1)}{e(1-e^{-T})}.
$$
Finally, (\ref{E:ak7}) implies the estimate in (\ref{E:oho}).

This completes the proof of Lemma \ref{L:use}.

Our next goal is to obtain an estimate from above for the rate function in the toy model.
\begin{corollary}\label{C:abo}
Let $T> 0$ and $k> 0$. Then, 
\begin{equation}
\widetilde{I}_T(k)\le\frac{k^2}{1-e^{-T}}.
\label{E:do1}
\end{equation}
\end{corollary}

\it Proof. \rm The estimate in (\ref{E:do1}) can be derived from (\ref{E:oioi}) by plugging the function $f$ that is identically equal to zero into the expression on the right-hand side of (\ref{E:oioi}).

\begin{theorem}\label{T:kont}
For $T> 0$ and $0< k<\sqrt{\frac{1-e^{-T}}{2Te}}$, the following two-sided estimates hold for the rate function in the toy model:
\begin{equation}
\frac{k^2(e-1)}{2e(1-e^{-T})}\le\widetilde{I}_T(k)\le\frac{k^2}{1-e^{-T}}.
\label{E:hor}
\end{equation}
\end{theorem}

Theorem \ref{T:kont} can be derived from (\ref{E:oho}) and (\ref{E:do1}).

The next assertion provides two-sided estimates by constant functions for the small-noise limit of the implied volatility in the toy model. 
\begin{theorem}\label{T:on}
For all $T> 0$ and $0< k<\sqrt{\frac{1-e^{-T}}{2Te}}$, the following inequalities are valid:
\begin{equation}
\frac{\sqrt{1-e^{-T}}}{\sqrt{2T}}\le\lim_{\varepsilon\rightarrow 0}V(\varepsilon,T,k)\le\frac{\sqrt{e(1-e^{-T})}}{\sqrt{T(e-1)}}.
\label{E:proz}
\end{equation}
\end{theorem}

\it Proof. \rm The estimates in (\ref{E:proz}) follow from Theorem \ref{T:kont} and the equality in (\ref{E:proo}).
\\
\\
This completes our analysis of the toy model.

\subsection{Asian Options}\label{SS:AA}
In this subsection, we study a small-noise asymptotic behavior of the price of an Asian option in the stochastic volatility model defined by
(\ref{E:form}). A small-noise Asian option is a path-dependent option with the payoff $(\frac{1}{T}\int_0^TS_t^{(\varepsilon)}dt-K)^{+}$. 
Here $K> 0$ is the strike, $T> 0$ is the maturity of the option, $\varepsilon\in(0,1]$ is the scaling parameter, and it is assumed that 
$K$ and $T$ are fixed. The price of the Asian option in the small-noise setting is given by 
$$
A(\varepsilon,K,T)=e^{-rT}\mathbb{E}\left[\left(\frac{1}{T}\int_0^TS_t^{(\varepsilon)}dt-K\right)^{+}\right].
$$

For every $a> 0$, define a subset of $\mathbb{C}$ by
$$
G(a,T)=\left\{f\in\mathbb{C}:\frac{1}{T}\int_0^T\exp\{f(t)\}dt> a\right\}.
$$
It is easy to see that the set $G(a,T)$ is open in $\mathbb{C}[0,T]$, and its closure in $\mathbb{C}[0,T]$ is given by
$$
\overline{G(a,T)}=\left\{f\in\mathbb{C}:\frac{1}{T}\int_0^T\exp\{f(t)\}dt\ge a\right\}.
$$

It is easy to see that for the model in (\ref{E:form}), the representation of the rate function $\widetilde{Q}_T$ given in (\ref{E:QT3}) and
(\ref{E:qeer}) is as follows:
\begin{align}
&\widetilde{Q}_T(g) 
=\frac{1}{2}\inf_{f\in\mathbb{H}_0^1}\int_0^T\left[\frac{(\dot{g}(s)-r-\rho\sigma(s,\widehat{f}(s))\dot{f}(s))^2}
{(1-\rho^2)\sigma(s,\widehat{f}(s))^2}
+\dot{f}(s)^2\right]ds
\label{E:tral}
\end{align}
for all $g\in\mathbb{H}_0^1$, and $\widetilde{Q}_T(g)=\infty$ for all $g\in\mathbb{C}\backslash\mathbb{H}^1_0$.

We will next characterize small-noise asymptotics of the price of the Asian option. 
\begin{theorem}\label{T:ee}
Suppose Assumption A and Assumptions (C1) -- (C7) hold true for the model in (\ref{E:form}) defined on the canonical set-up. Further 
suppose the condition in (\ref{E:rew}) holds and Assumption $B$ is satisfied.
Then, the following asymptotic formula is valid:
\begin{equation}
\varepsilon\log A(\varepsilon,K,T)=-\inf_{f\in\overline{G({\cal K},T)}}\widetilde{Q}_T(f)+o(1)
\label{E:strog}
\end{equation}
as $\varepsilon\rightarrow 0$ where the good rate function $\widetilde{Q}_T$ is defined in (\ref{E:QT3}) and (\ref{E:qeer}), 
and ${\cal K}$ is the moneyness 
given by ${\cal K}=\frac{K}{s_0}$.
\end{theorem}

\it Proof. \rm We will first obtain a lower large deviation estimate
for $A(\varepsilon,K,T)$. Let $\delta> 0$. It is easy to see that
\begin{align*}
&A(\varepsilon,K,T)\ge\delta e^{-rT}\mathbb{P}\left(\frac{1}{T}\int_0^TS_t^{(\varepsilon)}dt> K+\delta\right) \\
&=\delta e^{-rT}\mathbb{P}\left(\frac{1}{T}\int_0^T\exp\{X_t^{(\varepsilon)}-x_0\}dt>(K+\delta)s_0^{-1}\right) \\
&=\delta e^{-rT}\mathbb{P}\left(X^{(\varepsilon)}-x_0\in G((K+\delta)s_0^{-1},T)\right).
\end{align*}
Next, using the LDP in Theorem \ref{T:nmn},
we obtain
\begin{equation}
\liminf_{\varepsilon\rightarrow 0}\varepsilon\log A(\varepsilon,K,T)\ge-\inf_{\delta> 0}
\inf_{f\in G((K+\delta)s_0^{-1},T)}\widetilde{Q}_T(f).
\label{E:tt}
\end{equation}

Let $f\in G(Ks_0^{-1},T)$. Then, $f\in G((K+\delta)s_0^{-1},T)$ for some $\delta> 0$. Here $\delta$ depends on $f$.
Next, using (\ref{E:tt}), we see that
$$
\liminf_{\varepsilon\rightarrow 0}\varepsilon\log A(\varepsilon,K,T)\ge-
\inf_{f\in G(Ks_0^{-1},T)}\widetilde{Q}_T(f).
$$
Recall that the moneyness is defined by ${\cal K}=\frac{K}{s_0}$. Then, we have
\begin{equation}
\liminf_{\varepsilon\rightarrow 0}\varepsilon\log A(\varepsilon,K,T)\ge-
\inf_{f\in G({\cal K},T)}\widetilde{Q}_T(f).
\label{E:logm}
\end{equation}
The formula in (\ref{E:logm}) provides a lower LDP-style estimate for the price of the Asian option.

We will next obtain the corresponding upper estimate. Let $p> 1$ and $q=\frac{p}{p-1}$. It is not hard to see that
\begin{align*}
A(\varepsilon,K,T)&\le e^{-rT}T^{-\frac{1}{p}}\left\{\int_0^T\mathbb{E}\left[|S_t^{(\varepsilon)}|^p\right]dt\right\}^{\frac{1}{p}}
\mathbb{P}\left(\frac{1}{T}\int_0^TS_t^{(\varepsilon)}\ge K\right)^{\frac{1}{q}} \\
&=e^{-rT}T^{-\frac{1}{p}}\left\{\int_0^T\mathbb{E}\left[|S_t^{(\varepsilon)}|^p\right]dt\right\}^{\frac{1}{p}}
\mathbb{P}(X^{(\varepsilon)}-x_0\in\overline{G({\cal K},T)})^{\frac{1}{q}}.
\end{align*}
Therefore,
\begin{align}
\limsup_{\varepsilon\rightarrow 0}\varepsilon\log A(\varepsilon,K,T)&\le\frac{1}{p}\limsup_{\varepsilon\rightarrow 0}\varepsilon    
\log\int_0^T\mathbb{E}\left[|S_t^{(\varepsilon)}|^p\right]dt \nonumber \\
&\quad+\frac{1}{q}\limsup_{\varepsilon\rightarrow 0}\varepsilon\log
\mathbb{P}(X^{(\varepsilon)}-x_0\in\overline{G({\cal K},T)}).
\label{E:ldl}
\end{align}
Using the LDP in Theorem \ref{T:nmn} and the estimate in (\ref{E:ldl}), 
we obtain
\begin{align}
\limsup_{\varepsilon\rightarrow 0}\varepsilon\log A(\varepsilon,K,T)&\le\frac{1}{p}\limsup_{\varepsilon\rightarrow 0}\varepsilon    
\log\int_0^T\mathbb{E}\left[|S_t^{(\varepsilon)}|^p\right]dt \nonumber \\
&\quad-\frac{1}{q}\inf_{f\in\overline{G({\cal K},T)}}\widetilde{Q}_T(f).
\label{E:al1}
\end{align}

Our next goal is to estimate the first term on the right-hand side of (\ref{E:al1}). Using (\ref{E:trtt}), we get
$$
\int_0^T\mathbb{E}\left[|S_t^{(\varepsilon)}|^p\right]dt\le s_0^pe^{prT}TM_p.
$$
It follows that
\begin{equation}
\limsup_{\varepsilon\rightarrow 0}\varepsilon    
\log\int_0^T\mathbb{E}\left[|S_t^{(\varepsilon)}|^p\right]dt=0.
\label{E:pw}
\end{equation}
Now, (\ref{E:al1}) and (\ref{E:pw}) imply
\begin{equation}
\limsup_{\varepsilon\rightarrow 0}\varepsilon\log A(\varepsilon,K,T)\le-\inf_{f\in\overline{G({\cal K},T)}}\widetilde{Q}_T(f).
\label{E:ees}
\end{equation}

Our next goal is to prove that 
\begin{equation}
\inf_{f\in\overline{G({\cal K},T)}}\widetilde{Q}_T(f)=\inf_{f\in G({\cal K},T)}\widetilde{Q}_T(f).
\label{E:lds}
\end{equation}
It follows from (\ref{E:qeer}) that (\ref{E:lds}) can be rewritten as follows:
\begin{equation}
\inf_{f\in\overline{G({\cal K},T)}\cap\mathbb{H}_0^1}\widetilde{Q}_T(f)=\inf_{f\in G({\cal K},T)\cap\mathbb{H}_0^1}\widetilde{Q}_T(f).
\label{E:ty}
\end{equation}
Let us recall that since the volatility function $\sigma$ in the model in (\ref{E:mood}) satisfies the condition in (\ref{E:rew}), the rate function defined in (\ref{E:QT3}) is continuous on the space $\mathbb{H}_0^1$
(see Lemma \ref{L:dff}). It is not hard to prove that the set $G({\cal K},T)\cap\mathbb{H}_0^1$ is dense in the set 
$\overline{G({\cal K},T)}\cap\mathbb{H}_0^1$. It follows that the equalities in (\ref{E:ty}) and (\ref{E:lds}) hold. Finally, 
(\ref{E:logm}), (\ref{E:ees}), and (\ref{E:lds}) imply (\ref{E:strog}).

The proof of Theorem \ref{T:ee} is thus completed.

Our next goal is to describe what happens if we remove the condition in (\ref{E:rew}) from the formulation of Theorem \ref{T:ee}.
In such an environment, the formula in (\ref{E:tral}) may not hold for all $f\in\mathbb{H}^1_0$. However, the set of functions
for which (\ref{E:tral}) does not hold consists of at most one function. Indeed, suppose the set $B_2$ defined in Remark \ref{R:begin}
is not empty. Let $f\in B_2$. Then, the equation $\Phi(l,f,\widehat{f}(t))=g(t)$ in (\ref{E:QT1}) becomes
$rt=g(t)$. Here we use (\ref{E:disc}). It follows that the representation in (\ref{E:tral}) holds for all functions $g\in\mathbb{H}^1_0$,
except, maybe, the function $\tilde{g}(t)=rT$, $t\in[0,T]$. Moreover, by reasoning as in the proof of Lemma \ref{L:dff}, 
we can establish that the rate function in formula (\ref{E:tral}) is continuous on the set $\mathbb{H}_0^1\backslash\{\tilde{g}\}$.

\begin{theorem}\label{T:fe}
Suppose Assumption A and Assumptions (C1) -- (C7) hold true for the model in (\ref{E:form}) defined on the canonical set-up. Further 
suppose that Assumption $B$ is satisfied.
Then, for all $K>\frac{s_0}{rT}(e^{rT}-1)$, the following asymptotic formula holds:
\begin{equation}
\varepsilon\log A(\varepsilon,K,T)=-\inf_{f\in\overline{G({\cal K},T)}}\widetilde{Q}_T(f)+o(1)
\label{E:stroga}
\end{equation}
as $\varepsilon\rightarrow 0$ where the good rate function $\widetilde{Q}_T$ is defined in (\ref{E:QT3}) and (\ref{E:qeer}), 
and ${\cal K}$ is the moneyness given by ${\cal K}=\frac{K}{s_0}$.
\end{theorem}
\begin{remark}\label{R:zam}
The restriction on the strike $K$ in Theorem \ref{T:fe} is needed in order to take into account the existence of the exceptional function 
$\tilde{g}(t)=rT$. This will be explained in the proof below. If $r=0$, then this restriction becomes $K> s_0$.
\end{remark}

\it Proof of Theorem \ref{T:fe}. \rm Using the discussion before the formulation of Theorem \ref{T:fe}, we see that Theorem \ref{T:fe}
can be established exactly as Theorem \ref{T:ee}. The only difference is that the strike $K$ should be restricted so that the exceptional function $\tilde{g}$ does not belong to the set $\overline{G({\cal K},T)}$. The previous condition means that
$\frac{1}{T}\int_0^Te^{rt}dt<\frac{K}{s_0}$. The previous inequality is equivalent to $K>\frac{s_0}{rT}(e^{rT}-1)$. 

This completes the proof of Theorem \ref{T:fe}.
\subsection{Assumption $B$, revisited}\label{SS:re} In the present subsection, we provide examples of stochastic volatility models for which the asymptotic formulas for call options, the implied
volatility, and Asian options obtained in Subsections \ref{SS:call} -- \ref{SS:AA} hold true. Recall that the
results obtained in those subsections use Assumption B together with Assumption A and
Assumptions (C1) -- (C7). Our main goal in this subsection is to study scaled volatility
processes $\widehat{B}^{(\varepsilon)}$ for which Assumption B holds. However, our knowledge here is far from being complete.

In the following definition, we introduce a sub-linear growth condition for the volatility
function in (\ref{E:form}).
\begin{definition}\label{D:817}
It is said that the volatility function $\sigma:[0,T]\times\mathbb{R}^d\mapsto\mathbb{R}$ in (\ref{E:form}) satisfies a sub-linear growth
condition uniformly with respect to $t\in[0,T]$ provided that there exists a constant $C> 0$ such that 
$\sigma(t,x)^2\le C(1+||x||_d^2)$ for all $(t,x)\in[0,T]
\times\mathbb{R}^d$.
\end{definition}

The sub-linear growth condition in Definition \ref{D:817} allows us to work with a simpler inequality than that in Assumption $B$.
\begin{remark}\label{R:11}
Suppose the sub-linear growth condition holds for the volatility function $\sigma$. Suppose also that the following condition is satisfied:
For every $\beta> 0$ there exists $\varepsilon_1\in(0,1]$ depending only on $\beta$  and such that for all $0<\varepsilon<\varepsilon_1$,
the following estimate holds true:
\begin{equation}
\mathbb{E}\left[\exp\left\{\beta\int_0^T||\widehat{B}_s^{(\varepsilon)}||_d^2ds\right\}\right]\le L
\label{E:99}
\end{equation}
where $L> 0$ is a constant independent of $\varepsilon$. Then, Assumption $B$ holds. The previous statement can be easily established.
\end{remark}

\begin{remark}\label{R:LL}
It is not hard to see that if the condition in (\ref{E:99}) holds for a finite family of scaled
volatility processes, then the same condition holds for the sum of those processes. In addition, if
the condition in (\ref{E:99}) is satisfied for every component of the process $\widehat{B}^{(\varepsilon)}$, then it is satisfied
for the process $\widehat{B}^{(\varepsilon)}$.
\end{remark}
\underline{\it Example 1. Multivariate Gaussian Models.} \\
The scaled volatility process in a multivariate Gaussian stochastic volatility model is
given by $\widehat{B}^{(\varepsilon)}=(\widehat{B}^{1,(\varepsilon)},\cdots,\widehat{B}^{d,(\varepsilon)})$ where
$$
\widehat{B_t}^{i,(\varepsilon)})=x_i+\sqrt{\varepsilon}\sum_{j=1}^m\int_0^tK_{ij}(t,s)dB_s^{(j)},\quad t\in[0,T],\quad 1\le i\le d.
$$
In the previous formula, $\{K_{ij}\}$, with $1\le i\le d$ and $1\le j\le m$, is a family of admissible Volterra type
Hilbert-Schmidt kernels for which Fernique's condition is satisfied (see (\ref{E:Fer}).
\begin{lemma}\label{L:8.20}
The inequality in (\ref{E:99}) holds for a multivariate Gaussian model defined on any
set-up.
\end{lemma}

\it Proof.\rm It follows from Remark \ref{R:LL} that it suffices to prove the following: For every $\beta> 0$ there exists
$\varepsilon_1\in(0,1]$ depending only on $\beta$ and such that for all $0<\varepsilon<\varepsilon_1$, $1\le i\le d$, and
$1\le j\le m$,
$$
\mathbb{E}\left[\exp\left\{\beta\varepsilon\int_0^T|\int_0^sK_{ij}(s,u)dB_u^{(j)}|^2ds\right\}\right]\le L
$$
for some $L> 0$ independent of $\varepsilon$. The previous inequality was established in Lemma 34
in \cite{GVol}. Note that the Gaussian Volterra processes used in \cite{GVol} satisfy a stronger condition
than Assumption F used in the present paper (see Subsection \ref{SS:GM} where Assumption F is
introduced). However, the proof of Lemma 34 in\cite{GVol} does not change if we assume that
the Gaussian Volterra process satisfies Assumption F.

The proof of Lemma \ref{L:8.20} is thus completed. \\
\\
\underline{\it Example 2. Generalized Fractional Heston Models}. \\
\\
The Heston model was introduced in \cite{Hes}.
The asset price process S and the variance process V in the Heston model satisfy the following system of stochastic differential equations:
\begin{align}
&dS_t=S_t[rdt+\sqrt{V_t}(\sqrt{1-\rho^2}dW_t+\rho dB_t)] \nonumber \\
&dV_t=\kappa(\theta-V_t)dt+\eta\sqrt{V_t}dB_t.
\label{E:8.43}
\end{align}
In (\ref{E:8.43}), $W$ and $B$ are independent Brownian motions, $\kappa> 0$, $\theta > 0$, and $\eta > 0$ are positive
parameters, $r \ge 0$ is the interest rate, and $\rho\in(-1,1)$ is the correlation coefficient. The
initial conditions for the processes $S$ and $V$ are denoted by $s_0$ and $v_0$, respectively, and it
is assumed that $s_0 > 0$ and $v_0> 0$. The variance process $V$ in the Heston model is the
Cox-Ingersoll-Ross process (CIR-process). More information about the CIR-process can
be found in \cite{G}. A scaled version of the Heston model is as follows:
\begin{align}
&dS_t^{(\varepsilon)}=S_t^{(\varepsilon)}[rdt+\sqrt{\varepsilon}\sqrt{V_t^{(\varepsilon)}}(\sqrt{1-\rho^2}dW_t+\rho dB_t)] \nonumber \\
&dV_t^{(\varepsilon)}=\kappa(\theta-V_t^{(\varepsilon)})dt+\sqrt{\varepsilon}\eta\sqrt{V_t^{(\varepsilon)}}dB_t.
\label{E:8.44}
\end{align}
It is assumed in (\ref{E:8.44}) that $S_0^{(\varepsilon)}=s_0$ and $V_0^{(\varepsilon)}=v_0$ for all $\varepsilon\in(0,1]$.
\begin{remark}\label{R:8.21}
Note that for every $\varepsilon\in(0,1]$, the scaled variance process $V^{(\varepsilon)}$ is a $CIR$ process, with the same parameters 
$\kappa$ and $\theta$ as the process $V$, and with the parameter $\eta$ replaced by $\sqrt{\varepsilon}\eta$.
\end{remark}
Our next goal is to introduce generalized fractional Heston models. We will consider
such models defined on the canonical set-up. Let $K$ be a nonnegative admissible Volterra
type kernel satisfying Assumption F (see Subsection (\ref{SS:GM})). An important example of such
a kernel is the kernel of the Riemann-Liouville fractional Brownian motion. Recall that
this kernel is given by the following formula:
$K_H(t,s)=\Gamma(H+1/2)^{-1}(t-s)^{H-\frac{1}{2}}\mathbb{1}_{s< t}$ (see (\ref{E:RL}). Another interesting example is the kernel
$\widetilde{K}(t,s)=|K(t,s)|$ where $K$ is an admissible Volterra type Hilbert-Schmidt kernel satisfying Assumption F. It is not hard
to prove that the kernel $\widetilde{K}$ is also an admissible Volterra type Hilbert-Schmidt kernel satisfying Assumption F. 
A special example of such a kernel is given by 
$\widetilde{K}_H(t,s)=|K_H(t,s)|$ where $K_H$ is the Volterra type kernel of fractional Brownian motion (see Subsection \ref{SS:GM}).

Let $V$ be a $CIR$-process with parameters $\kappa$, $\theta$, and $\eta$ defined in (\ref{E:8.43}), and let $V^{(\varepsilon)}$ be a scaled version of this process defined in (\ref{E:8.44}).  Consider the following stochastic model
for the volatility process:
\begin{equation}\label{E:8.45}
\widehat{B}_t=x+\int_0^tK(t,s)V_sds,\quad x\ge 0.
\end{equation}
The scaled version of the process in (\ref{E:8.45}) is given by
\begin{equation}\label{E:8.46}
\widehat{B}_t^{(\varepsilon)}=x+\int_0^tK(t,s)V_s^{(\varepsilon)}ds,\quad x\ge 0, \quad\varepsilon\in(0,1].
\end{equation}
\begin{definition}\label{D:8.22}
The stochastic model given by
\begin{align}
&dS_t=S_t[rdt+\sqrt{V_t}(\sqrt{1-\rho^2}dW_t+\rho dB_t)] \nonumber \\
&dV_t=\kappa(\theta-V_t)dt+\eta\sqrt{V_t}dB_t. \nonumber \\
&\widehat{B}_t=x+\int_0^tK(t,s)V_sds
\label{E:8.48}
\end{align} 
will be called a generalized fractional Heston model. The scaled version of the model in (\ref{E:8.48}) is as follows:
\begin{align}
&dS_t^{(\varepsilon)}=S_t^{(\varepsilon)}[rdt+\sqrt{\widehat{B}_t^{(\varepsilon)}}(\sqrt{1-\rho^2}dW_t
+\rho dB_t)] \nonumber \\
&dV_t^{(\varepsilon)}=\kappa(\theta-V_t^{(\varepsilon)})dt+\sqrt{\varepsilon}\eta\sqrt{V_t^{(\varepsilon)}}dB_t \nonumber \\
&\widehat{B}_t^{(\varepsilon)}=x+\int_0^tK(t,s)V_s^{(\varepsilon)}ds.
\label{E:8.49}
\end{align}
\end{definition}
\begin{remark}\label{R:8.23}
If $K=K_H$ where $K_H$ is the kernel of the Riemann-Liouville fractional Brownian motion, then the model in (\ref{E:8.48}) is the fractional Heston model studied in
\cite{GJRS}.
\end{remark}
Using Theorem \ref{T:cxc} and Remark \ref{R:rrro}, we see that Assumptions (C1) -- (C7) hold true
for the volatility model introduced in (\ref{E:8.45}). Therefore, Theorems \ref{T:beg1} and \ref{T:beg11} hold for the
scaled volatility process $\widehat{B}^{(\varepsilon)}$ appearing in (\ref{E:8.46}). It follows that the scaled log-price process in
the generalized fractional Heston model satisfies the LDP in Theorem \ref{T:nmn}, provided that
the model is defined on the canonical set-up.
\begin{lemma}\label{L:8.24}
Assumption $B$ holds for the generalized fractional Heston model.
\end{lemma}

\it Proof. \rm For the sake of simplicity, we assume that $T=1$. The statement in Lemma (\ref{L:8.24}) means
the following: For every $\alpha> 0$ there exists $\varepsilon_0\in(0,1]$ depending only on $\alpha$ and such that for all $0<\varepsilon
<\varepsilon_0$, the following estimate holds true:
\begin{equation}
\mathbb{E}\left[\exp\left\{\alpha\int_0^1ds\int_0^sK(s,u)V_u^{(\varepsilon)}du\right\}\right]\le M
\label{E:88.49}
\end{equation}
where $M> 0$ is a constant depending only on $\alpha$. 

The scaled $CIR$-process $V^{(\varepsilon)}$ appearing in (\ref{E:88.49}) is given by $V_t^{(\varepsilon)}=\varepsilon V_t$, $t\in(0,1]$,
$\varepsilon\in(0,1]$
(see the last line on page 2 in \cite{GHP}). Let us denote the expression on the left-hand side of (\ref{E:88.49}) by 
$H(\varepsilon,\alpha)$. 
It is not hard to show that 
$H(\varepsilon,\alpha)\le\mathbb{E}\left[\exp\left\{\alpha C\varepsilon\max_{u\in[0.1]}V_u\right\}\right]$ where
$C=\sup_{s\in[0,1]}\int_0^sK(s,u)^2du\}^{\frac{1}{2}}$. We have $C<\infty$ (see the proof of (\ref{E:pre})). 

Next, analyzing the reasoning above, we see that in order to complete the proof of Lemma \ref{L:8.24}, it suffices to show that 
there exists $\delta> 0$ such that 
\begin{equation}
\mathbb{E}\left[\exp\left\{\delta\max_{u\in[0,1]}V_u\right\}\right]<\infty.
\label{E:8.50}
\end{equation}

We will establish the inequality in (\ref{E:8.50}) using Proposition 2.1 in \cite{D-MRYZ}.
The authors of \cite{D-MRYZ} use a scaled process $X^{(\varepsilon)}$ satisfying $dX_t^{(\varepsilon)}=b(X_t^{(\varepsilon)})dt
+2\varepsilon\sqrt{|X_t^{(\varepsilon)}|}dB_t$, $X_0^{(\varepsilon)}=a\ge 0$ for all $\varepsilon\in(0,1]$. Set $b(u)=\kappa(\theta-u)$,
$\varepsilon=\frac{\eta}{2}$, and $a=v_0$. Then, we have 
\begin{equation}
X_t^{(\frac{\eta}{2})}=V_t,\quad 0\le t\le 1.
\label{E:8.51}
\end{equation}

It was established in \cite{D-MRYZ}, Proposition 2.1 that there exists $\lambda> 0$ such that for every $\varepsilon> 0$,
\begin{equation}
\mathbb{E}\left[\exp\left\{\lambda\varepsilon^{-2}\max_{t\in[0,1]}X_t^{(\varepsilon)}\right\}\right]
\le\exp\left\{\varepsilon^{-2}+1\right\}.
\label{E:8.52}
\end{equation}

Next, using (\ref{E:8.51}) and (\ref{E:8.52}), with $\varepsilon=\frac{\eta}{2}$, we obtain
\begin{align*}
&\mathbb{E}\left[\exp\left\{\frac{4\lambda}{\eta^2}\max_{t\in[0,1]}V_t\right\}\right] 
\le\exp\left\{\frac{4}{\eta^2}+1\right\}.
\end{align*}

Therefore, the estimate in (\ref{E:8.50}) holds, with $\delta=\frac{4\lambda}{\eta^2}$.

The proof of Lemma \ref{L:8.24} is thus completed.

\section{Proof of Theorem \ref{T:beg1}}\label{S:proofs}
The proof of Theorem \ref{T:beg1} is similar in structure to that of Theorem 2.1 in \cite{CF}. However, there are certain differences between the proofs since the models that we use in the present paper are more general than 
those studied in \cite{CF}.

The next assertion (Theorem \ref{T:LP}) states that Laplace's principle holds for the process $Y^{(\varepsilon)}$. It is known that if the rate function is good, 
then the LDP in Theorem \ref{T:beg1} and Laplace's principle in Theorem \ref{T:LP} are equivalent. Since the rate function $I_y$ 
is good (see Remark \ref{R:rate}), Theorem \ref{T:beg1}
can be derived from Theorem \ref{T:LP}.
\begin{theorem}\label{T:LP}
Suppose the conditions in Theorem \ref{T:beg1} hold. Then, for all bounded and continuous functions $F:{\cal W}^d\mapsto\mathbb{R}$,
$$
\lim_{\varepsilon\rightarrow 0}-\varepsilon\log\mathbb{E}[\exp\{-\frac{1}{\varepsilon}F(Y^{(\varepsilon)})\}]
=\inf_{\phi\in{\cal W}^d}[I_y(\phi)+F(\phi)]
$$
where the rate function $I_y$ is defined by (\ref{E:LDP})
\end{theorem}

\it Proof. Lower bound in Laplace's principle. \rm We have to show that
\begin{equation}
\lim_{\varepsilon\rightarrow 0}-\varepsilon\log\mathbb{E}\left[\exp\{-\frac{1}{\varepsilon}F(Y^{(\varepsilon)})\}\right]
\ge\inf_{\phi\in{\cal W}^d}[I_y(\phi)+F(\phi)].
\label{E:lb1}
\end{equation} 
It suffices to prove that for any sequence $\varepsilon_n\in (0,1]$, $n\ge 1$, such that $\varepsilon_n\rightarrow 0$ as $n\rightarrow\infty$,
there exists a subsequence along which the inequality in (\ref{E:lb1}) holds. Let $Y^{(\varepsilon_n)}$, with $n\ge 1$, be the strong solution 
to the equation in (\ref{E:2ooo}) with $\varepsilon=\varepsilon_n$ (see Assumption (C3)). Then, there exists a map 
$h^{(n)}:{\cal W}^m\mapsto{\cal W}^d$ such that it is $\widetilde{{\cal B}}_t^m/{\cal B}_t^d$-measurable for all $t\in[0,T]$, and, moreover,
\begin{equation}
Y^{(\varepsilon_n)}=h^{(n)}(B)\quad\mathbb{P}-\mbox{a.s.} 
\label{E:oig}
\end{equation}
(see Remark \ref{R:oneof}). 

We will use a variational representation of functionals of Brownian motion (see, e.g., Theorem 3.6 in the paper \cite{BDa} of Budhiraja and Dupuis).
\begin{theorem}\label{T:vr}
Let $f$ be a bounded Borel measurable real function on ${\cal W}^m$. Then,
$$
-\log\mathbb{E}[\exp\{-f(B)\}]=\inf_{v\in{\cal M}^2[0,T]}\mathbb{E}\left[\frac{1}{2}\int_0^T||v_s||_m^2ds+f\left(B+\int_0^{\cdot}v_sds\right)
\right].
$$
\end{theorem}
For the function $F:{\cal W}^d\mapsto\mathbb{R}$ appearing in the formulation of Theorem \ref{T:LP}, the representation 
formula in Theorem \ref{T:vr} 
implies the following:
\begin{align}
&-\varepsilon_n\log\mathbb{E}\left[\exp\{-\frac{1}{\varepsilon_n}F(Y^{(\varepsilon_n)})\}\right]=
-\varepsilon_n\log\mathbb{E}\left[\exp\{-\frac{1}{\varepsilon_n}F\circ h^{(n)}(B)\}\right] \nonumber \\
&=\varepsilon_n\inf_{v\in{\cal M}^2[0,T]}\mathbb{E}\left[\frac{1}{2}\int_0^T||v_s||_m^2ds+\frac{1}{\varepsilon_n}
F\circ h^{(n)}\left(B+\int_0^{\cdot}v_sds\right)\right] \nonumber \\
&=\inf_{v\in{\cal M}^2[0,T]}\mathbb{E}\left[\frac{1}{2}\int_0^T||v_s||_m^2ds+
F\circ h^{(n)}\left(B+\frac{1}{\sqrt{\varepsilon_n}}\int_0^{\cdot}v_sds\right)\right].
\label{E:b1}
\end{align}
The last equality in (\ref{E:b1}) is obtained by passing from $v$ to $\frac{v}{\sqrt{\varepsilon_n}}$.

Fix $\delta> 0$. It can be shown exactly as in \cite{CF} 
that there exists $N> 0$ such that for every $n\ge 1$, a control $v^{(n)}$ can be found satisfying
\begin{equation}
v^{(n)}\in M^2_N[0,T]
\label{E:b3}
\end{equation}
and
\begin{align}
&-\varepsilon_n\log\mathbb{E}\left[\exp\{-\frac{1}{\varepsilon_n}F(Y^{(\varepsilon_n)})\}\right] \nonumber \\
&\ge\mathbb{E}\left[\frac{1}{2}\int_0^T||v_s^{(n)}||_m^2ds
F\circ h^{(n)}\left(B+\frac{1}{\sqrt{\varepsilon_n}}\int_0^{\cdot}v_s^{(n)}ds\right)\right]-\delta.
\label{E:b2}
\end{align} 
The controls $v^{(n)}$ in (\ref{E:b2}) depend on $\delta$. 
Let $N$ and $v_n$ be such as in (\ref{E:b3}) and (\ref{E:b2}).
For every $n\ge 1$, consider the process $B^{\varepsilon_n, v^n}$ defined by (\ref{E:brow}). For the sake of shortness, we will use the symbol
$B^{(n)}$ instead of $B^{\varepsilon_n, v^n}$. The process $B^{(n)}$ is an $m$-dimensional Brownian motion on ${\cal W}^m$ with respect to a measure $\mathbb{P}^{(n)}$ on ${\cal B}_T^m$ that is equivalent to the 
measure $\mathbb{P}$. This process is given by 
$B^{(n)}_s=B_s+\frac{1}{\sqrt{\varepsilon_n}}\int_0^sv_u^{(n)}du$, $s\in[0,T]$.
The process $B^{(n)}$ is adapted to the filtration $\{{\cal B}_t^m\}$.

Fix $n\ge 1$, and consider the following scaled controlled stochastic integral equation:
\begin{align}
Y_t^{\varepsilon_n,v^{(n)}}&=y+\int_0^ta(t,s,V^{1,\varepsilon_n,v^{(n)}},Y^{\varepsilon_n,v^{(n)}})ds
+\int_0^tc(t,s,V^{2,\varepsilon_n,v^{(n)}},
Y^{\varepsilon_n,v^{(n)}})v_s^{(n)}ds
\nonumber \\
&\quad+\sqrt{\varepsilon_n}\int_0^tc(t,s,V^{2,\varepsilon_n,v^{(n)}},Y^{\varepsilon_n,v^{(n)}})dB_s
\label{E:rt}
\end{align}
where the processes $V^{i,\varepsilon_n,v^{(n)}}$, with $i=1,2$, satisfy 
\begin{align}
V_s^{i,\varepsilon_n,v^{(n)}}&=V_0^{(i)}+\int_0^s\bar{b}_i(r,V^{i,\varepsilon_n,v^{(n)}})dr
+\int_0^s\bar{\sigma}_i(r,V^{i,\varepsilon_n,v^{(n)}})v_r^{(n)}dr \nonumber \\
&\quad+\sqrt{\varepsilon_n}\int_0^s\bar{\sigma}_i(r,V^{i,\varepsilon_n,v^{(n)}})dB_r.
\label{E:b5}
\end{align}
We will next show that for every $n\ge 1$, the equation in (\ref{E:rt}) possesses 
a strong solution. This solution is unique, by 
Assumption (C3)(b). We will also provide a representation formula for the unique solution. 
\begin{lemma}\label{L:uni}
Let $v\in{\cal M}^2[0,T]$ be such that (\ref{E:cve}) holds, and suppose Assumptions (C1) - (C3) and (C6) are satisfied. Then, 
the equation in (\ref{E:3ooo}) has 
a strong solution $Y^{(v)}$ with $Y_0^{(v)}=y$. Moreover, the following formula holds 
$\mathbb{P}$-a.s.: $Y^{(v)}=h\left(B^{(v)}\right)$ where $h$ and $B^{(v)}$
are defined in Remark \ref{R:ghg} and (\ref{E:brew}), respectively.
\end{lemma}
\begin{remark}\label{R:unk}
Lemma \ref{L:uni} is similar to Lemma A.1 in \cite{CF}.
\end{remark}
\begin{remark}\label{R:rr}
The existence of the strong solution $Y^{\varepsilon_n,v^{(n)}}$ to the equation in (\ref{E:rt}) and the representation formula
\begin{equation}
Y^{\varepsilon_n,v^{(n)}}=h^{(n)}(B^{\varepsilon_n,v^{(n)}}),\quad n\ge 1
\label{E:rem}
\end{equation}
where $h^{(n)}$ and $B^{\varepsilon,v}$ are defined in (\ref{E:oig}) and (\ref{E:brow}), respectively,
can be established as follows. Fix $n\ge 1$, and replace the control $v$ by $\frac{v^{(n)}}{\sqrt{\varepsilon_n}}$, the maps 
$\bar{\sigma}_i$
and $c$ by $\sqrt{\varepsilon_n}\bar{\sigma}_i$ and $\sqrt{\varepsilon_n}c$, respectively, and also replace $h$ by $h^{(n)}$.
After making such replacements, we can apply Lemma \ref{L:uni} to establish the formula in (\ref{E:rem}). 
\end{remark}

\it Proof of Lemma \ref{L:uni}. \rm Using (\ref{E:2a}) and (\ref{E:2b}), we obtain
\begin{equation}
h(B)=y+\int_0^{\cdot}a(\cdot,s,V^{(1)},h(B))ds
+\int_0^{\cdot}c(\cdot,s,V^{(2)},h(B))dB_s.
\label{E:op}
\end{equation}
In (\ref{E:op}), $V^{(i)}$ are the unique solutions to the equations
\begin{equation}
V_s^{(i)}=V_0^{(i)}+\int_0^s\bar{b}_i(r,V^{(i)})dr+\int_0^s\bar{\sigma}_i(r,V^{(i)})dB_r,\quad i=1,2.
\label{E:777}
\end{equation}
Since the conditions in Theorem 10.4 in 
\cite{RW} are satisfied for the equations in (\ref{E:777}), there exist maps $g^{(i)}:{\cal W}^m\mapsto{\cal W}^{k_i}$, with $i=1,2$, such that: 
(i) The map $g^{(i)}$ is
$\widetilde{{\cal B}}_t^m/{\cal B}_t^{k_i}$-measurable for every $t\in[0,T]$; (ii) $V^{(i)}=g^{(i)}(B)$ for $i=1,2$; (iii) For every $i=1,2$,
\begin{align}
&g^{(i)}(B^{(v)})=V_0^{(i)}+\int_0^s\bar{b}_i(r,g^{(i)}(B^{(v)}))dr+\int_0^s\bar{\sigma}_i(r,g^{(i)}(B^{(v)}))dB^{(v)}_s \nonumber \\
&=V_0^{(i)}+\int_0^s\bar{b}_i(r,g{(i)}(B^{(v)}))dr+\int_0^s\bar{\sigma}_i(r,g^{(i)}(B^{(v)}))v_sds+\int_0^s\bar{\sigma}_i(r,g^{(i)}
(B^{(v)}))dB_s.
\label{E:murr}
\end{align}
Next, using (\ref{E:murr}), we see that for every $i=1,2$, the equality $V^{i,v}=g^{(i)}(B^{(v)})$ holds.
Therefore, the equation in (\ref{E:op}) can be rewritten as follows:
\begin{equation}
h(B)=y+\int_0^{\cdot}a(\cdot,s,g^{(1)}(B),h(B))ds
+\int_0^{\cdot}c(\cdot,s,g^{(2)}(B),h(B))dB_s.
\label{E:opp}
\end{equation}
It is not hard to see that if we could replace $B$ by $B^{(v)}$ in (\ref{E:opp}), then 
Lemma \ref{L:uni} will be established. Here we take into account Assumption (C3)(b). There is no problem in replacing $B$ by $B^{(v)}$ in the expression on the left-hand side of 
(\ref{E:opp}) and in the first integral on the right-hand side of (\ref{E:opp}). We will next explain how to deal with the stochastic integral appearing in (\ref{E:opp}). By the second inequality in (\ref{E:c2}),
\begin{equation}
\int_0^tc_{lj}(t,s,g^{(2)}(B),h(B))^2ds<\infty\quad\mathbb{P}-\mbox{a.s.} 
\label{E:rrc}
\end{equation}
for all $t\in[0,T]$, $1\le l\le d$, and $1\le j\le m$. It is clear that the paths of 
the process $h(B)$ are continuous. Therefore, (\ref{E:opp}) and the first restriction on the function $a$ in Assumption (C2)(b) imply
that the paths of the process 
\begin{equation}
t\mapsto\int_0^tc(t,s,g^{(2)}(B),h(B))dB_s
\label{E:also}
\end{equation} 
are continuous as well. 

Using Lemma 10.1 in \cite{RW}, we can show that for every rational number $r\in[0,T]$ there exists a functional 
$u_r:{\cal W}^m\mapsto{\cal W}^d$ satisfying the following conditions: \\
(i)\,For every $r$, the functional $u_r$ 
is $\widetilde{{\cal B}}_t^m/{\cal B}_t^d$-measurable for all $t\in[0,T]$; \\
(ii)\,$\mathbb{P}$-a.s. on the space ${\cal W}^m$, the equalities
$$
u_r(B)_{\cdot}=\int_0^{\cdot}c(r,s,g^{(2)}(B),h(B))\mathbb{1}_{\{s\le r\}}dB_s
$$
and 
$$
u_r(B^{(v)})_{\cdot}
=\int_0^{\cdot}c(r,s,g^{(2)}(B^{(v)}),h(B^{(v)}))
\mathbb{1}_{\{s\le r\}}dB^{(v)}_s
$$
hold for all $r$. The existence of the previous stochastic integrals follows from the second inequality in (\ref{E:c2}).

Using (ii), we see that $\mathbb{P}$-a.s. on ${\cal W}^m$, $u_r(B)(r)=\int_0^rc(r,s,g^{(2)}(B),h(B))dB_s$ and, in addition,
$
u_r(B^{(v)})(r)=\int_0^rc(r,s,g^{(2)}(B^{(v)}),h(B^{(v)}))dB^{(v)}_s
$
for all $r\in[0,T]$. It follows from the continuity of the process in (\ref{E:also}) and Assumption (C6) that the 
following statement holds true
$\mathbb{P}$-a.s. on ${\cal W}^m$:
For every $t\in[0,T]$ and every sequence $r_i$ of rational numbers in $[0,T]$ such that $r_i\rightarrow t$, 
$$
\lim_{r_i\rightarrow t}u_{r_i}(B)(r_i)=\int_0^tc(t,s,g^{(2)}(B),h(B))dB^{(v)}_s
$$
and
\begin{equation}
\lim_{r_i\rightarrow t}u_{r_i}(B^{(v)})(r_i)=\int_0^tc(t,s,g^{(2)}(B^{(v)}),h(B^{(v)}))dB^{(v)}_s.
\label{E:prof}
\end{equation}
Therefore, (\ref{E:opp}) implies that
\begin{equation}
h(B)=y+\int_0^{\cdot}a(\cdot,s,g^{(1)}(B),h(B))ds+u(B)\quad\mathbb{P}-\rm a.s.
\label{E:oppo}
\end{equation}
where $u(B)=\lim_{r_i\rightarrow t}u_{r_i}(B)(r_i)$. It follows from (\ref{E:oppo}) that
$u:{\cal W}^m\mapsto{\cal W}^d$ is an $\widetilde{{\cal B}}_t^m/{\cal B}_t^d$-measurable functional for every $t\in[0,T]$. Since $B$ is the coordinate process on 
${\cal W}^m$, the mapping $s\mapsto B_s$ is the identity mapping on ${\cal W}^m$. Therefore,  
$u(\eta)(t)=\lim_{r_i\rightarrow t}u_{r_i}(\eta)(r_i)$ for all $t\in[0,T]$\,\,$\mathbb{P}$-a.s.
It follows that
$u(B^{(v)})=\lim_{r_i\rightarrow t}u_{r_i}(B^{(v)})(r_i)$\,\,
$\mathbb{P}$-a.s. The equality in (\ref{E:oppo}) is a statement about some measurable functional of $B$ and the measure $\mathbb{P}$. Since the distribution of $B$ with respect to $\mathbb{P}$ is the same as the distribution of $B^{(v)}$ with respect to $\mathbb{P}^v$,
we can replace $B$ by $B^{(v)}$ and $\mathbb{P}$ by $\mathbb{P}^{(v)}$ in (\ref{E:oppo}). Next, using (\ref{E:prof}) we see that
\begin{equation}
h(B^{(v)})=y+\int_0^{\cdot}a(\cdot,s,g^{(1)}(B^{(v)}),h(B^{(v)}))ds+\int_0^tc(t,s,g^{(2)}(B^{(v)}),h(B^{(v)}))dB^{(v)}_s.
\label{E:opposum}
\end{equation}

Finally, by comparing (\ref{E:3ooo}) and (\ref{E:opposum}) and using Assumption (C3)(b), we see
that the following equality holds: $Y^{(v)}=h\left(B^{(v)}\right)$.

This completes the proof of Lemma \ref{L:uni}.

We will next return to the proof of Theorem \ref{T:LP}. 
Using formula \ref{E:rem} we can rewrite the estimate in (\ref{E:b2})
as follows:
\begin{align}
&-\varepsilon_n\log\mathbb{E}\left[\exp\{-\frac{1}{\varepsilon_n}F(Y^{(\varepsilon_n)})\}\right]
\ge\mathbb{E}\left[\frac{1}{2}\int_0^T||v_s^{(n)}||_m^2ds+
F(Y^{\varepsilon_n,v^{(n)}})\right]-\delta.
\label{E:b21}
\end{align} 
Next, reasoning as in \cite{CF} we can show that the sequence $(Y^{\varepsilon_n,v^{(n)}},V^{1,\varepsilon_n,v^{(n)}},
V^{2,\varepsilon_n,v^{(n)}},v^{(n)})$, 
$n\ge 1$, is tight as a family of random variables with values in 
${\cal W}^d\times{\cal W}^{k_1}\times{\cal W}^{k_2}\times D_N$ for some $N> 0$ where $D_N$ is defined in Assumption (C5). The tightness of 
$\{Y^{\varepsilon_n,v^{(n)}}\}$ follows from Assumption (C7), while the tightness of $\{V^{1,\varepsilon_n,v^{(n)}}\}$ and 
$\{V^{2,\varepsilon_n,v^{(n)}}\}$ follows from (H6) in \cite{CF}.
Hence, possibly taking a subsequence, we see that the sequence of random variables $(Y^{\varepsilon_n,v^{(n)}},V^{1,\varepsilon_n,v^{(n)}},
V^{2,\varepsilon_n,v^{(n)}},v^{(n)})$ converges in distribution as $n\rightarrow\infty$ to a $({\cal W}^d\times{\cal W}^{k_1}\times{\cal W}^{k_2}\times D_N)$-valued random variable 
$(\widehat{Y},\widehat{V}^{(1)},\widehat{V}^{(2)},v)$ that is defined on some probability space 
$(\widehat{\Omega},\widehat{{\cal F}},\widehat{\mathbb{P}})$.
\begin{lemma}\label{L:sat}
The processes $\widehat{V}^{(1)}$, $\widehat{V}^{(2)}$, and $\widehat{Y}$ satisfy the following system of equations:
\begin{align}
&\widehat{V}^{(i)}_s=V_0^{(i)}+\int_0^s\bar{b}_i(r,\widehat{V}^{(i)})dr+\int_0^s\bar{\sigma}_i(r,\widehat{V}^{(i)})v_rdr,\quad i=1,2, \nonumber \\
&\widehat{Y}_t=y+\int_0^ta(t,s,\widehat{V}^{(1)},\widehat{Y})ds+\int_0^tc(t,s,\widehat{V}^{(2)},\widehat{Y})v_sds
\label{E:rtr}
\end{align}
$\widehat{P}$-a.s. on $\widehat{\Omega}$.
\end{lemma}

\it Proof. \rm In the proof of Lemma \ref{L:sat}, we borrow some ideas from the proof on p. 1132 in \cite{CF}. However, there are 
significant differences between those proofs because of the difficulties caused by the extra processes $V^{(1)}$ and $V^{(2)}$. 

For $\widehat{V}^{(i)}$, with $i=1,2$, the validity of the statement in Lemma \ref{L:sat} was established in \cite{CF} (see the proof of formula (12) in \cite{CF}). We will next prove 
the same statement for $\widehat{Y}$. 
Let $t\in[0,T]$, and consider the map
$\Psi_t:{\cal W}^d\times{\cal W}^{k_1}\times{\cal W}^{k_2}\times D_N\mapsto\mathbb{R}$ defined by
$$
\Psi_t(\varphi,\tau_1,\tau_2,f)=||\varphi(t)-y-\int_0^ta(t,s,\tau_1,\varphi)ds-\int_0^tc(t,s,\tau_2,\varphi)f(s)ds||_d\wedge 1.
$$
It is clear that this map is bounded. Our next goal is to show that it is continuous. Let $\varphi_n\rightarrow\varphi$ in ${\cal W}^d$,  
$\tau^{(n)}_1\rightarrow\tau_1$ in ${\cal W}^{k_1}$, $\tau^{(n)}_2\rightarrow\tau_2$ in ${\cal W}^{k_2}$. Suppose also that $f_n\in D_N$,
$f\in D_N$, and
$f_n\rightarrow f$ in the weak topology of $L^2([0,T],\mathbb{R}^m)$ (see Assumption (C5) for the definition of $D_N$. We have
\begin{align}
|\Psi_t(\varphi_n,\tau^{(n)}_1,\tau^{(n)}_2,f_n)-\Psi_t(\varphi,\tau_1,\tau_2,f)|&\le||\varphi_n(t)-\varphi(t)||_d \nonumber \\
&+||\int_0^ta(t,s,\tau^{(n)}_1,\varphi_n)ds-\int_0^ta(t,s,\tau_1,\varphi)ds||_d \nonumber \\
&+\int_0^t||c(t,s,\tau^{(n)}_2,\varphi_n)-c(t,s,\tau_2,\varphi)||_{d\times m}||f_n(s)||_mds \nonumber \\
&+||\int_0^tc(t,s,\tau_2,\varphi)(f_n(s)-f(s))ds||_d. 
\label{E:ali1}
\end{align}
The first term on the right-hand side of (\ref{E:ali1}) tends to zero as $n\rightarrow\infty$ because 
$\varphi_n\rightarrow\varphi$ in ${\cal W}^d$. To prove the same for the second term, we use Assumption (C2)(b) for the vector function
 $a$.
The third term can be handled by 
using H\"{o}lder's inequality, the boundedness of the family $f_n$ in the space $L^2([0,T],\mathbb{R}^m)$, and Assumption (C2)(c) for the matrix function $c$. Finally, the fourth term on the right-hand side of (\ref{E:al1}) tends to zero as $n\rightarrow\infty$ by the $L^2$-condition for the matrix function $c$ in Assumption (C2)(a) and because $f_n\rightarrow f$ weakly in 
$L^2([0,T],\mathbb{R}^m)$. This establishes the continuity of the map $\Psi_t$.
Next, we use the continuous mapping theorem for the weak convergence to show that
\begin{equation}
\lim_{n\rightarrow\infty}\mathbb{E}[\Psi_t(Y^{\varepsilon_n,v^n},V^{1,\varepsilon_n,v^n},V^{2,\varepsilon_n,v^n},v^n)]
=\mathbb{E}_{\widehat{\mathbb{P}}}[\Psi_t(\widehat{Y},\widehat{V}^1,\widehat{V}^2,v)].
\label{E:kut}
\end{equation}
It follows from the definition of $\Psi_t$ that
\begin{align*}
&\mathbb{E}[\Psi_t(Y^{\varepsilon_n,v^n},V^{1,\varepsilon_n,v^n},V^{2,\varepsilon_n,v^n},v^n)] \\
&\le\mathbb{E}\left[||Y^{\varepsilon_n,v^n}_t-y-\int_0^ta(t,s,V^{1,\varepsilon_n,v^n},Y^{\varepsilon_n,v^n})ds
-\int_0^tc(t,s,V^{2,\varepsilon_n,v^n},Y^{\varepsilon_n,v^n})v^n_sds||_d\right].
\end{align*}
Now, using (\ref{E:3ooo}) we obtain
\begin{align}
\mathbb{E}[\Psi_t(Y^{\varepsilon_n,v^n},V^{1,\varepsilon_n,v^n},V^{2,\varepsilon_n,v^n},v^n)]
&\le\sqrt{\varepsilon_n}\mathbb{E}\left[||\int_0^tc(t,s,V^{2,\varepsilon_n,v^n},Y^{\varepsilon_n,v^n})dB_s||_d\right] \nonumber \\
&\le\sqrt{\varepsilon_n}\sqrt{\int_0^t\mathbb{E}[||c(t,s,V^{2,\varepsilon_n,v^n},Y^{\varepsilon_n,v^n})||_d^2]ds}.
\label{E:gisht}
\end{align}
For every $t\in[0,T]$, the last expression in (\ref{E:gisht}) tends to zero as $n\rightarrow\infty$ by
(\ref{E:b3}) and Assumption (C7). Therefore, (\ref{E:kut}) implies that for every $t\in[0,T]$, $\Psi_t(\widehat{Y},\widehat{V},v)=0$\,\,
$\widehat{\mathbb{P}}$-a.s., and hence for every $t\in[0,T]$, $\widehat{Y}$ satisfies the second equation in (\ref{E:rtr})\,\,
$\widehat{\mathbb{P}}$-a.s. Since $\widehat{Y}$ maps $\widehat{\Omega}$ into ${\cal W}^d$, the equation in (\ref{E:rtr}) holds 
for all $t\in[0,T]$\,\,$\widehat{\mathbb{P}}$-a.s.

This completes the proof of Lemma \ref{L:sat}.

We are finally ready to finish the proof of the lower bound in Laplace's principle in Theorem \ref{T:LP}. We will follow the proof at the bottom of p. 1133 in \cite{CF}. Using Fatou's lemma for the convergence in distribution,
we obtain 
\begin{equation}
\liminf_{n\rightarrow\infty}\mathbb{E}[\int_0^T||v_s^n||_m^2ds]\ge\mathbb{E}_{\widehat{\mathbb{P}}}[\int_0^T||v_s||_m^2ds].
\label{E:llp}
\end{equation}
It follows from Lemma \ref{L:sat} and (\ref{E:b21}) that
\begin{equation}
\liminf_{n\rightarrow\infty}-\varepsilon_n\log\mathbb{E}\left[\exp\{-\frac{1}{\varepsilon_n}F(Y^{\varepsilon_n})\}\right]
\ge\liminf_{n\rightarrow\infty}\mathbb{E}\left[\frac{1}{2}\int_0^T||v_s^n||_m^2ds+
F(Y^{\varepsilon_n,v^n})\right]-\delta.
\label{E:bb2}
\end{equation} 
Since $F$ is continuous and bounded, we can use (\ref{E:bb2}), the continuous mapping theorem, and (\ref{E:llp}) to get
\begin{align}
&\liminf_{n\rightarrow\infty}-\varepsilon_n\log\mathbb{E}\left[\exp\{-\frac{1}{\varepsilon_n}F(Y^{\varepsilon_n})\}\right]
\ge\mathbb{E}_{\widehat{\mathbb{P}}}\left[\frac{1}{2}\int_0^T||v_s||_m^2ds+
F(\widehat{Y})\right]-\delta.
\label{E:tryy}
\end{align}

Let us recall that under the restrictions imposed on the functions $\bar{b}_i$ and $\bar{\sigma}_i$ in \cite{CF}, the functional equations
\begin{equation}
\psi_i(s)=V_0^i+\int_0^s\bar{b}_i(r,\psi_i)dr+\int_0^s\bar{\sigma}_i(r,\psi_i)f(r)dr,\quad i=1,2,
\label{E:psik}
\end{equation}
are uniquely solvable. Moreover, for every $i=1,2$, the solution $\psi_{i,f}$ belongs to the space ${\cal W}^{k_i}$, 
and if $f_n\rightarrow f$ weakly in 
$L^2([0,T],\mathbb{R}^m)$, then $\psi_{i,f_n}\rightarrow\psi_{i,f}$ in ${\cal W}^{k_i}$. Therefore, the solutions to the equations in 
(\ref{E:rty}) are deterministic, and they coincide with the functions $\psi_{i,f}$ (see Remark \ref{R:redi}). 

Let $v$ be the control process appearing in Lemma \ref{L:sat}. Then, $\widehat{\mathbb{P}}$-almost all paths of $v$ belong to the space 
$L^2([0,T],\mathbb{R}^m)$. For every such path $f$, (\ref{E:psik}) and the first equation in (\ref{E:rtr}) show that 
$\widehat{V}^i=\psi_{i,v}$\,\,$\widehat{\mathbb{P}}$-a.s. Moreover, the previous equality and the second equation in (\ref{E:rtr}) imply
that $\widehat{Y}=\Gamma_y(v)$\,\,$\widehat{\mathbb{P}}$-a.s. Now, we can estimate the expectation on the right-hand side of 
(\ref{E:tryy}) by the essential greatest lower bound of the integrand and use the equality $\widehat{Y}=\Gamma_y(v)$. This gives
\begin{align*}
&\liminf_{n\rightarrow\infty}-\varepsilon_n\log\mathbb{E}\left[\exp\{-\frac{1}{\varepsilon_n}F(Y^{\varepsilon_n})\}\right] \\
&\ge\inf_{\{(\varphi,f)\in{\cal W}^d\times L^2:\varphi=\Gamma_y(f)\}}\left\{\frac{1}{2}\int_0^T||f(s)||_m^2ds+F(\varphi)\right\}-\delta
\ge\inf_{\varphi\in{\cal W}^d}\{I_y(\varphi)+F(\varphi)\}-\delta
\end{align*}
where $I_y$ is the rate function defined in (\ref{E:LDP}). 

This completes the proof of the lower bound in the Laplace's principle in Theorem \ref{T:LP}, since the previous estimate holds for an arbitrary $\delta> 0$. 
\begin{remark}\label{R:tyu}
It can be shown that the equation in (\ref{E:t2}) is solvable for any deterministic control 
$f\in L^2([0,T],\mathbb{R}^m)$ by using the same ideas as above. Therefore, only the uniqueness condition should be included in 
Assumption (C4). A similar remark can be found on p. 1133 in \cite{CF}.
\end{remark}

\it Upper bound in Laplace's principle. \rm We have to show that
\begin{equation}
\limsup_{\varepsilon\rightarrow 0}-\varepsilon\log\mathbb{E}\left[\exp\{-\frac{1}{\varepsilon}F(Y^{\varepsilon})\}\right]
\le\inf_{\varphi\in{\cal W}^d}[I_y(\varphi)+F(\varphi)]
\label{E:lbc1}
\end{equation}
for all bounded and continuous functions $F:{\cal W}^d\mapsto\mathbb{R}$. Exactly as in the proof 
on p. 1134 in \cite{CF}, we assume that the greatest lower bound in (\ref{E:lbc1}) is finite and prove that for every fixed $\delta> 0$ 
there exists $\varphi\in{\cal W}^d$ for which
\begin{equation}
I_y(\varphi)+F(\varphi)\le\inf_{\zeta\in{\cal W}^d}(I_y(\zeta)+F(\zeta))+\frac{\delta}{2}<\infty.
\label{E:odd}
\end{equation}
We can also choose, for the function $\varphi$ defined above, a control $f\in L^2([0,T],\mathbb{R}^m)$ such that
$$
\frac{1}{2}\int_0^T||f(s)||_m^2ds\le I_y(\varphi)+\frac{\delta}{2}
$$
and $\varphi=\Gamma_y(f)$. It follows from (\ref{E:odd}) that
\begin{equation}
I_y(\Gamma_y(f))+F(\Gamma_y(f))\le\inf_{\zeta\in{\cal W}^d}(I_y(\zeta)+F(\zeta))+\frac{\delta}{2}<\infty.
\label{E:oddi}
\end{equation}

Let us suppose that $\varepsilon_n\in(0,1]$ and 
$\varepsilon_n\rightarrow 0$. For every $n\ge 1$, denote by $Y^{\varepsilon_n,f}$ the unique strong solution to the equation in (\ref{E:2ooo}).
The existence of the solution follows from Lemma \ref{L:uni}. Next, reasoning exactly as in the proof of the lower bound, we establish that the family 
$(Y^{\varepsilon_n,f},V^{1,\varepsilon_n,f},V^{2,\varepsilon_n,f},f)$, with $n\ge 1$, 
is tight in ${\cal W}^d\times{\cal W}^{k_1}\times{\cal W}^{k_2}\times D_N$ for some number $N> 0$ where $D_N$ is defined in 
Assumption (C5). Hence, possibly taking a subsequence, $(Y^{\varepsilon_n,f},
V^{1,\varepsilon_n,f},V^{2,\varepsilon_n,f},f)$ converges in distribution to a 
${\cal W}^d\times{\cal W}^{k_1}\times{\cal W}^{k_2}\times D_N$-valued random variable 
$(\widehat{Y},\widehat{V}^1,\widehat{V}^2,f)$ that is defined on some probability space 
$(\widehat{\Omega},\widehat{{\cal F}},\widehat{\mathbb{P}})$. 
We can also prove that $\widehat{\mathbb{P}}$-a.s. on $\widehat{\Omega}$ we have
\begin{align*}
&\widehat{V}_s^i=V_0^i+\int_0^s\bar{b}_i(r,\widehat{V}^i)dr+\int_0^s\bar{\sigma}_i(r,\widehat{V}^i)f(r)dr 
\end{align*}
and
\begin{align}
&\widehat{Y}_t=y+\int_0^ta(t,s,\widehat{V}^1,\widehat{Y})ds+\int_0^tc(t,s,\widehat{V}^2,\widehat{Y})f(s)ds
\label{E:eee}
\end{align}
for all $t\in[0,T]$. The solutions to the equations for $\widehat{V}^i$ in (\ref{E:eee}), with $i=1,2$, are deterministic, 
and we have $\widehat{V}^i=\psi_{i,f}$\,\, 
$\widehat{\mathbb{P}}$-a.s. on $\widehat{\Omega}$. It follows that the second equation in (\ref{E:eee}) can be rewritten as follows::
$$
\widehat{Y}_t=y+\int_0^ta(t,s,\psi_{1,f},\widehat{Y})ds+\int_0^tc(t,s,\psi_{2,f},\widehat{Y})f(s)ds.
$$
Therefore, $\widehat{Y}=\Gamma_y(f)$\,\,$\widehat{\mathbb{P}}$-a.s. on $\widehat{\Omega}$. Next, using (\ref{E:b1}) 
and Remark \ref{R:rr} and reasoning as at the end of the proof of the upper bound in Laplace's principle in \cite{CF}
(see p. 1134 in \cite{CF}), we see that
\begin{equation}
\limsup_{\varepsilon\rightarrow 0}-\varepsilon\log\mathbb{E}\left[\exp\{-\frac{1}{\varepsilon}F(Y^{\varepsilon})\}\right]
\le I_y(\Gamma_y(f))+\frac{\delta}{2}+\lim_{n\rightarrow\infty}\mathbb{E}[F(Y^{\varepsilon_n,f})].
\label{E:rtu}
\end{equation}
Recall that the sequence $Y^{\varepsilon_n,f}$ converges in distribution to $\Gamma_y(f)$. Now, using the continuity theorem for the convergence in distribution and (\ref{E:oddi}), we obtain from (\ref{E:rtu}) that the following estimate holds:
\begin{align*}
\limsup_{\varepsilon\rightarrow 0}-\varepsilon\log\mathbb{E}\left[\exp\{-\frac{1}{\varepsilon}F(Y^{\varepsilon})\}\right] 
&\le I_y(\Gamma_y(f))+\frac{\delta}{2}+F(\Gamma_y(f)]) \\
&\le\inf_{\zeta\in{\cal W}^d}(I_y(\zeta)+F(\zeta))+\delta.
\end{align*}
The upper estimate in Laplace's principle easily follows from the previous inequalities.

The proof of Theorems \ref{T:beg1} and \ref{T:LP} is thus completed.
\section{Proof of Theorem \ref{T:nmn}}\label{S:prov}
We will use the extended contraction principle in the proof of Theorem \ref{T:nmn} 
(see Theorem 4.2.23 in \cite{DZ}, or Lemma 2.1.4 in \cite{DS}, 
see also \cite{Gar}). The main idea is to consider a family of discrete approximations to the functional $\Phi$ given by (\ref{E:disc}). Let $n\ge 2$, and define the functional
$\Phi_n:\mathbb{C}_0^m\times\mathbb{C}_0^m\times{\cal W}^d\mapsto\mathbb{C}_0^m$ as follows.
For $(l,f,h)\in\mathbb{C}_0^m\times\mathbb{C}_0^m\times{\cal W}^d$ and $\frac{jT}{n}<t\le\frac{(j+1)T}{n}$, with $1\le j\le n-1$, set
\begin{align}
\Phi_n(l,f,h)(t)&=\frac{T}{n}\sum_{k=0}^{j-1}b\left(\frac{kT}{n},h\left(\frac{kT}{n}\right)\right)+\left(t-\frac{jT}{n}\right)
b\left(\frac{jT}{n},h\left(\frac{jT}{n}\right)\right) \nonumber \\
&\quad+\sum_{k=0}^{j-1}
\sigma\left(\frac{kT}{n},h\left(\frac{kT}{n}\right)\right)\bar{C}\left[l\left(\frac{(k+1)T}{n}\right)-l\left(\frac{kT}{n}\right)\right] 
\nonumber \\
&\quad+\sigma\left(\frac{jT}{n},h\left(\frac{jT}{n}\right)\right)\bar{C}\left[l\left(t\right)-l\left(\frac{jT}{n}\right)\right] 
\nonumber \\
&\quad+\sum_{k=0}^{j-1}
\sigma\left(\frac{kT}{n},h\left(\frac{kT}{n}\right)\right)C\left[f\left(\frac{(k+1)T}{n}\right)-f\left(\frac{kT}{n}\right)\right] 
\nonumber \\
&\quad+\sigma\left(\frac{jT}{n},h\left(\frac{jT}{n}\right)\right)C\left[f\left(t\right)-f\left(\frac{jT}{n}\right)\right],
\label{E:lopi}
\end{align}
and for $0\le t\le\frac{T}{n}$, put
$
\Phi_n(l,f,h)(t)=tb(0,h(0))+\sigma(0,h(0))\bar{C}l(t)+\sigma(0,h(0))C f(t).
$
It is not hard to see that for every $n\ge 2$, the map $\Phi_n$ is continuous. 

The rest of the proof consists of three parts. First, we will show that the term
\begin{equation}
-\frac{1}{2}\varepsilon\int_0^t\mbox{diag}(\sigma(s,\widehat{B}^{(\varepsilon)}_s)\sigma(s,\widehat{B}^{(\varepsilon)}_s)^{\prime})ds
\label{E:55}
\end{equation}
appearing in (\ref{E:iu}) can be removed, and the resulting process $\widehat{X}^{\varepsilon}$ given by (\ref{E:gyt}) satisfies the same large deviation principle
as the process $X^{\varepsilon}$. It suffices to establish that the processes  $X^{\varepsilon}-x_0$ and 
$\widehat{X}^{\varepsilon}-x_0$ are exponentially equivalent, that is, for every $\delta> 0$,
\begin{align}
\limsup_{\varepsilon\rightarrow 0}\varepsilon\log
\mathbb{P}\left(\varepsilon\sup_{t\in[0,T]}||\int_0^t\mbox{diag}\,(\sigma(s,\widehat{B}^{\varepsilon}_s)
\sigma(s,\widehat{B}_s)^{\prime})ds||_d\ge\delta\right)=-\infty.
\label{E:L}
\end{align}

The second step in the proof of Theorem \ref{T:nmn} is as follows. Recall that the function $\widetilde{I}_y$ is defined 
on the space $\mathbb{C}_0^m\times\mathbb{C}_0^m\times{\cal W}^d$ as follows:
$$
\widetilde{I}_y(l,f,h)=\frac{1}{2}\int_0^T||\dot{l}(t)||_m^2dt
+\frac{1}{2}\int_0^T||\dot{f}(t)||_m^2dt
$$
when $l,f\in (H_0^1)^m$ and $h=\widehat{f}$, and by 
$\widetilde{I}_y(l,f,h)=\infty$ otherwise 
(see (\ref{E:LDPa})). It will be established that the sequence of functionals $\Phi_n$ given by (\ref{E:lopi}) approximates the 
functional $\Phi$ in (\ref{E:disc}) in the following sense: For every $\alpha> 0$,
\begin{align}
&\lim_{n\rightarrow\infty}\sup_{\{(l,f,h):\widetilde{I}_y(l,f,h)\le\alpha\}}||\Phi(l,f,h)-\Phi_n(l,f,h)||_{{\cal W}^m}=0.
\label{E:alii1}
\end{align}

Finally, we will prove that $\Phi_n(\sqrt{\varepsilon}W,\sqrt{\varepsilon}
B,\widehat{B}^{\varepsilon})$ are exponentially good approximations of $\widehat{X}^{\varepsilon}$. The latter statement means the 
following: For all $\delta> 0$,
\begin{align}
&\lim_{n\rightarrow\infty}\limsup_{\varepsilon\rightarrow 0}\varepsilon\log\mathbb{P}(\sup_{t\in[0,T]}||\widehat{X}^{\varepsilon}_t-
\Phi_n(\sqrt{\varepsilon}W,\sqrt{\varepsilon}B,\widehat{B}^{\varepsilon})(t)||_m\ge\delta)=-\infty.
\label{E:jus}
\end{align}
It is not hard to see that Theorem \ref{T:nmn} can be derived from (\ref{E:L}), (\ref{E:alii1}), (\ref{E:jus}), Theorem \ref{T:beg11}, and the extended contraction principle.
\\
\\
\it Proof of (\ref{E:L}). \rm The proof is similar to that used in Section 5 of \cite{GVol}. It follows from the continuity of 
the functions $\sigma_{ij}$ that there exists a positive even function $\eta$ defined on $[0,\infty)$ and satisfying the following conditions: $\eta$ is strictly increasing and continuous on $[0,\infty)$; $\eta(u)\rightarrow\infty$ as 
$u\rightarrow\infty$; $\sigma_{ij}(t,z)^2\le\eta(||z||_d)$ for all $1\le i,j\le m$, $t\in[0,T]$, and $z\in\mathbb{R}^d$. The inverse function of the function $\eta$ will be denoted by $\eta^{-1}$. The latter function is defined on $[\eta(0),\infty)$.

We have
\begin{align}
&\sup_{t\in[0,T]}||\int_0^t\mbox{diag}\,(\sigma(s,\widehat{B}^{\varepsilon}_s)\sigma(s,\widehat{B}_s)^{\prime})ds||_d
\le\int_0^T\sum_{i,j=1}^m\sigma_{ij}(s,\widehat{B}_s^{\varepsilon})^2ds \nonumber \\
&\le Tm^2\eta(\max_{t\in[0,T]}||\widehat{B}_t^{\varepsilon}||_d).
\label{E:fgh}
\end{align} 
It follows from (\ref{E:fgh}) that for $\varepsilon<\varepsilon_0$,
\begin{align}
&\mathbb{P}\left(\varepsilon\sup_{t\in[0,T]}
||\int_0^t\mbox{diag}\,(\sigma(s,\widehat{B}^{\varepsilon}_s)\sigma(s,\widehat{B}_s)^{\prime})ds||_d\ge\delta\right)
\le\mathbb{P}\left(\max_{t\in[0,T]}||\widehat{B}_t^{\varepsilon}||_d\ge\eta^{-1}\left(\frac{\delta}{Tm^2\varepsilon}\right)\right).
\label{E:121}
\end{align}

Let us denote by $U$ the expression on the left-hand side of (\ref{E:L}). Using (\ref{E:121}), we get
\begin{equation}
U\le\limsup_{\varepsilon\rightarrow 0}\varepsilon\log\mathbb{P}\left(\max_{t\in[0,T]}||\widehat{B}_t^{\varepsilon}||_d
\ge\eta^{-1}\left(\frac{\delta}{Tm^2\varepsilon}\right)\right).
\label{E:str}
\end{equation}
Fix $N> 0$. Then, for $\varepsilon<\varepsilon_N$, we have $\eta^{-1}\left(\frac{\delta}{Tm^2\varepsilon}\right)\ge N$. Therefore,
(\ref{E:str}) implies that
$$
U\le\limsup_{\varepsilon\rightarrow 0}\varepsilon\log\mathbb{P}\left(\max_{t\in[0,T]}||\widehat{B}_t^{\varepsilon}||_d\ge N\right).
\label{E:uoi}
$$
Set $A_N=\{\varphi\in{\cal W}^d:\max_{t\in[0,T]}||\varphi(t)||_d\ge N\}$. The set $A_N$ is a closed subset of ${\cal W}^d$. Next, using Corollary \ref{C:67} we obtain
\begin{equation}
U\le-\inf_{\{\varphi\in A_N\}}J_y(\varphi).
\label{E:56}
\end{equation}

It remains to prove that 
\begin{equation}
K_N=\inf_{\{\varphi\in A_N\}}J_y(\varphi)\rightarrow\infty\quad\mbox{as}\quad N\rightarrow \infty.
\label{E:plo}
\end{equation}
The sequence $K_N$, $N\ge 1$, is positive and nondecreasing. In order to prove (\ref{E:plo}), it suffices to show 
that the sequence $K_N$ is unbounded. We will reason by contradiction. Suppose $K_N\le C_1$ for some $C_1> 0$ and all $N\ge 1$. 
Then, there exist sequences $\varphi_N\in A_N$ and $f_N\in L^2([0,T],\mathbb{R}^m)$ such that $\varphi_N={\cal A}f_N$
and 
$
\sup_N||f_N||_{L^2([0,T],\mathbb{R}^m)}< C_2
$
for some $C_2> 0$. It follows that there exists a weakly convergent subsequence $f_{N_k}\in L^2([0,T],\mathbb{R}^m)$. Next, using
Assumption (C5), the continuity of the map $G$ on the space ${\cal W}^d$, and the equality $\varphi_N={\cal A}f_N$, 
we see that the sequence $\varphi_{N_k}$ converges in the space ${\cal W}^d$. Hence, $\max_k||\varphi_{N_k}||_{{\cal W}^d}<\infty$.
The previous inequality contradicts the assumption that $\varphi_N\in A_N$ for all $N\ge 1$. Therefore, (\ref{E:plo}) holds true. Finally,
(\ref{E:56}) implies that $U=-\infty$.

This completes the proof of (\ref{E:L}).\\
\\
\it Proof of (\ref{E:ali1}). \rm It is not hard to see that
\begin{equation}
\Phi_n(l,f,\widehat{f})(t)=\int_0^tb_n(s,\widehat{f}(s))ds+\int_0^t\sigma_n(s,\widehat{f}(s))\bar{C}\dot{l}(s)ds
+\int_0^t\sigma_n(s,\widehat{f}(s))C\dot{f}(s)ds
\label{E:foot}
\end{equation}
where the maps $b_n$ and $\sigma_n$ are given by
$$
b_n(s,\widehat{f}(s))=\sum_{k=0}^{n-1}b\left(\frac{kT}{n},\widehat{f}\left(\frac{kT}{n}\right)\right)
\mathbb{1}_{\frac{kT}{n}< s\le\frac{(k+1)T}{n}}
$$
and 
$$
\sigma_n(s,\widehat{f}(s))=\sum_{k=0}^{n-1}\sigma\left(\frac{kT}{n},\widehat{f}\left(\frac{kT}{n}\right)\right)
\mathbb{1}_{\frac{kT}{n}< s\le\frac{(k+1)T}{n}}.
$$
Recall that 
\begin{align} 
\Phi(l,f,\widehat{f})(t)&=\int_0^tb(s,\widehat{f}(s))ds+\int_0^t\sigma(s,\widehat{f}(s))\bar{C}\dot{l}(s)ds
\nonumber \\
&\quad+\int_0^t\sigma(s,\widehat{f}(s))C\dot{f}(s)ds
\label{E:discc} 
\end{align}
(see (\ref{E:disc})). For every $\alpha> 0$, set 
$$
D_{\alpha}=\{\tau\in(\mathbb{H}_0^1)^m:\int_0^T||\dot{\tau}(t)||_m^2dt\le 2\alpha\}.
$$ 
Note that $l$ and $f$ appearing 
in (\ref{E:ali1}) belong to $D_{\alpha}$. Next, we take into account that $f\mapsto\widehat{f}$ is a compact map from any closed ball 
in $(\mathbb{H}_0^1)^m$ into ${\cal W}^d$. The previous statement follows from Assumption (C5) and the definition of the map 
$f\mapsto\widehat{f}$. Applying the Arzel\`{a}-Ascoli theorem, we obtain the following formulas:
$$
\sup_{h\in D_{\alpha}}\sup_{\{t\in[0,T]\}}||\widehat{f}(t)||_d< M_{\alpha}
$$
and
\begin{equation}
q_{\alpha,n}=\sup_{f\in D_{\alpha}}\sup_{\{t,u\in[0,T]:|t-u|\le\frac{T}{n}\}}||\widehat{f}(t)-\widehat{f}(u)||_d\rightarrow 0
\label{E:AA}
\end{equation}
as $n\rightarrow\infty$. Finally, using Assumption A, (\ref{E:foot}), (\ref{E:discc}), and reasoning as in the proofs of Lemmas 6.23 and 6.24 
in \cite{Gul1}, we obtain
\begin{align}
&\sup_{\{(l,f,\widehat{f}):\widetilde{I}_y(l,f,\widehat{f})\le\alpha\}}||\Phi(l,f,\widehat{f})(t)-\Phi_n(l,f,\widehat{f)}(t)||_{{\cal W}^m} \nonumber \\
&\le\sup_{\{l,f\in D_{\alpha}\}}||\Phi(l,f,\widehat{f})(t)-\Phi_n(l,f,\widehat{f})(t)||_{{\cal W}^m} 
\le\zeta\omega\left(\frac{T}{n}+q_{\alpha,n}\right)
\label{E:dli}
\end{align}
where the constant $\zeta> 0$ does not depend on $n$, and $\omega$ is the modulus of continuity in Assumption A. Now, (\ref{E:ali1}) follows from
(\ref{E:AA}) and (\ref{E:dli}). \\
\\
\it Proof of (\ref{E:jus}). \rm The structure of the proof of the formula in (\ref{E:jus}) is the same as in a similar proof of 
Lemma 6.25 in \cite{Gul1}. Note that in \cite{Gul1}, the case where $d=1$ is considered, while in the present paper, we deal with multivariate models. Moreover, the volatility process in \cite{Gul1} is Gaussian. It will be explained below how to take into account these differences. We will need the following two lemmas.
\begin{lemma}\label{L:lobs}
Let $0< r< r_0$ where $r_0$ is a small number. Suppose $q$ is the function on $(0,r_0]$ defined in the proof of Lemma 6.26 in \cite{Gul1}.
Then, the following equality holds:
\begin{align}
&\limsup_{r\rightarrow 0}\limsup_{\varepsilon\rightarrow 0}\varepsilon\log\mathbb{P}(\sup_{t\in[0,T]}||\widehat{B}_t^{\varepsilon}||_d
\ge 2^{-1}q(r))=-\infty.
\label{E:lll}
\end{align}
\end{lemma}
\begin{lemma}\label{L:otter}
For every $\beta> 0$, the following formula is valid:
\begin{align}
&\limsup_{n\rightarrow\infty}\limsup_{\varepsilon\rightarrow 0}\varepsilon\log\mathbb{P}(\max_{t,s\in[0,T]:|t-s|\le\frac{T}{n}}
||\widehat{B}_t^{\varepsilon}-\widehat{B}_s^{\varepsilon}||_d\ge\beta)=-\infty.
\label{E:ott}
\end{align}
\end{lemma}
\begin{remark}\label{R:rrrs}
The function $q$ appearing in the formulation of Lemma \ref{L:lobs} is positive, strictly decreasing, and continuous 
(see the proof of Lemma 6.26 in \cite{Gul1}). Moreover, $\lim_{r\rightarrow 0}q(r)=\infty$. 
For the statements similar to those in Lemmas \ref{L:lobs} and \ref{L:otter}, see (6.50) and (6.51) in \cite{Gul1}.
\end{remark}

\it Proof of Lemma \ref{L:lobs}. \rm The proof of (\ref{E:lll}) is similar to that of (\ref{E:L}). For every $r\in(0,r_0]$, define a closed subset of ${\cal W}^d$ by
$
E_r=\{\varphi\in{\cal W}^d:\sup_{t\in[0,T]}||\varphi(t)||_d\ge 2^{-1}q(r)\}
$
and set 
$
N_r=\limsup_{\varepsilon\rightarrow 0}\varepsilon\log\mathbb{P}(\sup_{t\in[0,T]}||\widehat{B}_t^{\varepsilon}||_d
\ge 2^{-1}q(r)).
$ 
The large deviation principle in Corollary \ref{C:67} implies that $N_r=-\inf_{\varphi\in E_r}J_y(\varphi)$
where $J_y$ is the rate function defined in (\ref{E:LDPs}). Set $K_r=-N_r$. Then, $K_r$ is a nonnegative nonincreasing function on $(0,r_0]$.
It remains to prove that $\lim_{r\rightarrow 0}K_r=\infty$. We will reason by contradiction. Suppose $K_r\le M$ for all $r\in(0,r_0]$.
Then, exactly as in the proof of the formula in (\ref{E:L}), we see that there exist a sequence $r_k\rightarrow 0$ 
and a sequence $\varphi_k\in N_{r_k}$ such that the set $\{\varphi_k\}$ is compact in ${\cal W}^d$. It follows that the set 
$\{\varphi_k\}$ is bounded in ${\cal W}^d$. This means that we arrived at a contradiction since $||\varphi_k||_{{\cal W}^d}
\ge 2^{-1}q(r_k)$ and $q(r_k)\rightarrow\infty$ as $k\rightarrow\infty$. It follows that $\lim_{r\rightarrow 0}N_r=-\infty$.

This completes the proof of Lemma \ref{L:lobs}.

\it Proof of Lemma \ref{L:otter}. \rm For every $n\ge 1$, define a closed subset of ${\cal W}^d$ by
$$
V_n=\{g\in{\cal W}^d:\max_{t,s\in[0,T]:|t-s|\le\frac{T}{n}}||g(t)-g(s)||_d\ge\beta\}
$$
and set
$$
U_n=\limsup_{\varepsilon\rightarrow 0}\varepsilon\log\mathbb{P}(\max_{t,s\in[0,T]:|t-s|\le\frac{T}{n}}
||\widehat{B}_t^{\varepsilon}-\widehat{B}_s^{\varepsilon}||_d\ge\beta).
$$
Using the large deviation principle established in Corollary \ref{C:67}, we can prove the following equality:
$
U_n=-\inf_{\varphi\in V_n}J_y(\varphi)
$
where $J_y$ is the rate function defined in (\ref{E:LDPs}). Set $K_n=-U_n$. Then, $K_n$ is a nonnegative nondecreasing sequence.

Our next goal is to prove that $K_n\rightarrow\infty$ as $n\rightarrow\infty$. 
It suffices to show that the sequence $K_n$ is unbounded. To prove the previous statement, we reason by contradiction like 
in the proof of the formula in (\ref{E:L}). Suppose $K_n\le C$ for all $n\ge 1$. Then, there exists a sequence
$\varphi_{n_k}\in V_{n_k}$ that is compact in ${\cal W}^d$ (see the proof of (\ref{E:L})). Applying the Arzel\`{a}-Ascoli theorem,
we obtain
$$
\sup_{k\ge 1}\max_{t,s\in[0,T]:|t-s|\le\frac{T}{m}}||\varphi_{n_k}(t)-\varphi_{n_k}(s)||_d\rightarrow 0
$$
as $m\rightarrow\infty$. Therefore, there exists $m_0$ such that 
\begin{equation}
\max_{t,s\in[0,T]:|t-s|\le\frac{T}{m}}||\varphi_{n_k}(t)-\varphi_{n_k}(s)||_d<\frac{\beta}{2}
\label{E:er}
\end{equation}
for all $k\ge 1$ and $m\ge m_0$. Let $k_0$ be such that $n_{k}> m_0$ for all $k\ge k_0$. Then, (\ref{E:er}) implies that
\begin{equation}
\max_{t,s\in[0,T]:|t-s|\le\frac{T}{n_k}}||\varphi_{n_k}(t)-\varphi_{n_k}(s)||_d<\frac{\beta}{2}
\label{E:era}
\end{equation}
for all $k\ge k_0$. It is not hard to see that (\ref{E:era}) contradicts the condition $\varphi_{n_k}\in V_{n_k}$ for all $k\ge 1$.

This completes the proof of Lemma \ref{L:otter}.

Let us return to the proof of (\ref{E:jus}).
Using (\ref{E:iu}) and (\ref{E:foot}), we obtain
\begin{align}
\widehat{X}^{\varepsilon}_t-
\Phi_n(\sqrt{\varepsilon}W,\sqrt{\varepsilon}B,\widehat{B}^{\varepsilon})(t)
&=\int_0^t[b(s,\widehat{B}_s^{\varepsilon})-b_n(s,\widehat{B}_s^{\varepsilon})]ds
+\sqrt{\varepsilon}\int_0^t[\sigma(s,\widehat{B}_s^{\varepsilon})-\sigma_n(s,\widehat{B}_s^{\varepsilon})]\bar{C}dW_s \nonumber \\
&\quad+\sqrt{\varepsilon}\int_0^t[\sigma(s,\widehat{B}_s^{\varepsilon})-\sigma_n(s,\widehat{B}_s^{\varepsilon})]CdB_s,\quad t\in[0,T].
\label{E:pf}
\end{align}
By analyzing the formula in (\ref{E:pf}), we observe that in order to prove (\ref{E:jus}), it suffices to show that for all $1\le i,j,k\le m$ and $\kappa> 0$,
\begin{align}
&\lim_{n\rightarrow\infty}\limsup_{\varepsilon\rightarrow 0}\varepsilon\log\mathbb{P}(\sup_{t\in[0,T]}
|\int_0^t[b_i(s,\widehat{B}_s^{\varepsilon})-b_{i}^{(n)}(s,\widehat{B}_s^{\varepsilon})]ds|\ge\kappa)=-\infty,
\label{E:ht1}
\end{align}
\begin{align}
&\lim_{n\rightarrow\infty}\limsup_{\varepsilon\rightarrow 0}\varepsilon\log\mathbb{P}(\sqrt{\varepsilon}\sup_{t\in[0,T]}
|\int_0^t[\sigma_{ij}(s,\widehat{B}_s^{\varepsilon})-\sigma_{ij}^{(n)}(s,\widehat{B}_s^{\varepsilon})]dW_k(s)|\ge\kappa)=-\infty,
\label{E:ht2}
\end{align}
and
\begin{align}
&\lim_{n\rightarrow\infty}\limsup_{\varepsilon\rightarrow 0}\varepsilon\log\mathbb{P}(\sqrt{\varepsilon}\sup_{t\in[0,T]}
|\int_0^t[\sigma_{ij}(s,\widehat{B}_s^{\varepsilon})-\sigma_{ij}^{(n)}(s,\widehat{B}_s^{\varepsilon})]dB_k(s)|\ge\kappa)=-\infty.
\label{E:ht3}
\end{align}

The proofs of the equalities in (\ref{E:ht1}) -- (\ref{E:ht3}) are similar to those of (6.38) -- (6.40) in \cite{Gul1}. These proofs are long and involved, but all the necessary details are given in \cite{Gul1}. However, there are also certain differences in the proofs in \cite{Gul1} and in the present paper. First of all, the process $\varepsilon\mapsto\sqrt{\varepsilon}\widehat{B}$ in \cite{Gul1} should be replaced by the process
$\varepsilon\mapsto\widehat{B}^{\varepsilon}$. We also use Lemmas \ref{L:lobs} and \ref{L:otter} instead of (6.50) and Corollary 6.22 
in \cite{Gul1}. The stopping time $\xi^{\varepsilon,m,r}$ defined in (6.42) of \cite{Gul1} is replaced by the following stopping time:
For $\varepsilon\in(0,1]$, $n\ge 2$, and $0< r< r_0$, we set
$$
\xi^{\varepsilon,m,r}=\inf_{s\in[0,T]}\left\{\frac{r}{q(r)}||\widehat{B}_s^{\varepsilon}||_d+||\widehat{B}_s^{\varepsilon}
-\widehat{B}_{\frac{[nsT^{-1}]T}{n}}^{\varepsilon}||_d
> r\right\}.
$$
Note that everywhere in the proof of (\ref{E:ott}) we assume that the filtration $\{{\cal F}_t^B\}$ is hidden in the background.
By taking into account the remarks mentioned above, we can complete the proof of (\ref{E:jus}) by reasoning as in the proof of 
Lemma 6.25 in \cite{Gul1}. We do not provide more details here, and leave finding a detailed proof as an exercise for the interested reader.

Finally, we can finish the proof of Theorem \ref{T:nmn} by using Theorem \ref{T:beg1}, the formulas in (\ref{E:L}), (\ref{E:alii1}), and 
(\ref{E:jus}), and applying the extended contraction principle.

\section{Acknowledgements}\label{S:ac}
I am indebted to Peter Friz and Stefan Gerhold for valuable remarks. I also thank the anonymous referee for reading the paper
and for providing useful comments which significantly contributed to improving the paper.

\end{document}